\documentclass[12pt,leqno]{article}
\let\stemptyset\emptyset
\let\stcup\cup
\let\stcap\cap
\let\stnot\not
\let\sttimes\times
\let\stin\in
\let\stneg\neg
\let\stwedge\wedge
\let\stvee\vee

\usepackage{latexsym}
\usepackage{amstext}
\usepackage{amssymb}
\usepackage{amsfonts}
\usepackage{amsthm}
\usepackage{amsopn}
\usepackage{amsbsy}
\usepackage{layout}
\usepackage[usenames,dvipsnames,svgnames,table]{xcolor}
\definecolor{mygreen}{RGB}{0, 160, 0}
\usepackage{amsmath}
\usepackage{graphicx}
\usepackage{tikz}
\usetikzlibrary{positioning,shapes,shadows}
\usetikzlibrary{arrows,automata}
\usepackage{bm}
\usepackage{yfonts}
\usepackage{mathtools}
\usepackage{mathabx} 
\usepackage{esvect}

\def\newtheorems{\newtheorem{theorem}{Theorem}[section]
                 \newtheorem{cor}[theorem]{Corollary}
                 \newtheorem{prop}[theorem]{Proposition}
                 \newtheorem{lemma}[theorem]{Lemma}
                 \newtheorem{definition}[theorem]{Definition}
                 \newtheorem{notation}[theorem]{Notations}

                 \newtheorem{remark}[theorem]{Remark}
                 
                 \newtheorem{observation}[theorem]{Observation}
                 \newtheorem{question}[theorem]{Question}

}

\newtheorems

\font\nmini=cmr10 scaled 550

\font\tinier=cmsy6 scaled 730
\font\sbf=cmbx8 scaled 850
\font\sbfreg=cmbx8 scaled 1200
\font\srbf=cmbx8 scaled 920
\font\ssbf=cmbx8 scaled 700
\font\scmu=cmu10 scaled 600

\font\lcmssifont=lcmssi8 scaled 1000 

\font\blcmssifont=lcmssi8 scaled 1350 

\def\fs#1{\mbox{\it #1\kern 1.3pt}}
\def\rfs#1{\mbox{\rm #1\kern 1.3pt}}
\def\nrfs#1{\mbox{\rm #1\kern 0.0pt}}
\def\bffs#1{\mbox{\bf #1\kern 1.3pt}}
\def\bfs#1{\mbox{\bf #1\kern 1.3pt}}
\def\fss#1{\mbox{\scriptsize\it #1\kern 1.3pt}}
\def\fst#1{\mbox{\tiny\it #1\kern 1.1pt}}
\def\sifs#1{\mbox{\scriptsize\it #1\kern 1.3pt}}
\def\srfs#1{\mbox{\kern0.7pt\scriptsize\rm #1\kern 1.3pt}}
\def\oversrfs#1{\mbox{\kern2pt$\overline{\kern-2pt\mbox{\kern0.7pt\scriptsize\rm #1\kern 1.3pt}\kern-2pt}$\kern2pt}}
\def\bfsrfs#1{\mbox{\kern0.7pt\scriptsize\bf #1\kern 1.3pt}}
\def\tbfs#1{\mbox{\kern0.7pt\ssbf #1\kern -0.6pt}}
\def\sbfs#1{\mbox{\kern0.7pt\srbf #1\kern -0.6pt}}
\def\srbfs#1{\mbox{\kern0.7pt\srbf #1\kern -0.6pt}}
\def\spfs#1{\mbox{\kern0.7pt\scmu #1\kern 1.3pt}}
\def\sspfs#1{\mbox{\kern0.5pt\sscmu #1\kern 1.1pt}}
\def\ssbfs#1{\mbox{\kern0.7pt\ssbf #1\kern 1.3pt}}
\def\fsm#1{\mbox{\tiny\it #1\kern 1.0pt}}
\def\trfs#1{\mbox{\tiny\rm #1\kern 1.0pt}}
\def\utrfs#1{\raise1.8pt\hbox{\mbox{\kern1.4pt\tiny\rm #1\kern 1.0pt}}}
\def\tbffs#1{\mbox{\tiny\bf #1\kern 1.0pt}}
\def\ttrfs#1{\raise2pt\hbox{\kern0.7pt\mbox{\tiny\rm #1\kern 1.0pt}}}
\def\rtrfs#1{\raise 1.7pt\hbox{\tiny\rm\kern0.5pt #1\kern 1.0pt}}
\def\tfs#1{\mbox{\tiny\it #1\kern 1.0pt}}

\def\({\left(}
\def\){\right)}
\def\what0{\widehat{0}}
\def\w={\kern3pt=\kern3pt}

\newcommand{\comp}{\hbox{$<\kern -3pt >$}}
\newcommand{\ncomp}
		{\;\hbox{\hbox{/}\kern -9.5pt \hbox{$<\kern -3pt >$}}}

\newcommand{\meet}
	       {\hbox{$\wedge \kern -5.75pt \raise 1.5pt \hbox{$.$}\,$}}
\newcommand{\Meet}
	     {\hbox{$\bigwedge \kern -8pt \raise 0.75pt \hbox{$.$}\:$}}
\newcommand{\ld}
	       {\hbox{$< \kern -6pt \raise 2pt \hbox{$.$}\,$}}
\newcommand{\sss}{\: \hbox{$
\underline{\hbox{$\subset$}}\kern -4pt\raise -2pt \hbox{$\tiny |$}
$}\: }
\newcommand{\almostcontained}{\: \hbox{$
\raise 1.5pt \hbox{\scriptsize $\subset$}\kern -6.3pt\raise -3.5pt
\hbox{\scriptsize $\sim$}
$}\: }
\newcommand{\rraro}[2]{\hbox{$\kern 3pt\raise 2pt \hbox{$\raro$}
 \kern -14pt \raise
-3.5pt\hbox{\tiny{$#1\raro #2$}}$}}

\newcommand{\frc}{\hbox{$\parallel \kern -5.7pt \hbox{$-$}$}}

\newcommand{\nrestriction}{\kern-2.5pt\upharpoonright\kern-2.5pt}
\newcommand{\nnrestriction}{\kern-4.5pt\upharpoonright\kern-4.5pt}
\newcommand{\wrestriction}{\kern3pt\restriction\kern3pt}
\newcommand{\wequality}{\kern3pt=\kern3pt}

\newcommand{\nfrc}{\not \kern -5pt \frc}
\newcommand{\rest}{\vbox{\hbox{$\:\kern -2pt\mathbin{\vert\kern-3.1pt\lower-1pt
   \hbox{$\mathsurround=0pt\mathchar"0012$}\kern-4pt}\:$}}}

\newcommand{\drest}{\vbox{\hbox{$\:\kern -2pt\mathbin{\rest\kern-4.6pt\lower1.7pt
   \hbox{$\mathsurround=0pt\mathchar"0012$}\kern-4pt}\:$}}}
\newcommand{\vdrest}{\vbox{\hbox{$\:\kern -2pt\mathbin{\rest\kern-3.95pt\lower1.7pt
   \hbox{$\mathsurround=0pt\mathchar"0012$}\kern-4pt}\:$}}}
   \newcommand{\smdrest}{\vbox{\hbox{$\:\kern -2pt\mathbin{\smrest\kern-2.7pt\lower1.7pt
      \hbox{$\mathsurround=0pt\mathchar"0012$}\kern-4pt}\:$\kern 1pt}}}

\newcommand{\rests}{\vbox{\hbox{\scriptsize$\:\kern
-1.4pt\mathbin{\vert\kern-2.4pt\lower-1pt
\hbox{\scriptsize$\mathsurround=0pt\mathchar"0012$}\kern-
3.0pt}\:$}}}

\newcommand{\cntd}{\subseteq}

\newcommand{\pcntda}{\lower5pt\hbox{$\stackrel{\subset}{\neq}$}}
\newcommand{\pcntdb}{\lower5pt\hbox{$\stackrel{\supset}{\neq}$}}
\newcommand{\pcntdc}{\lower5pt\hbox{$\stackrel{\subseteq}{\neq}$}}

\newcommand{\m}[1]{\hbox{$ #1 $}}

\newcommand{\itPhi}{{\it \Phi}}
\newcommand{\itPsi}{{\it \Psi}}

\newcommand{\bfitGamma}{{\it \Gamma \kern-9.5pt \Gamma}}

\newcommand{\bfitDelta}{{\it \Delta \kern-9.2pt \Delta}}

\newcommand{\sKappa}{\lower1.5pt\hbox{\small$\kappa$}}

\newcommand{\longvec}[1]{\stackrel{\hbox{\scriptsize$\rightarrow$}}{#1}\kern-2.5pt}
\newcommand{\llongvec}[1]{\stackrel{\hbox{\scriptsize$\longrightarrow$}}{#1}\kern-2.5pt}

\newcommand{\gammaprime}{\gamma\kern1pt'}

\newcommand{\bfB}{\hbox{\bf B}}

\newcommand{\superbfL}{\tiny\bf \kern1ptL}

\newcommand{\superbfM}{\tiny\bf \kern1ptM}

\newcommand{\bfO}{{\bf O}}

\newcommand{\superbfS}{\tiny\bf \kern1ptS}

\newcommand{\bfy}{\hbox{\bf y}}

\newcommand{\bfX}{{\bf X}}
\newcommand{\bfY}{{\bf Y}}
\newcommand{\bfZ}{\hbox{\bf Z}}
\newcommand{\bfIL}{\hbox{\bf I\kern-2ptL}}
\newcommand{\bfJI}{\hbox{\bf J\kern-2ptI}}

\newcommand{\calL}{{\cal L}}

\newcommand{\swhatT}
{
\mbox{\tiny$\widehat{\rule{0pt}{9.0pt}\kern-1.0pt\mbox{\normalsize$T$}\kern-1.0pt}$}
}

\newcommand{\dwhatT}{\mbox{\swhatT\kern-8.9pt$\widehat{\kern1.5pt\rule{0pt}{8.5pt} \kern1.5pt}$}\rule{4pt}{0pt}}

\newcommand{\hatI}{\hat{I}}

\newcommand{\vecrmD}{\hbox{$\stackrel{\rightarrow}{\rule{2pt}{0pt}\rfs{D}}\kern-0.5pt$}}

\newcommand{\swhatX}{\kern6pt\widehat{\rule{0pt}{8pt}}\kern-5.9pt X}

\newcommand{\calJet}{{\cal J}{\kern-2pt\it et}}

\newcommand{\hatoverU}{\hat{\overU\kern2pt}\kern-2pt}

\newcommand{\veca}{\vec{a}}
\newcommand{\vecb}{\vec{b}}

\newcommand{\vecf}{\vec{f}}

\newcommand{\svecu}{\vec{u}\kern1.5pt}

\newcommand{\svecv}{\vec{v}\kern1.5pt}

\newcommand{\svecw}{\vec{w}\kern1.5pt}
\newcommand{\vecx}{\vec{x}}
\newcommand{\svecx}{\vec{x}\kern1.5pt}

\newcommand{\svecy}{\vec{y}\kern1.5pt}

\newcommand{\svecz}{\vec{z}\kern1.5pt}

\newcommand{\svecA}{\vec{A}\kern1.5pt}

\newcommand{\vecphi}{\vec{\phi}}

\newcommand{\swtildeX}{\kern6pt\widetilde{\rule{0pt}{8pt}}\kern-5.9pt X}

\newcommand{\raro}{\rightarrow}

\newcommand{\itmb}[1]{\item[[\kern 1.5pt #1\kern -4pt]]}
\newcommand{\itms}[1]{\item[[#1\kern -5pt]]}

\newcommand{\II}{{\bf I\kern -1pt I}}

\newcommand{\eqdf}{\hbox{\bf \,:=\,}}
\newcommand{\wideeqdf}{\hbox{\kern5pt\bf \,:=\,\kern5pt}}

\newcommand{\surjectionarrow}{\stackrel{\mbox{\tiny onto}}{\longrightarrow}}

\newcommand{\surj}[3]{#1 \kern-1.5pt:\kern-1.5pt #2 \surjectionarrow #3}

\newcommand{\bijectionarrow}{\vbox{\hbox{$\longrightarrow $
                  \kern -22pt \hbox{\lower 2.5pt  \hbox{\tiny onto}}
                  \kern -16pt \hbox{\raise 5pt  \hbox{\tiny 1-1}}
                  \kern 3pt}}}
\newcommand{\uarrow}[2]{\vbox{\hbox{$\longrightarrow $
                  \kern -16pt \hbox{\raise 5pt  \hbox{\tiny $#1$}}
                  \kern 10pt }}}

\newcommand{\spri}{\vbox{\hbox{\raise 2pt \hbox{\tiny $\|$}}}}
\newcommand{\tlr}{\vbox{\hbox{\raise 2pt \hbox{\tiny $\leftarrow$}}}}
\newcommand{\trr}{\vbox{\hbox{\raise 2pt \hbox{\tiny $\rightarrow$}}}}
\newcommand{\spr}{\mathrel{\vbox{\hbox{\tlr \kern -1.7pt \spri \kern -1.7pt \trr}}}}

\newcommand{\mspri}{\vbox{\hbox{\raise 2pt \hbox{$\scriptscriptstyle \|$}}}}

\newcommand{\mtlr}{\vbox{\hbox{\raise 2pt \hbox{$\scriptscriptstyle \leftarrow$}}}}

\newcommand{\mtrr}{\vbox{\hbox{\raise 2pt \hbox{$\scriptscriptstyle \rightarrow$}}}}

\newcommand{\mspr}{\mathrel{\vbox{\hbox{\mtlr \kern -1.7pt \mspri \kern -1.7pt \mtrr}}}}

\newcommand{\Nm}{\hbox{\kern -1.3pt \em I\kern-.2300em N\,}}
\newcommand{\Na}{\hbox{\kern -1.3pt \it I\kern-.2300em N\,}}
\newcommand{\Aa}{\hbox{\it A\kern -6.8pt\lower 1.0pt\hbox{-}\,}}
\newcommand{\Ba}{\hbox{\kern -1.3pt \it I\kern-.2300em B\,}}
\newcommand{\Da}{\hbox{\kern -1.3pt \it I\kern-.2300em D\,}}
\newcommand{\Ka}{\hbox{\kern -1.3pt \it I\kern-.2300em K\,}}
\newcommand{\La}{\hbox{\kern -1.3pt \it I\kern-.2300em L\,}}
\newcommand{\Ta}{\hbox{\kern 0.5pt \it T\kern-.5550em T\,}}
\newcommand{\jTa}{\hbox{\kern 0.5pt \it T\kern-.5550em T\,}}
\newcommand{\nTa}{\hbox{\kern 0.5pt \it T\kern-.6300em T\,}}
\newcommand{\Nmt}{\hbox{\kern -1.3pt I\kern-.2300em N\,}}
\newcommand{\N}{\hbox{I\kern-.1500em \hbox{\sf N}}}
\newcommand{\As}{\hbox{\scriptsize\it A\kern -5.3pt\lower 0.7pt\hbox{\it -}\,}}
\newcommand{\Ns}{\hbox{{\scriptsize\it I\kern-.1500em N}}}
\newcommand{\Nss}{\hbox{{\tiny\it I\kern-.1500em N}}}
\newcommand{\Ls}{\hbox{{\scriptsize\it I\kern-.1700em L}}}
\newcommand{\Rm}{\hbox{\kern -1.3pt \em I\kern-.1950em R\,}}
\newcommand{\Ra}{\hbox{\kern -1.3pt \it I\kern-.1950em R\,}}
\newcommand{\R}{\hbox{I\kern-.1500em \hbox{\sf R}}}
\newcommand{\sR}{\hbox{\tiny \hbox{I\kern-.1500em \hbox{\sf R}}}}
\newcommand{\msR}{\hbox{\tiny \hbox{\it I\kern-.2200em \hbox{\it R}}}}
\newcommand{\Rs}{\hbox{\scriptsize\it \hbox{I\kern-.2100em \hbox{R}}}}

\newcommand{\Q}
   {\hbox{${\rm Q} \kern -7.5pt \raise 2pt \hbox{\tiny$|$}\kern 7.5pt$}}
\newcommand{\C}
   {\hbox{${\rm C} \kern -6.5pt \raise 2pt \hbox{\tiny$|$}\kern 6.5pt$}}
\newcommand{\sC}
   {\hbox{\tiny \hbox{${\rm C} \kern -5.0pt \raise 1pt \hbox{\tiny$|$}\kern 7.5pt$}}}
\newcommand{\Z}{\hbox{\sf Z\kern-0.720em\hbox{ Z}}}
\newcommand{\sZ}{\hbox{\tiny\hbox{ \sf Z\kern-0.720em\hbox{ Z}}}}

\newcommand{\myqed}{\kern 5pt\vrule height7.5pt width6.9pt depth0.2pt \rule{1.3pt}{0pt}}
\newcommand{\solidqed}{\kern 5pt\vrule height7.5pt width6.9pt depth0.2pt \rule{1.3pt}{0pt}}

\newcommand{\oldcrossedqed}{$\raise0.52pt\hbox{\large $\sttimes$}\kern-10.10pt\square$}
\newcommand{\oldcrossedsquare}{$\raise0.52pt\hbox{\large $\sttimes$}\kern-10.10pt\square$}

\newcommand{\proofend}
{
\hbox{
\rule{0.8pt}{7pt}
\kern-3.9pt
\raise6.3pt
\hbox{\rule{5.0pt}{0.8pt}}
\kern-3.8pt
\rule{0.8pt}{7pt}
\kern-9.8pt
\rule{5.8pt}{0.8pt}
}
\kern-06pt}

\newcommand{\proofends}{\proofend\kern-07pt}

\newcommand{\proofendeol}
{
\hbox{
\rule{0.8pt}{7pt}
\kern-3.9pt
\raise6.3pt
\hbox{\rule{4.1pt}{0.8pt}}
\kern-3.8pt
\rule{0.8pt}{7pt}
\kern-9.0pt
\rule{4.9pt}{0.8pt}
\kern-5pt
}
\kern-8pt
}

\newcommand{\bd}{\begin{description}}
\newcommand{\ed}{\end{description}}

\newcommand{\ben}{\begin{enumerate}}
\newcommand{\een}{\end{enumerate}}

\newcommand{\sngltn}[1]{\{ #1 \}}

\newcommand{\dbltn}[2]{\{ #1, #2 \}}
\newcommand{\threetn}[3]{\{ #1, #2, #3 \}}

\newcommand{\fourtn}[4]{\{ #1, #2, #3, #4 \}}

\newcommand{\sixtn}[6]{\{ #1, #2, #3, #4, #5, #6 \}}

\newcommand{\setm}[2]{\{#1\:|\:#2\}}

\newcommand{\seqn}[2]{\langle#1,\ldots ,#2\rangle}

\newcommand{\fsetn}[2]{\{\,#1,\ldots ,#2\}}
\newcommand{\fseqn}[2]{\langle\,#1,\ldots ,#2\rangle}

\newcommand{\boldnorm}[1]
{
\kern1pt\lower3pt\hbox{\rule{1.0pt}{12pt}\kern1.4pt}
#1
\lower3pt\hbox{\kern1.6pt\rule{1.0pt}{12pt}\kern2.0pt}
                                            }

\newcommand{\bboldnorm}[1]
{
\kern1pt\lower3pt\hbox{\rule{1.2pt}{12pt}\kern1.4pt}
#1
\lower3pt\hbox{\kern1.6pt\rule{1.2pt}{12pt}\kern2.0pt}
                                            }

\newcommand{\sboldnorm}[1]
{
\kern1pt\lower1.6pt\hbox{\rule{0.9pt}{07pt}\kern1.1pt}
#1
\lower1.6pt\hbox{\kern1.1pt\rule{0.9pt}{07pt}\kern1.8pt}
                                            }

\newcommand{\sbboldnorm}[1]
{
\kern1pt\lower0.5pt\hbox{\rule{1.0pt}{06pt}\kern1.1pt}
#1
\lower0.5pt\hbox{\kern1.1pt\rule{1.0pt}{06pt}\kern1.8pt}
                                            }
\newcommand{\abs}[1]{| #1 |}

\newcommand{\onetpl}[1]{\langle#1\kern1.4pt\rangle}
\newcommand{\pair}[2]{\langle#1 ,#2\kern1.4pt\rangle}
\newcommand{\wpair}[2]{\langle\rule{1pt}{0pt}#1 ,#2\kern1.4pt\rangle}
\newcommand{\twotpl}[2]{\langle#1 ,#2\kern1.4pt\rangle}
\newcommand{\smpair}[2]{\langle#1 ;#2\kern1.4pt\rangle}

\newcommand{\rvlpar}{\kern2pt(\kern-1.3pt\lower2.225pt\hbox{\rule{0.5pt}{10.65pt}}\kern1.9pt}
\newcommand{\srvlpar}{(\kern-1.3pt\lower1.45pt\hbox{\rule{0.4pt}{6.8pt}}\kern1.7pt}
\newcommand{\rvrpar}{\kern1.9pt\lower2.20pt\hbox{\rule{0.5pt}{10.6pt}}\kern-1.2pt)\kern2pt}
\newcommand{\srvrpar}{\kern1.7pt\lower1.4pt\hbox{\rule{0.4pt}{6.8pt}}\kern-1.2pt)}
\newcommand{\rvpair}[2]{\rvlpar#1 ,#2\rvrpar}
\newcommand{\srvpair}[2]{\srvlpar#1 ,#2\srvrpar}
\newcommand{\rvtrpl}[3]{\rvlpar#1, #2, #3\kern1pt\rvrpar}
\newcommand{\drvpair}[4]{\rvlpar#1, #2 \kern1pt\vert\kern1pt #3 ,#4\kern1pt\rvrpar}

\newcommand{\trpl}[3]{\langle#1 ,#2 ,#3 \rangle}

\newcommand{\semitrpl}[3]{\langle#1 ,#2 \kern1pt;#3 \rangle}
\newcommand{\qdrpl}[4]{\langle#1 ,#2 ,#3 ,#4 \rangle}
\newcommand{\fourtpl}[4]{\langle#1 ,#2 ,#3 ,#4\kern1.4pt \rangle}

\newcommand{\fivetpl}[5]{\langle #1 ,#2 ,#3 ,#4 ,#5 \rangle}

\newcommand{\change}[1]{\lower 0.7pt \hbox{\mbox{\large $#1$}}}

\newcommand{\fnn}[3]{#1:#2 \raro #3}

\newcommand{\fnnp}[3]{#1\hbox{\kern3pt$:$\kern0.6pt\raise1.5pt\hbox{\tiny$\subseteq$ \kern2pt}}#2 \raro #3}

\newcommand{\tendin}[3]
{\lim_{#1\kern 1pt \ni\kern 1.5pt #2\kern 1pt \rightarrow\kern 1.5pt #3}}
\newcommand{\tend}[2]
{\lim_{#1\kern 1pt \rightarrow\kern 1.5pt #2}}

\newcommand{\iso}[3]{\m{ #1 : #2 \cong #3 }}
\newcommand{\isop}[3]{#1\hbox{\kern3pt$:$\kern0.6pt\raise1.5pt\hbox{\tiny$\subseteq$ \kern2pt}}#2 \cong #3}

\newcommand{\inverse}{^{-1}}

\newcommand{\tsubseteq}{\lower2.1pt\hbox{\tiny$\subseteq$}}
\newcommand{\temptyset}{\lower1.6pt\hbox{\tiny$\stemptyset$}}
\newcommand{\tSml}{\lower1.6pt\hbox{\tiny\it Sml}}
\newcommand{\tSpprtd}{\lower1.6pt\hbox{\tiny\it Spprtd}}
\newcommand{\tprec}{\lower1.6pt\hbox{\tiny$\prec$}}
\newcommand{\tcong}{\lower1.6pt\hbox{\tiny$\cong$}}
\newcommand{\tspr}{\lower3.5pt\hbox{\tiny$\spr$}}
\newcommand{\indexentry}[2]{#1 \hfill #2\newline}

\newcommand{\onetoone}{\hbox{1 \kern-3.2pt -- \kern-3.0pt 1}}

\newcommand{\onetoonen}
{\hbox
{1 \kern-3.2pt \rule[3.7pt]{4pt}{0.6pt} \kern-2.8pt 1}}

\newcommand{\overbbA}
               {\kern3.1pt
	       \overline{\kern-1.1pt\bbA\kern-0.6pt}\kern0.6pt}
\newcommand{\oversbbA}
               {\kern2.6pt
	       \overline{\kern-2.6ptsboldbbA\kern-0.3pt}\kern0.3pt}
\newcommand{\fbbA}{\kern1.3pt\bbA}
\newcommand{\fbbL}{\kern1.3pt\bbL}
\newcommand{\fbbT}{\kern1.3pt\bbT}
\newcommand{\overB}
               {\kern3.1pt\overline{\kern-3.1ptB\kern-0.6pt}\kern0.6pt}
\newcommand{\overrmB}
               {\kern3.1pt\overline{\kern-1.3pt{\rule{0pt}{8.7pt}\rm B\rule{0.7pt}{0pt}}\kern-1.7pt}\kern1.7pt}
\newcommand{\soverrmB}
               {\kern3.1pt\overline{\kern-1.3pt{\rule{0pt}{6.7pt}\rm B\rule{0.7pt}{0pt}}\kern-1.7pt}\kern1.7pt}
\newcommand{\overrmN}
               {\kern3.1pt\overline{\kern-1.3pt{\rule{0pt}{8.7pt}\rm N\rule{0.7pt}{0pt}}\kern-1.7pt}\kern1.7pt}
\newcommand{\overrmU}
               {\kern3.1pt\overline{\kern-1.3pt{\rule{0pt}{8.7pt}\rm U\rule{0.7pt}{0pt}}\kern-1.7pt}\kern1.7pt}
\newcommand{\overC}
               {\kern3.1pt\overline{\kern-2.5ptC\kern-0.6pt}\kern0.6pt}
\newcommand{\doverC}
               {\kern3.1pt\overline{\kern-3.1pt\overC\kern-0.6pt}\kern0.6pt}
\newcommand{\overE}
               {\kern3.1pt\overline{\kern-2.5ptE\rule{0pt}{8.5pt}\kern-0.6pt}\kern0.6pt}
\newcommand{\oversE}
               {\kern2.6pt\overline{\kern-2.6ptE\rule{0pt}{5.8pt}\kern-0.3pt}\kern0.3pt}
\newcommand{\overssE}
               {\kern2.1pt
	       \overline{\kern-1.2ptE\rule{0.3pt}{0pt}}
	       \kern-0.3pt}
\newcommand{\overF}
               {\kern3.1pt\overline{\kern-3.1ptF\kern-0.6pt}\kern0.6pt}
\newcommand{\oversF}
               {\kern2.6pt\overline{\kern-2.6ptF\kern-0.3pt}\kern0.3pt}
\newcommand{\overG}
               {\kern3.1pt\overline{\kern-3.1ptG\kern-0.6pt}\kern0.6pt}
\newcommand{\overH}
               {\kern3.1pt\overline{\kern-3.1ptH\kern-0.6pt}\kern0.6pt}
\newcommand{\overBH}
               {\kern1pt\overline{\kern-1pt\rule{0pt}{10pt}\rfs{BH}\kern-2.3pt}\kern2.3pt}
\newcommand{\overK}
               {\kern3.1pt\overline{\kern-3.1ptK\kern-0.6pt}\kern0.6pt}
\newcommand{\overL}
               {\kern3.1pt\overline{\kern-3.1ptL\kern-0.6pt}\kern0.6pt}
\newcommand{\overN}
               {\kern3.1pt\overline{\kern-3.1pt\rule{0pt}{8.8pt}N\kern-0.6pt}\kern0.6pt}
\newcommand{\overmN}
               {\kern1.6pt$\overline{\kern-1.3pt{\rule{0pt}{8.7pt}\rm N\rule{0.7pt}{0pt}}\kern-1.7pt}$\kern1.7pt}
\newcommand{\overT}
               {\kern3.1pt\overline{\kern-0.9pt\rule{0pt}{8.5pt}T\kern-0.3pt}\kern0.6pt}
\newcommand{\overU}
               {\kern2.3pt\overline{\kern-2.3pt\rule{0pt}{8.8pt}\rule{0.5pt}{0pt}U\kern-0.6pt}\kern0.6pt}
\newcommand{\overmU}
               {\kern1.6pt$\overline{\kern-1.3pt{\rule{0pt}{8.7pt}\rm U\rule{0.7pt}{0pt}}\kern-1.7pt}$\kern1.7pt}
\newcommand{\overV}
               {\kern3.1pt\overline{\kern-0.9ptV\kern-0.3pt}\kern0.6pt}
\newcommand{\overW}
               {\kern3.1pt\overline{\kern-0.9ptW\kern-0.3pt}\kern0.6pt}
\newcommand{\overX}
               {\kern3.1pt\overline{\kern-3.1ptX\kern-0.6pt}\kern0.6pt}
\newcommand{\oversG}
               {\kern2.6pt\overline{\kern-2.6ptG\kern-0.3pt}\kern0.3pt}
\newcommand{\oversX}
               {\kern2.6pt\overline{\kern-2.6ptX\kern-0.3pt}\kern0.3pt}

\newcommand{\overcalO}
               {\kern3.1pt\overline{\kern-3.1pt{\cal O}\kern-0.6pt}\kern0.6pt}
\newcommand{\overitPhi}
               {\kern3.1pt\overline{\kern-3.1pt\rule{0pt}{8.5pt}\itPhi\kern-0.6pt}\kern0.6pt}
\newcommand{\oversitPhi}
               {\kern3.1pt\overline{\kern-3.1pt\rule{0pt}{5.8pt}\itPhi\kern-0.6pt}\kern0.6pt}
\newcommand{\overitPsi}
               {\kern3.1pt\overline{\kern-3.1pt\rule{0pt}{8.5pt}\itPsi\kern-0.6pt}\kern0.6pt}
\newcommand{\oversitPsi}
               {\kern3.1pt\overline{\kern-3.1pt\rule{0pt}{5.8pt}\itPsi\kern-0.6pt}\kern0.6pt}
\newcommand{\overleq}{\mathrel{\hbox{\kern2pt\tiny $\overline{\rule{0pt}{5.1pt}\kern-3.0pt\leq\kern-0.5pt}$}\kern1.3pt}}

\newcommand{\sstar}{\hbox{\scriptsize \kern1pt$\star$\kern0pt}}
\newcommand{\rsstar}{\raise1pt\hbox{\scriptsize \kern1pt$\star$\kern1pt}}
\newcommand{\raisedsstar}{\kern0.5pt\raise1.0pt\hbox{\scriptsize$*$}}
\newcommand{\sperp}{\raise2pt\hbox{\kern0pt\tiny $\perp$}}
\newcommand{\ssperp}{\raise1.5pt\hbox{\kern0pt\tinier \symbol{63}}}

\newcommand{\nperp}
{
\hbox{\kern1.0pt
\hbox{\raise1.0pt \hbox{\rule{4.6pt}{0.3pt}}}\kern-5.1pt$\sperp$}
}
\newcommand{\snperp}
{
\hbox{\kern1.0pt
\hbox{\raise0.4pt \hbox{\rule{3.3pt}{0.2pt}}}\kern-3.65pt$\ssperp$}
}

\newcommand{\sor}{\raise2pt\hbox{\tiny\rm Or}}

\newcommand{\srs}[1]{\raise1.3pt\hbox{\tiny\rm #1}}
\newcommand{\ssrfs}[1]{\raise1.3pt\hbox{\tiny\rm #1}}
\newcommand{\ssrs}[1]{\raise1.3pt\hbox{\nmini #1}}
\newcommand{\scirc}
{\raise1pt\hbox{\scriptsize\kern1.5pt$\circ$\kern1.5pt}}
\newcommand{\sscirc}
{\hbox{\tiny\kern1.5pt$\circ$\kern0.3pt}}
\newcommand{\tcirc}
{\hbox{\tiny\kern1.5pt$\circ$\kern0.6pt}}

\newcommand{\bcirc}
{\mathop{\lower1pt\hbox{\large\kern1.5pt$\circ$}}}
\newcommand{\bbcirc}
{\mathop{\lower1.8pt\hbox{\Large\kern1.5pt$\circ$}}}

\newcommand{\raisedstar}
{\raise-3.1pt\hbox{$\rule{0.5pt}{0pt}^*$}}
\newcommand{\rstar}
{\raise-3.1pt\hbox{$\rule{0.5pt}{0pt}^*$}}
\newcommand{\raisedtinystar}
{\kern0.7pt\raise1.5pt\hbox{\tiny $*$}}
\newcommand{\rtstar}
{\kern0.7pt\raise1.0pt\hbox{\tiny $*$}\kern-2.0pt}
\newcommand{\bfia}{\hbox{\blcmssifont \symbol{97}\kern1pt}}
\newcommand{\vecbfia}{\vec{\bfia}}

\newcommand{\bfib}{\hbox{\blcmssifont \symbol{98}\kern1pt}}

\newcommand{\bfic}{\hbox{\blcmssifont \symbol{99}\kern1pt}}

\newcommand{\bfid}{\hbox{\blcmssifont \symbol{100}\kern1pt}}

\newcommand{\bfie}{\hbox{\blcmssifont \symbol{101}\kern1pt}}

\newcommand{\bfif}{\hbox{\blcmssifont \symbol{102}\kern1pt}}
\newcommand{\vecbfif}{\vec{\bfif\rule{1pt}{0pt}}}
\newcommand{\sbfif}{\hbox{\lcmssifont \symbol{102}}}
\newcommand{\bfig}{\hbox{\blcmssifont \symbol{103}}}
\newcommand{\vecbfig}{\vec{\bfig}}
\newcommand{\sbfig}{\hbox{\lcmssifont \symbol{103}}}
\newcommand{\bfih}{\hbox{\blcmssifont \symbol{104}\kern1pt}}
\newcommand{\vecbfih}{\vec{\bfih}}

\newcommand{\bfii}{\hbox{\blcmssifont \symbol{105}}}
\newcommand{\bfij}{\hbox{\kern1.7pt\blcmssifont \symbol{106}\kern1.7pt}}
\newcommand{\fbfij}{\hbox{\kern1.7pt\blcmssifont \symbol{106}\kern1.7pt}}

\newcommand{\bfik}{\hbox{\blcmssifont \symbol{107}\kern1pt}}
\newcommand{\fbfik}{\hbox{\kern1.7pt\blcmssifont \symbol{107}\kern1.7pt}}

\newcommand{\bfil}{\hbox{\kern0.5pt\blcmssifont \symbol{108}\kern1.7pt}}

\newcommand{\bfip}{\hbox{\blcmssifont \symbol{112}}}

\newcommand{\bfiu}{\hbox{\blcmssifont \symbol{117}}}
\newcommand{\bfiv}{\hbox{\blcmssifont \symbol{118}}}
\newcommand{\bfiw}{\hbox{\blcmssifont \symbol{119}}}
\newcommand{\bfispw}{\bfiw\kern2.0pt}
\newcommand{\bfix}{\hbox{\blcmssifont \symbol{120}}}
\newcommand{\vecbfix}{\vec{\hbox{\blcmssifont \symbol{120}}}}
\newcommand{\bfiy}{\hbox{\blcmssifont \symbol{121}}}
\newcommand{\vecbfiy}{\vec{\bfiy}}
\newcommand{\bfiz}{\hbox{\blcmssifont \symbol{122}}}

\newcommand{\bfiI}{\hbox{\blcmssifont \symbol{73}\kern2.3pt}}

\newcommand{\mcdot}{\kern-1.4pt\cdot\kern-1.4pt}
\newcommand{\lcdot}{\hbox{$\kern0.5pt\cdot\kern0.9pt$}}
\newcommand{\ncdot}{\hbox{$\kern-0.2pt\cdot\kern0.6pt$}}
\newcommand{\bcdot}{\mathbin{\lower1.2pt\hbox{\Large$\cdot$}}}
\newcommand{\fr}{\kern1.0ptr}
\newcommand{\fzero}{\kern1.3pt0}
\newcommand{\fone}{\kern1.0pt1}
\newcommand{\ftwo}{\kern1.3pt2}
\newcommand{\fprime}{\kern1.3pt'}
\newcommand{\Fprime}{\kern2pt'}
\newcommand{\fdprime}{\kern1.3pt''}
\newcommand{\fsscr}[1]{\kern1pt#1}
\newcommand{\ffsscr}[1]{\kern2pt#1}
\newcommand{\fstar}{\kern1.7pt*}
\newcommand{\ffprime}{\kern2.8pt'}
\newcommand{\farsu}[1]{^{\kern1.3pt#1}}
\newcommand{\sprime}{\raise1.3pt\hbox{\scriptsize\kern1.4pt$'$}}
\newcommand{\fprimesub}{\kern1.5pt'\kern-3.8pt}
\newcommand{\fprimew}{\kern0.7pt'}
\newcommand{\fprimei}{\kern1.5pt'\kern-3.0pt}
\newcommand{\fnprime}{\kern1.5pt'\kern-3.0pt}

\newcommand{\clubsign}{{\raise1.3pt\hbox{\tiny $\clubsuit$}}}
\newcommand{\sclub}{{\raise1.3pt\hbox{\tiny $\clubsuit$}}}
\newcommand{\srsp}[1]{\hbox{\kern1pt\tiny $(#1)$}}
\newcommand{\sminus}{\rule{5.6pt}{0.3pt}}
\newcommand{\ssim}{\hbox{\scriptsize $\sim$}}
\newcommand{\neweq}{
\mathrel{
\hbox{
\lower-0.0pt\hbox{\sminus}
\kern-9.95pt\ssim\kern-6.2pt\raise4.4pt\hbox{\sminus}
}}}

\newcommand{\dcup}{\mathbin{\hbox{$\cup$
\kern-10.43pt\raise1.4pt \hbox{\tiny $\cup$}}}}

\newcommand{\sprt}[1]
{\lower1.8pt\hbox{\kern1.5pt\rule{0.5pt}{10pt}}
\kern-0.5pt\underline{\kern1pt#1\kern1pt}
\lower1.8pt\hbox{\kern-0.2pt\rule{0.5pt}{10pt}}\kern1pt}

\newcommand{\sprtd}[2]
{#1\lower1.3pt\hbox{\kern1.5pt\rule{0.5pt}{10pt}}
\kern-0.5pt\underline{\kern1pt#2\kern1pt}
\lower1.3pt\hbox{\kern0.2pt\rule{0.5pt}{10pt}}\kern1pt}

\newcommand{\sprtlm}[1]
{\lower3.4pt\hbox{\kern1.5pt\rule{0.5pt}{12pt}}
\kern-0.5pt\underline{\kern1pt#1\kern1pt}
\lower3.4pt\hbox{\kern-0.5pt\rule{0.5pt}{12pt}}\kern1pt}

\newcommand{\sprtl}[1]
{\lower4.4pt\hbox{\kern1.5pt\rule{0.5pt}{14pt}}
\kern-0.5pt\underline{\kern1pt#1\kern1pt}
\lower4.4pt\hbox{\kern-0.5pt\rule{0.5pt}{14pt}}\kern1pt}

\newcommand{\sprtlk}[1]
{\lower4.8pt\hbox{\kern1.5pt\rule{0.5pt}{14pt}}
\kern-0.5pt\underline{\kern1pt#1\kern2.5pt}
\lower4.8pt\hbox{\kern-0.5pt\rule{0.5pt}{14pt}}\kern1pt}

\newcommand{\sprtll}[1]
{\lower5.7pt\hbox{\kern1.5pt\rule{0.5pt}{15.4pt}}
\kern-0.5pt\underline{\kern1pt#1\kern1pt}
\lower5.7pt\hbox{\kern-0.5pt\rule{0.5pt}{15.4pt}}\kern1pt}

\newcommand{\sprtln}[1]
{\lower4.8pt\hbox{\kern1.5pt\rule{0.5pt}{15.4pt}}
\kern-0.5pt\underline{\kern1pt#1\kern1pt}
\lower4.8pt\hbox{\kern-0.5pt\rule{0.5pt}{15.4pt}}\kern1pt}

\newcommand{\sprtdl}[2]
{#1\lower4.0pt\hbox{\kern1.5pt\rule{0.5pt}{14pt}}
\kern-0.5pt\underline{\kern1pt#2\kern1pt}
\lower4.1pt\hbox{\kern0.0pt\rule{0.5pt}{14pt}}\kern1pt}

\newcommand{\sprtm}[1]
{\lower3.3pt\hbox{\kern1.5pt\rule{0.5pt}{14pt}}
\kern-0.5pt\underline{\kern1pt#1\kern1pt}
\lower3.4pt\hbox{\kern0.0pt\rule{0.5pt}{14pt}}\kern1pt}

\newcommand{\sprtdm}[2]
{#1\lower3.3pt\hbox{\kern1.5pt\rule{0.5pt}{14pt}}
\kern-0.5pt\underline{\kern1pt#2\kern1pt}
\lower3.4pt\hbox{\kern0.0pt\rule{0.5pt}{14pt}}\kern1pt}

\DeclareSymbolFont{AMSb}{U}{msb}{m}{n}
\DeclareMathSymbol{\bbA}{\mathbin}{AMSb}{"41}
\DeclareMathSymbol{\bbB}{\mathbin}{AMSb}{"42}
\DeclareMathSymbol{\bbC}{\mathbin}{AMSb}{"43}
\DeclareMathSymbol{\bbD}{\mathbin}{AMSb}{"44}
\DeclareMathSymbol{\bbE}{\mathbin}{AMSb}{"45}
\DeclareMathSymbol{\bbF}{\mathbin}{AMSb}{"46}
\DeclareMathSymbol{\bbG}{\mathbin}{AMSb}{"47}
\DeclareMathSymbol{\bbH}{\mathbin}{AMSb}{"48}
\DeclareMathSymbol{\bbI}{\mathbin}{AMSb}{"49}
\DeclareMathSymbol{\bbJ}{\mathbin}{AMSb}{"4A}
\DeclareMathSymbol{\bbK}{\mathbin}{AMSb}{"4B}
\DeclareMathSymbol{\bbL}{\mathbin}{AMSb}{"4C}
\DeclareMathSymbol{\bbM}{\mathbin}{AMSb}{"4D}
\DeclareMathSymbol{\bbN}{\mathbin}{AMSb}{"4E}
\DeclareMathSymbol{\bbO}{\mathbin}{AMSb}{"4F}
\DeclareMathSymbol{\bbP}{\mathbin}{AMSb}{"50}
\DeclareMathSymbol{\bbQ}{\mathbin}{AMSb}{"51}
\DeclareMathSymbol{\bbR}{\mathbin}{AMSb}{"52}
\DeclareMathSymbol{\bbS}{\mathbin}{AMSb}{"53}
\DeclareMathSymbol{\bbT}{\mathbin}{AMSb}{"54}
\DeclareMathSymbol{\bbU}{\mathbin}{AMSb}{"55}
\DeclareMathSymbol{\bbV}{\mathbin}{AMSb}{"56}
\DeclareMathSymbol{\bbW}{\mathbin}{AMSb}{"57}
\DeclareMathSymbol{\bbX}{\mathbin}{AMSb}{"58}
\DeclareMathSymbol{\bbY}{\mathbin}{AMSb}{"59}
\DeclareMathSymbol{\bbZ}{\mathbin}{AMSb}{"5A}
\newcommand{\boldbbA}{\mathbb{A}\kern-9pt\mathbb{A}}
\newcommand{\boldbbL}{\mathbb{L}\kern-8.3pt\mathbb{L}}
\newcommand{\boldbbR}{\mathbb{R}\kern-8.3pt\mathbb{R}}

\newcommand{\sboldbbA}
                {\hbox{\kern0.7pt\scriptsize
		$\mathbb{A}\kern-6pt\mathbb{A}$}}
\newcommand{\sboldbbL}
                {\hbox{\kern0.7pt\scriptsize
		$\mathbb{L}\kern-5.68pt\mathbb{L}$}}
\newcommand{\sboldbbN}
                {\hbox{\kern0.7pt\scriptsize
		$\mathbb{N}\kern-5.60pt\mathbb{N}$}}
\newcommand{\sboldbbR}
                {\hbox{\kern0.7pt\scriptsize
		$\mathbb{R}\kern-5.68pt\mathbb{R}$}}
\newcommand{\sboldbbT}
                {\hbox{\kern0.7pt\scriptsize
		$\mathbb{T}\kern-5.60pt\mathbb{T}$}}

\newcommand{\itbfI}{\boldsymbol{I}}

\newcommand{\itbfT}{\boldsymbol{T}}

\newcommand{\hatitbfI}{\hat{\itbfI}}

\newcommand{\wpm}{\kern0.7pt\pm}
\newcommand{\doubledagger}{\hbox{$\dagger\kern-4pt\dagger$}}

\newcommand{\lip}{\hbox{\bfil\bfii\kern-0.5pt\bfip}}

\newcommand{\wtildebbU}{\kern1pt\widetilde{\bbU}}
\newcommand{\wbbU}{\kern1pt\bbU}
\newcommand{\bfbbU}{\bbU\kern-12pt\bbU}
\newcommand{\sbfbbU}{\bbU\kern-10pt\bbU}

\newcommand{\eclass}{\hbox{$\kern1.8pt/\kern-4.0pt$}}

\newcommand{\overrfs}[1]{\overline{\rfs{#1}}}
\def\overrfs#1{\mbox{\kern2pt$\overline{\kern-2pt\mbox{\kern0.7pt\rm #1\kern 1.3pt}\kern-2pt}$\kern2pt}}

\newcommand{\drestriction}
{\mathbin{
\lower0.3pt\hbox{\small$\upharpoonright$}
\kern-5.5pt\raise2.0pt\hbox{\small$\upharpoonright$}}}

\newcommand{\sprtk}[1]
{\lower1.9pt\hbox{\kern1.5pt\rule{0.5pt}{11pt}}
\kern-0.5pt\underline{\kern2.6pt#1\kern3.0pt}
\kern-0.4pt\lower1.9pt\hbox{\kern0.2pt\rule{0.5pt}{11pt}}\kern1pt}

\newcommand{\farslash}{\kern1.2pt\slash\kern-3pt}

\newcommand{\sNu}{\lower1.5pt\hbox{\small$\nu$}}

\newcommand{\fsup}[1]{^{\kern1.5pt #1}}
\newcommand{\Fsup}[1]{^{\kern2.0pt #1}}
\newcommand{\tinysup}[1]{^{\kern1pt\raise2pt\hbox{\tiny $#1$}}}
\newcommand{\fsub}[1]{_{\kern1.3pt #1}}
\newcommand{\Fsub}[1]{_{\kern2.3pt #1}}
\newcommand{\nsup}[1]{^{\kern-1.0pt #1}}
\newcommand{\Nsup}[1]{^{\kern-2.3pt #1}}
\newcommand{\nsub}[1]{_{\kern-1.0pt #1}}
\newcommand{\Nsub}[1]{_{\kern-2.3pt #1}}

\newcommand{\aoverb}[2]{\left(\kern-14pt\begin{array}{l}#1\\#2\end{array}\kern-16pt\right)}
\newcommand{\dvec}[1]{\stackrel{\lower1pt\hbox{\tiny$\Rightarrow$}}{#1}\rule{0.0pt}{0pt}\kern-4pt}

\newcommand{\concat}{\mathbin{\kern0.7pt\widehat{\ }\kern0.7pt}}

\newcommand{\iseg}{\mathrel{\vbox{\hbox{$< \kern-5.4pt \lower 2.5pt\hbox{\tiny $|$} \kern -2.3pt \raise 5.0pt\hbox{\tiny $|$}$}}\kern 0.5pt}}
\newcommand{\isegeq}{\mathrel{\vbox{\hbox{$\leq \kern-5.4pt \lower 2.5pt\hbox{\tiny $|$} \kern -2.3pt \raise 5.0pt\hbox{\tiny $|$}$}}\kern 0.5pt}}
\newcommand{\risegeq}{\mathrel{\vbox{\hbox{$\geq \kern-12.8pt \lower 2.5pt\hbox{\tiny $|$} \kern -2.3pt \raise 5.0pt\hbox{\tiny $|$}$}}\kern 5.9pt}}
\newcommand{\reverseisegeq}{\mathrel{\vbox{\hbox{$\geq \kern-12.8pt \lower 2.5pt\hbox{\tiny $|$} \kern -2.3pt \raise 5.0pt\hbox{\tiny $|$}$}}\kern 5.9pt}}
\newcommand{\riseg}{\mathrel{\vbox{\hbox{$> \kern-12.8pt \lower 2.5pt\hbox{\tiny $|$} \kern -2.3pt \raise 5.0pt\hbox{\tiny $|$}$}}\kern 5.9pt}}
\newcommand{\sriseg}
{\mathrel{\vbox{\hbox{\tiny$> \kern-8.8pt \lower 3.5pt\hbox{\tiny $\vert$} \kern -2.3pt \raise 3.3pt\hbox{\tiny $\vert$}$}}\kern 5.9pt}}
\def\Tiny={\mbox{\tiny $=$}}
\def\hollowsquare{\hbox{${\vcenter{\vbox{
   \hrule height 0.4pt\hbox{\vrule width 0.4pt height 6pt
   \kern5pt\vrule width 0.4pt}\hrule height 0.4pt}}}$}}
\def\oldhollowqed{\hbox{${\vcenter{\vbox{
   \hrule height 0.4pt\hbox{\vrule width 0.4pt height 6pt
   \kern5pt\vrule width 0.4pt}\hrule height 0.4pt}}}$}}
\newcommand{\smallhollowqed}
{\indent\indent$
\mathrel{
\hbox{
\hbox{
\rule{0.6pt}{7pt}
\kern-4.5pt
\raise6.3pt
\hbox{\rule{5.8pt}{0.6pt}}
\kern-4.5pt
\rule{0.6pt}{7pt}
\kern-09.5pt
\rule{5.8pt}{0.6pt}
\kern-06pt
}
}
}
$}
\newcommand{\hollowqed}{\hfill\hbox{\small$\Box$}}
\newcommand{\dhollowqed}{\hfill\hbox{\small$\Box$}\kern2pt\hbox{\small$\Box$}}

\newcommand{\tinytimes}{\mathbin{\raise1.2pt\hbox{\tiny $\sttimes$}}}

\newcommand{\nnot}[1]{\mathrel{\hbox{\kern-0.6pt$\not \kern-3pt #1 \kern0pt$}}}
\newcommand{\ubracketX}{\hbox{$\raise1.7pt\hbox{$\ulcorner$}\kern-4pt \overX \kern-3pt\raise1.7pt\hbox{$\urcorner$}$}}
\newcommand{\ubracket}[1]{\hbox{$\raise1.7pt\hbox{$\ulcorner$}\kern-4pt \overline{\rule{0pt}{8.5pt}#1} \kern-3pt\raise1.7pt\hbox{$\urcorner$}$}}
\newcommand{\uubracket}[1]{\hbox{$\raise2.4pt\hbox{$\ulcorner$}\kern-4pt \overline{\rule{0pt}{8.5pt}#1} \kern-3pt\raise2.4pt\hbox{$\urcorner$}$}}
\newcommand{\uuubracket}[1]{\hbox{$\raise2.7pt\hbox{$\ulcorner$}\kern-4pt \overline{\rule{0pt}{8.5pt}#1} \kern-3pt\raise2.7pt\hbox{$\urcorner$}$}}

\newcommand{\myboldminus}{\rule{0pt}{0pt}\raise2pt\hbox{\scriptsize$\boldsymbol{-}$}}
\newcommand{\vbminus}{\rule{0pt}{0pt}\raise2pt\hbox{\scriptsize$\boldsymbol{-}$}}
\newcommand{\vvbminus}{\rule{0pt}{0pt}\raise1.52pt\hbox{\scriptsize$\boldsymbol{-}$} \kern-7.5pt \raise2pt\hbox{\scriptsize$\boldsymbol{-}$}}

\newcommand{\vbcap}{\boldsymbol{\hbox{\small $\bigcap \kern-10.5pt \bigcap$}}}

\newcommand{\bfminus}{\boldsymbol{-}}
\newcommand{\veryboldminus}{\bfminus\kern-11.0pt\raise0.4pt\hbox{$\bfminus$}}

\newcommand{\boldexists}{\boldsymbol{\exists}}
\newcommand{\boldforall}{\boldsymbol{\forall}}
\newcommand{\veryboldexists}{\boldexists\kern-7.5pt\raise0.3pt\hbox{$\boldexists$}}
\newcommand{\vbexists}{\boldexists\kern-7.5pt\raise0.3pt\hbox{$\boldexists$}}
\newcommand{\vbforall}{\boldforall\kern-7.5pt\raise0.3pt\hbox{$\boldforall$}}

\newcommand{\smallvecone}{\lower0.3pt\hbox{$\vec{\hbox{\sbfreg 1}}$\rule{0pt}{06pt}}}

\newcommand{\fcap}{\kern0.6mm\cap\kern0.6mm}
\newcommand{\fY}{\kern0.6mmY}
\newcommand{\fequal}{\kern4.0pt=\kern4.0pt}

\newcommand{\mysetminus}{\mathbin {\raise1pt\hbox{\mbox{\scriptsize\kern2pt $\setminus$\kern2pt}}}}
\newcommand{\myover}[1]{\slash\kern-2pt #1}
\newcommand{\myoversim}{\slash\kern-2pt \sim}
\newcommand{\mynarrowover}[1]{\kern-2pt#1\kern-2pt \sim}
\newcommand{\mynarrowoversim}{\kern-2pt\slash\kern-2pt \sim}
\newcommand{\hcli}[2]{[\kern1pt #1,#2\kern1pt)}
\newcommand{\cli}[2]{[\kern1pt #1,#2\kern1pt]}

\newcommand{\lldot}{\mathrel{\hbox{$\ll\kern-7.3pt\cdot$}}}
\newcommand{\lessdotnear}{\lessdot^{\trfs{near}\kern-2.5pt }}
\newcommand{\equivnear}{\equiv^{\trfs{near}\kern-2.5pt }}

\newcommand{\supminus}{^{\kern1.5pt -}}
\newcommand{\subminus}{_{\kern1.5pt -}}

\newcommand{\underlinedC}{\kern1pt\underline{\kern-0ptC\kern-3pt}\kern2.5pt}
\newcommand{\underlinedU}{\kern1pt\underline{\kern-0ptU\kern-3pt}\kern2.5pt}

\newcommand{\hatoperation}{\kern3pt\widehat{\rule{3pt}{0pt}}\kern2pt}

\newcommand{\wcheck}[1]{\widecheck{#1}}
\newcommand{\bfMu}{\lower0.0pt\hbox{\large${\boldsymbol\mu}$\kern0.7pt}}
\newcommand{\bfRho}{\lower0.0pt\hbox{\large${\boldsymbol\rho}$\kern0.7pt}}
\newcommand{\bfUpsilon}{\lower0.0pt\hbox{\large${\boldsymbol\upsilon}$\kern0.7pt}}
\newcommand{\Mu}{\lower1.3pt\hbox{\mbox{\large${\mathbf\mu}$\kern0.7pt}}}
\newcommand{\mynot}[2]{\kern#1pt/\kern-#2pt}
\newcommand{\wcheckvec}[1]{(\vec#1)\kern1pt\wcheck{\ }}
\newcommand{\widesetminus}{\mathbin{\kern1.0pt\setminus\kern1.0pt}}
\newcommand{\widein}{\kern1.5pt\in\kern1.5pt}
\newcommand{\tinyrightarrow}{\raise1pt\hbox{\mbox{\tiny \kern1pt$\rightarrow$\kern1pt}}}
\newcommand{\supprtd}[1]
{
\kern2pt\lower2pt\hbox{\rule{0.7pt}{11.0pt}}
\underline{\rule{1pt}{0pt}#1\rule{1pt}{0pt}}
\lower2pt\hbox{\rule{0.7pt}{11.0pt}\kern1.7pt}
}
\newcommand{\supprtdsub}[1]
{
\kern2pt\lower3.8pt\hbox{\rule{0.7pt}{13.5pt}}
\underline{\rule{1pt}{0pt}#1\rule{1pt}{0pt}}
\lower3.8pt\hbox{\rule{0.7pt}{13.5pt}\kern1.7pt}
}
\newcommand{\vwidesetminus}{\mathbin{\kern2pt\setminus\kern2pt}}

\newcommand{\dhat}[1]{\kern12pt\widehat{\hat{#1}}}
\newcommand{\dhatI}{\hbox{$\hatI$\kern-1.6pt\raise5.0pt\hbox{$\hat{ }$}}}

\newcommand{\notmodels}{\mathrel{\models\kern-11pt\backslash \kern6pt}}

\newcommand{\binslash}{\mathbin{\kern2pt\slash\kern-3pt}}
\newcommand{\supunderbracket}[1]{\raise3pt\hbox{\kern1pt\scriptsize$\underbracket{#1}$}}
\def\boldf1{{\bf 1}}

\newcommand\fakeslant[1]{%
  \pdfliteral{1 0 0.167 1 0 0 cm}#1\pdfliteral{1 0 -0.167 1 0 0 cm}}
\newcommand\mathbbsl[1]{\mathbb{\fakeslant{#1}}}

\newcommand{\wsetminus}{\mathbin{\rule{2pt}{0pt}\setminus\rule{2pt}{0pt}}}
\newcommand{\swsetminus}{\mathbin{\rule{0.0pt}{0pt}\setminus\rule{0.0pt}{0pt}}}

\newcommand{\itbbI}{\mathbbsl{I}}
\newcommand{\hatitbbI}{\hat{\mathbbsl{I}\rule{2.6pt}{0pt}}\kern-1pt}
\newcommand{\baritbbI}{\bar{\mathbbsl{I}\rule{2.6pt}{0pt}}\kern-1pt}
\newcommand{\whatitbbI}{\widehat{\mathbbsl{I}\rule{2.6pt}{0pt}}\kern-1pt}

\newcommand{\dhatitbfI}{\hbox{$\hatitbfI$\kern-1.8pt\raise5.0pt\hbox{$\hat{ }$}}}
\newcommand{\dhatitbbI}{\hat{\itbbI\rule{2pt}{0pt}}\kern-5pt\raise4pt\hbox{$\widehat{\ }$}}

\newcommand{\adm}{\mathbin{+^{\kern-1.6pt\hbox{\tiny $\itbfT$}}}}
\newcommand{\sbm}{\mathbin{-^{\kern-1.2pt\hbox{\tiny $\itbfT$}}}}
\newcommand{\mlm}{\bcdot^{\kern-1.6pt\hbox{\tiny $\itbfT$}}}

\newcommand{\modminus}{-^{\kern-0.5pt\hbox{\tiny $\hatitbfI$}}}
\newcommand{\modplus}{+^{\kern-0.5pt\hbox{\tiny $\hatitbfI$}}}
\newcommand{\modbcdot}{\bcdot^{\kern-0.5pt\hbox{\tiny $\hatitbfI$}}}
\newcommand{\sdhatitbbI}{\hat{\itbbI\rule{2pt}{0pt}}\kern-4pt\raise3pt\hbox{\scriptsize$\widehat{\ }$}}

\newcommand{\dbarrestriction}{\mathbin{\lower2.5pt\hbox{\rule{0.5pt}{10.5pt}}\kern-5.5pt\upharpoonright\kern-2.5pt}}

\newcommand{\ff}[1]{\hbox{\it #1} \kern 0.7pt}
\newcommand{\sff}[1]{\hbox{\scriptsize\it #1} \kern 0.4pt}

\newcommand{\length}{{\scriptsize\rm lngth}}

\newcommand{\thinrd}[1]{\kern -1pt #1}
\newcommand{\thinld}[1]{#1  \kern -1pt}

\newcommand{\wrdbrack}{\kern1pt\rdbrack}
\newcommand{\wldbrack}{\kern1pt\ldbrack}

\newcommand{\fstneg}{\stneg\rule{1pt}{0pt}}

\usepackage{mathrsfs}

\makeindex
\makeatletter
\newcommand{\leqnomode}{\tagsleft@true}
\newcommand{\reqnomode}{\tagsleft@false}
\makeatother
\begin{document}
\baselineskip 18pt
\pagestyle{empty}
\centerline{\large\bf Reconstruction theorems for semigroups of functions}
\centerline{\large\bf which contain all transpositions of a set}
\centerline{\large\bf and for clones with the same property}
\vspace{1cm}
\centerline{\bf Jonah Maissel and Matatyahu Rubin}
\vspace{0.1cm}
\centerline{\bf Ben Gurion University, Beer Sheva ISRAEL}
\vspace{1cm}
\centerline{\bf August 2015}
\vspace{5mm}

\begin{abstract}
$A^B$ denotes the set of functions from $B$ to $A$.
A subset $S \subseteq A^A$ is called a {\it function semigroup},
if $S$ is closed under composition.
A permutation of $A$ which moves exactly to elements
of $A$ is called a {\it transposition} of $A$. Let $S \subseteq A^A$ be a
function semigroup. $S$ is {\it fully-transpositional}, if every transposition
of $A$ belongs to $S$.
\vspace{1.5mm}
\\
{\bf Theorem 1:} Let $A,B$ be sets such that $\abs{A},\abs{B} \neq 0,1,2,6$
and $S \subseteq A^A$, $T \subseteq B^B$
be fully-transpositional function semigroups.
Suppose that $\varphi$ is an isomorphism between
$\pair{S}{\circ}$ and $\pair{T}{\circ}$.
Then there is a bijection $\tau$ from $A$ to $B$ such that for every $g \stin S$,
$\varphi(g) = \tau \circ g \circ\tau\inverse$.%
\vspace{1.5mm}
\\
Theorem 1 strengthens a theorem of R. Mckenzie who proved the same theorem
for permutation groups.
\vspace{1.5mm}
\\
Let $f \stin A^{A^n}$, $g_1,\ldots,g_n \stin A^{A^k}$.
The notation $f \circ^A_{n,k}(g_1,\ldots,g_n)$
stands for the {\it $(n,k)$-composition-operation} on $A$, which is defined as
follows:
$f \circ^A_{n,k}(g_1,\ldots,g_n)(\vecx) \eqdf f(g_1(\vecx),\ldots,g_k(\vecx))$,
\ $\vecx \stin A^k$.
For $1 \leq i \leq n$ the {\it $(n,i)$-coordinate-function} on $A$ is defined by
$p^A_{n,i}(x_1,\ldots,x_n) \eqdf x_i$, \ $x_1,\ldots,x_n \stin A$.
A subset $C \subseteq \bigcup_{i = m}^{\infty} A^{A^m}$ containing all the
$p^A_{n,i}$'s and closed under all the $\circ^A_{n,k}$'s
is called a {\it clone} on $A$.
A clone $C$ on $A$ is {\it fully-transpositional},
if every transposition of $A$ belongs to $C$.
Let $g \in A^{A^n}$ and $\tau$ be a bijection from $A$ to $B$.
Then $g^{\tau}$ is defined by: $g^{\tau} \in B^{B^n}$ and
$g^{\tau}(\vecx) \eqdf
\tau(g \tau\inverse(\vecx),\ldots,g \tau\inverse(\vecx))$,\ 
where $\vecx \in B^n$ denotes $(x_1,\ldots,x_n)$ and
$\tau\inverse(\vecx) \eqdf (\tau\inverse(x_1),\ldots,\tau\inverse(x_n))$.
\vspace{1.5mm}
\\
{\bf Theorem 2:} Let $A,B$ be sets such that $\abs{A},\abs{B} \neq 0,1,2,6$
and $C,D$ be fully-transpositional clones on $A$ and $B$.
Suppose that $\varphi$ is an\break isomorphism between the structures
$\pair{C}{\sngltn{\circ^A_{n,k}}_{n,k \stin \bbN^+}}$
and\\
$\pair{D}{\sngltn{\circ^B_{n,k}}_{n,k \stin \bbN^+}}$.
Then there is a bijection $\tau$ from $A$ to $B$ such that for every $g \stin C$,
$\varphi(g) = g^{\tau}$.
\end{abstract}

\noindent
\vspace{3mm}
\tableofcontents
\baselineskip18pt
\pagestyle{plain}
\setcounter{page}{1}
\section{\bf Introduction}\label{s1}
A function semigroup on $A$ is a nonempty subset of $A^A$
closed under composition.
It is ``fully-transpositional'', if it contains all transpositions of $A$.

This work deals with two classes of objects:
``fully-transpositional function semigroups'',
and ``fully-transpositional clones''.
Each object $X$ of one of the above types gives rise to an
algebraic structure which is denoted by $\rfs{Alg}(X)$.
Also, each such $X$ acts on a certain set $A$, and this action is described by
a structure which is denoted by $\rfs{Act}(X)$.
The goal of this work is to show for each of these objects,
$\rfs{Act}(X)$ is recoverable from $\rfs{Alg}(X)$ ``using first order formulas''.
(This phrase will be explained later.)
\vspace{3mm}

\noindent
{\bf Where these questions come from?}

Our interest in the topic of reconstructing the a action of a semigroup or
a clone from its algebraic structure, originates from a theorem of
R. McKenzie \cite{McK} and from a work of
M. Bodirsky, M. Pinsker, and A. Pongra\'cz \cite{BPP}.
The two related theorems proved in this work
strengthen the theorem of McKenzie, and also seem to strengthen a theorem
from \cite{BPP}.

Recently, Robert Barham \cite{Ba}, too,
had progress on the topic of reconstructing the action of a clone
from its algebraic structure.
He proved a reconstruction theorem for a certain class of polymorphic
clones of locally moving Boolean algebras. His theorem has applications
to the polymorphic clones of some other types of structures.

Notwithstanding the reconstruction facts that have been already proved,
the general problem of reconstructing the action of a function semigroup
or a clone from its algebraic structure still remains widely open.
\vspace{2mm}

We start with the Theorem of McKenzie. (See \cite{McK} Corollary 2.)
\vspace{2mm}

\noindent
{\bf Theorem } (R. Mckenzie) Let $G$ be a group of permutations
of a set $A$ containing every transposition of $A$,
and $H$ be a group of permutations of a set $B$
containing every transposition of $B$.
Suppose $\abs{A} \neq 6$. Assume further that $\varphi$ is an isomorphism
between $\pair{G}{\circ^G}$ and $\pair{H}{\circ^H}$.
\underline{Then} there is a bijection $\alpha$ between $A$ and $B$ such that for
every $g \stin G$, $\varphi(g) = \alpha \circ g \circ \alpha\inverse$.
\vspace{3mm}

In fact, the above theorem is only a corollary
from what Mckenzie actually proved.
In the terminology of this work (Definition~\ref{d1.1-08-14}),
McKenzie proved that the permutation group structure of $G$
is first order strongly interpretable in $\pair{G}{\circ^G}$.
(The permutation group structure of $G$ is $\pair{G \stcup A}{\rfs{App}^G}$,
where $\rfs{App}^G$ is defined as follows:
(i) $\fnn{\rfs{App}^G}{G \sttimes A}{A}$,
and (ii) for every $g \stin G$ and $a \stin A$, $\rfs{App}^G(g,a) = g(a)$.)

The proof of Mckenzie's Theorem is short and elegant.
It amounts to proving the existence of three first order formulas:
(1) $\phi_{\srfs{Trns}}(\bfif\kern2.3pt)$ which defines in $\pair{G}{\circ^G}$
the set of transpositions of $A$;
(2) $\phi_{\srfs{Rep}}(\bfif,\bfig\kern1pt)$
which expresses the fact that the supports of the
transpositions $\bfif$ and $\bfig$ intersect in a singleton;
and (3) $\phi_{\srfs{Eq}}(\bfif\fsub{1},\bfig\fsub{1},
\bfif\fsub{2},\bfig\fsub{2}\kern1pt)$
which expresses in $\pair{G}{\circ^G}$
the fact that the intersection of the supports
of the transpositions $\bfif\fsub{1}$ and $\bfig\fsub{1}$
is equal to the intersection of the supports
of the transpositions $\bfif\fsub{2}$ and $\bfig\fsub{2}$.
(The formula $\phi_{\srfs{Rep}}(\bfif,\bfig\kern1pt)$ is trivial.
It just says the $\bfif$ and $\bfig$ do not commute.)

McKenzie's Theorem raises two questions.
\begin{itemize}
\addtolength{\parskip}{-11pt}
\addtolength{\itemsep}{06pt}
\item[(1)]
Is the analogue of Mckenzies's Theorem for function semigroups true?
\item[(2)]
Is the analogue of Mckenzies's Theorem for clones true?
\vspace{-05.7pt}
\end{itemize}
The notion of a clone $C$ on a set $A$ will be defined later in the introduction.
$C$ is ``fully transpositional'' if it contains all transposions of $A$.
\noindent

This work answers positively questions (1) and (2).
The proof that we found is not short,
and we cannot be sure that there are no simpler proofs.
\vspace{2mm}

We now turn to the extensive work of 
M. Bodirsky, M. Pinsker, and A. Pongra\'cz. \cite{BPP}.
Bodirsky et.\ al., too, deal with the reconstruction of the
action of function semigroups and clones from their algebraic structure.
(They also consider reconstruction from the algebraic structure
when it is equipped with pointwise convergence topology
which come from their action.)
They consider several notions of reconsructiblity.
One of them is the notion of ``automatic homeomorphicity''.
Automatic homeomorphicity means
that every isomorphism between to two function semigroups
or two clones is a homeomorphism
with respect to their pointwise convergence topologies.
``Reconstruction up-to a homeomorphism'' is weaker than
``Reconstruction up-to a bijection''.

Here is one of their results
that seems related to the theorems proved in this work.
\vspace{2mm}

\noindent
{\bf Theorem } (M. Bodirsky, M. Pinsker, and A. Pongra\'cz)
Any closed subclone of $\bfO \eqdf \bigcup_{n = 1}^{\infty} \bbN^{\bbN^n}$
containing the symmetric group of $\bbN$
has automatic homeomorphicity with respect to all closed subclones of $\bfO$
(\cite{BPP} Section 5.1).

The work of Bodirsky et.\ al.\ concetrates more on the polymorphic clones
of $\aleph_0$-categorical structures.
Many questions in this area are still open.
\vspace{2mm}

\noindent
{\bf Function semigroups}
\vspace{-3mm}

\begin{definition}\label{d1.1}
\begin{rm}
(a)
Let $A$ be a nonempty set. A nonempty subset of $A^A$ closed under
composition is called a \underline{function semigroup}. The set $A$
is called the \underline{domain of $S$} and is denoted by $\bfs{Dom}(S)$.
\index{function semigroup@@Function semigroup}%
\index{N@dom1@@$\bfs{Dom}(S)$. If $S \subseteq A^A$, then $\bfs{Dom}(S) \eqdf A$}%

(b)
The structure $\pair{S}{\circ^S}$ is denoted by $\rfs{Alg}(S)$.
Note that $S$ is a set and not a structure.
{\bf We regard $\pmb{\circ^S}$ as a three-place relation},
rather than as a binary operation.
In this work, this will turn out to be more convenient.
Nevertheless, we keep writing $f \circ^S g = h$,
even though the correct thing to write is $\trpl{f}{g}{h} \stin \circ^S$.
\index{N@alg0@@$\rfs{Alg}(S) \eqdf \pair{S}{\circ^S}$. ($S$ is a function semigroup)}%

(c)
Let $S \subseteq A^A$ be a function semigroup.
We say that $S$ is a \underline{fully-}\break
\underline{transposional (FT)} semigroup,
if every transposition of $A$ belongs to $S$.
\index{fully-transpositional semigroup@@Fully-transpositional semigroup. Abbreviated by FT semigroup}%
\index{ft semigroup@@FT semigroup. A function semigroup on a set $A$ which contains all\\\indent transposions of $A$}%
\end{rm}
\end{definition}

We make the following conventions.

\noindent
{\bf Conventions } 
\begin{itemize}
\addtolength{\parskip}{-11pt}
\addtolength{\itemsep}{06pt}
\item[(1)]
$S$ always denotes a function semigroup.
\item[(2)]
$\bfs{Dom}(S)$ is always denoted by $A$. \hfill\hollowqed
\vspace{0pt}
\end{itemize}

Let $K_{\srfs{FT}}$ be the class of all FT semigroups $S$
such that $\abs{\bfs{Dom}(S)} \neq 1,2,6$,
and let
$K^{\srfs{Alg}}_{\srfs{FT}} \eqdf
\setm{\rfs{Alg}(S)}{S \stin K_{\srfs{FT}}}$.%
\index{N@kft-smgr@@$K_{\srfs{FT}}$. The class of all FT semigroups $S$ such that $\abs{\bfs{Dom}(S)} \neq 1,2,6$}%
\index{N@kft-alg@@$K^{\srfs{Alg}}_{\srfs{FT}} \eqdf \setm{\rfs{Alg}(S)}{S \stin K_{\srfs{FT}}}$}
\vspace{2mm}

Let $S \subseteq A^A$ be a function semigroup.
We define the structure $\rfs{Act}(S)$ as follows.
$$
\rfs{Act}(S)\ \ \eqdf\ \ %
\pair{S \stcup A}{\rfs{App}^S},
$$
\index{N@app0@@$\rfs{App}^S$. The application function of a function semigroup $S$}%
\index{N@act0@@$\rfs{Act}(S) \eqdf \pair{S \stcup A}{\rfs{App}^S}$. ($S$ is a function semigroup)}%
where $\fnn{\rfs{App}^S}{S \sttimes A}{A}$
is the application function.
Namely, for every $f \stin S$ and $a \stin A$,
then $\rfs{App}^S(f,a) = f(a)$.

Let
$$
K^{\srfs{Act}}_{\srfs{FT}} \eqdf
\setm{\rfs{Act}(S)}{S \stin K_{\srfs{FT}}}.
$$
\index{N@kact0@@$K^{\srfs{Act}}_{\srfs{FT}} \eqdf \setm{\rfs{Act}(S)}{S \stin K_{\srfs{FT}}}$}%

Theorem A below is our first goal.
Intuitively it says that for members $S \stin K_{\srfs{FT}}$,
$\rfs{Act}(S)$ is recoverable from $\rfs{Alg}(S)$.
\\
\noindent
{\bf Theorem A }
$K^{\srfs{Act}}_{\srfs{FT}}$ is first order strongly interpretable in
$K^{\srfs{Alg}}_{\srfs{FT}}$.
\vspace{2mm}
Theorem A restated as Corollary~\ref{c9.3-08-18}(d).

The notion of ``first order strong interpretability'' has yet to be defined
(Section 2 Definition~\ref{d2.3-08-14}(d)).
In the meantime, let us state the main conclusion of Theorem~A.
This conclusion can be understood without knowing the definition of the above
notion.

\noindent
{\bf Corollary B }
Let $S,T \stin K_{\srfs{FT}}$ and $\varphi$ be an isomorphism between
$\pair{S}{\circ^S}$ and $\pair{T}{\circ^T}$.
Then there is a bijection $\tau$ between $\bfs{Dom}(S)$ and $\bfs{Dom}(T)$
such that for every $f \stin S$, $\varphi(f) = \tau f \tau\inverse$.
(See Theorem~\ref{t2.4-08-14}.)
\vspace{2mm}

\noindent
{\bf Remark }
Whereas the fact ``$\rfs{Act}(S)$ is recoverable from $\rfs{Alg}(S)$''
is a theorem,
the converse - ``$\rfs{Alg}(S)$ is recoverable from $\rfs{Act}(S)$''
is a triviality.

Let $\calL$ denote the language of structures of the
form $\rfs{Act}(S)$, where $S$ is a function semigroup.
There are formulas $\phi_{\bm{F}}(\bfix\kern1pt)$
and $\phi_{\pmb{\circ}}(\bfif,\bfig,\bfih\kern1pt)$ in $\calL$
such that for every function semigroup $S$:

(1) For every $x \stin \abs{\rfs{Act}(S)}$:
$x \stin S$ iff $\rfs{Act}(S) \models \phi_{\bm{F}}[x]$, and

(2) For every $f,g,h \stin S$:
$f \circ g = h$ iff
$\pair{S \stcup A}{\rfs{App}^S} \models \phi_{\pmb{\circ}}[f,g,h]$.
\\
Finding $\phi_{\bm{F}}$ and $\phi_{\pmb{\circ}}$ is trivial.

Facts (1) and (2) mean that $\rfs{Alg}(S)$ is recoverable from
$\rfs{Act}(S)$.
With the terminology introduced in Definition~\ref{d2.3-08-14}(d),
it means that $\setm{\rfs{Alg}(S)}{S \mbox{ is a}\break\mbox{function semigroup}}$
is first order strongly interpretable in
{\thickmuskip=2mu \medmuskip=1mu \thinmuskip=1mu
$\setm{\rfs{Act}(S)}{S \mbox{ is}\break\mbox{a function semigroup}}$.}
Note that this interpretability fact is true for the class of all function
semigroups, and not just for fully transpositional semigroups.
\vspace{2mm}

\noindent
{\bf Clones}

In the context of this work,
reconstruction questions can be asked for permutation groups,
for function semigroups and for clones.
Indeed, the counterpart of Theorem A for clones is also true,
and is proved in this work.

We now define the notion of a clone and state the counterparts of Theorem~A
and Corollary B for clones.

Let $\bbN^+$ denote the set of positive natural numbers.
For a set $A$ and $1 \leq i \leq n \stin \bbN$
let $\fnn{p^A_{n,i}}{A^n}{A}$ be the projection onto the $i$'th coordinate.
That is,
$$
p^A_{n,i}(x_1,\ldots,x_n) = x_i.
$$
A subset $C$ of\ \kern1pt$\bigcup_{1 \leq n \stin \bbN} A^{A^n}$
is called a \underline{clone} on $A$, if
\begin{itemize}
\addtolength{\parskip}{-11pt}
\addtolength{\itemsep}{06pt}
\item[(1)]
For every $1 \leq i \leq n \stin \bbN$, \ $p^A_{n,i} \stin C$,
\item[] and
\item[(2)] For every $f(y_1,\ldots,y_n) \stin C$,
$k \stin \bbN^+$
and $g_1(\vecx\kern2pt),\ldots,g_n(\vecx\kern2pt) \stin C \stcap A^{A^k}$,
the function
$h(\vecx\kern2pt) \eqdf f(g_1(\vecx\kern2pt),\ldots,g_n(\vecx\kern2pt))$
belongs to $C$.
\vspace{-05.7pt}
\end{itemize}

Let $C$ be a clone on $A$.
Set $\rfs{Fnc}_n(C) \eqdf C \stcap A^{A^n}$.
\index{N@fnc@@$\rfs{Fnc}_n(C) \eqdf C \stcap A^{A^n}$}%
For every $n,k \stin \bbN^+$ we define an
$(n + 1)$-place partial operation $\rfs{Cmp}^C_{n,k}$.
\begin{itemize}
\addtolength{\parskip}{-11pt}
\addtolength{\itemsep}{06pt}
\item[(1)]
$\rfs{Dom}(\rfs{Cmp}^C_{n,k}) =
\rfs{Fnc}_n(C) \sttimes \left(\rule{0pt}{12pt}\rfs{Fnc}_k(C)\right)^n$,
and
\item[(2)]
For every $f \stin \rfs{Fnc}_n(C)$, $g_1,\ldots,g_n \stin \rfs{Fnc}_k(C)$
and $\vecx \stin A^k$,
$$
\left(\rule{0pt}{12pt}\rfs{Cmp}^C_{n,k}(f,g_1,\ldots,g_n)\right)(\vecx\kern1pt) =
f(g_1(\vecx\kern1pt),\ldots,g_n(\vecx\kern1pt)).
$$
\vspace{-8.0mm}
\end{itemize}
We also define the $n$-place application operation $\rfs{App}^C_n$,
which is an $(n + 1)$-place partial operation.
\begin{itemize}
\addtolength{\parskip}{-11pt}
\addtolength{\itemsep}{06pt}
\item[(1)]
$\rfs{Dom}(\rfs{App}^C_n) = \rfs{Fnc}_n(C) \times A^n$, and
\item[(2)]
For every $f \stin \rfs{Fnc}_n(C)$,
and $a_1,\ldots,a_n \stin A$,
$$
\rfs{App}^C_n(f,a_1,\ldots,a_n) = f(a_1,\ldots,a_n).
$$
\vspace{-08.0mm}
\end{itemize}
\index{N@cmp@@$\rfs{Cmp}^C_{n,k}$. The $(n,k)$'th composition function of a clone $C$}%
\index{N@app1@@$\rfs{App}^C_n$. The $n$-place application function of a clone $C$}%

Let $C$ be a clone on $A$.
We define two structures based on $C$.
The \underline{algebraic structure} of $C$ is the following structure:
\vspace{-1mm}
$$
\rfs{Alg}(C) \eqdf
\pair{C}{\kern3pt\sngltn{\rfs{Cmp}^C_{n,k})}_{n,k \stin \bbN^+}}.
\vspace{-1mm}
$$
\index{N@alg1@@$\rfs{Alg}(C) \eqdf \pair{C}{\kern3pt\sngltn{\rfs{Cmp}^C_{n,k})}_{n,k \stin \bbN^+}}$. ($C$ is a clone)}%
The \underline{action structure} of $C$ is the following structure:
\vspace{-1mm}
$$
\rfs{Act}(C) =
\pair{C \stcup A}{\kern3pt\sngltn{\rfs{App}^C_n)}_{n \stin \bbN^+}}.
\vspace{-5mm}
$$
\index{N@act1@@$\rfs{Act}(C) \eqdf \pair{C \stcup A}{\kern3pt\sngltn{\rfs{App}^C_n)}_{n \stin \bbN^+}}$. ($C$ is a clone)}%

Since partial operations do not conform with the standard definition
of a structure, we shall regard any of these $(n + 1)$-place partial operations
as an $(n + 2)$-place relation.
Nevertheless we shall continue writing them as functions
rather than as relations.
For example,
$\rfs{Cmp}^C_{n,k}(f,g_1,\ldots,g_n) = h$ stands for
$\fivetpl{f}{g_1}{\ldots}{g_n}{h} \stin \rfs{Cmp}^C_{n,k}$.

Let $C$ be a clone on $A$.
$C$ is called a \underline{fully transpostional clone (FT clone)},
\index{ft clone@@FT clone. A clone on a set $A$ which contains all transpositions of $A$}%
if every transposition of $A$ belongs to $C$.

If $C$ is a clone on $A$, then $A$ is denoted by $\bfs{Dom}(C)$.
\index{N@dom2@@$\bfs{Dom}(C)$. If $C$ is a clone on $A$, then $\bfs{Dom}(C)$ denotes $A$}%
Let $K_{\srfs{FT-cln}}$ be the class of all FT clones $C$ such that
$\abs{\bfs{Dom}(C)} \neq 1,2,6$.
\index{N@kft-cln@@$K_{\srfs{FT-cln}}$. The class of all FT clones $C$ such that $\abs{\bfs{Dom}(C)} \neq 1,2,6$}%
Set
$$\mbox{
$K^{\srfs{Alg}}_{\srfs{FT-cln}} \eqdf
\setm{\rfs{Alg}(C)}{C \stin K_{\srfs{FT-cln}}}$%
\ \ and\ \ %
$K^{\srfs{Act}}_{\srfs{FT-cln}} \eqdf
\setm{\rfs{Act}(C)}{C \stin K_{\srfs{FT-cln}}}$.
}
\hspace{-2mm}
$$
\index{N@kalg-ft-cln@@$K^{\srfs{Alg}}_{\srfs{FT-cln}} \eqdf \setm{\rfs{Alg}(C)}{C \stin K_{\srfs{FT-cln}}}$}%
\index{N@kact-ft-cln@@$K^{\srfs{Act}}_{\srfs{FT-cln}} \eqdf \setm{\rfs{Act}(C)}{C \stin K_{\srfs{FT-cln}}}$}%

\noindent
{\bf Theorem C }
$K^{\srfs{Act}}_{\srfs{FT-cln}}$
is first order strongly interpretable in $K^{\srfs{Alg}}_{\srfs{FT-cln}}$.
\vspace{2mm}

Finally, we mention the main corollary. It is the analogue of Corollary~B.%
\vspace{2mm}

\noindent
{\bf Corollary D }
Let $C,D \stin K_{\srfs{FT-cln}}$ and $\varphi$ be an isomorphism between
$\rfs{Alg}(C)$ and $\rfs{Alg}(D)$.
Then there is a bijection $\tau$ between $\bfs{Dom}(C)$ and $\bfs{Dom}(D)$
such that for every $n \stin \bbN^+$,
$f \stin \rfs{Fnc}_n(C)$ and $a_1,\ldots,a_n \stin \bfs{Dom}(C)$,
\vspace{-3mm}
$$
\varphi(f)(\tau(a_1),\ldots,\tau(a_n)) = \tau(f(a_1,\ldots,a_n)).
$$
%
%

\section{\bf Preliminaries concerning interpretations}\label{s2}

This work is to a large degree self-contained.
Though it uses the terminology of mathematical logic,
it does not rely on any real theorems from mathematical logic.
The only facts being used are facts which are usually called
``general nonsense''.

In this section we present two notions of interpretability:
``first order interpretability'' and ``first order strong interpretability''.
It is the latter notion which will be used. But the first notion is needed
in the definition of the second.
For each notion, we state a reconstruction theorem related to that notion.
We do not give the proofs of these theorems.
However, all the proofs except for \ref{t2.4-08-14} are trivial.
Most of the statements presented in this section are proved in \cite{Ru}.
The ones which are not proved there, are completely trivial.

The terms a ``language", a ``structure" or a ``model", and a ``first 
order formula" are used as in mathematical logic.
A systematic definition of these notions can be found in \cite{CK} 
Section 1.3 or in many other introductory texts in mathematical logic.
(See e.g. Herbert Enderton - A Mathematical Introduction to Logic,
or Wilfrid Hodges - Model Theory.)

Let $\calL$ be a first order language.
We use bold slanted letters $\bfia,\bfib,\bfif,\bfig,\bfix,\bfiy\kern2pt$
as the individual variables of $\calL$.
An $\calL$-formula is a first order formula in the language $\calL$.
If $M$ is a structure, then its language is denoted by $\calL(M)$.
\index{N@l@@$\calL(M)$. The language of the structure $M$}%
Let $M$ be a structure, and $\phi(\bfix\fsub{1},\ldots,\bfix\fsub{n})$
be an $\calL(M)$-formula.
Then $\phi[M^n] \eqdf\break
\setm{\vec{a} \stin |M|^n\ }{\ M \models \phi[\vec{a}\kern2pt]}$.
\index{N@AAAA@@$\phi[M^n] \eqdf \setm{\vec{a} \stin |M|^n\ }{\ M\models \phi[\vec{a}]}$}%
For simplicity,
we deal only with languages which do not contain function symbols.

\begin{definition}\label{d1.1-08-14}
\begin{rm}
Let $K$ and $K^*$ be classes of structures in the languages\break
$\calL$ and $\calL^*$ respectively and $R \subseteq K \sttimes K^*$.
$K^*$ is \underline{first order interpretable}
\underline{(FO-interpretable)} in $K$ with respect to $R$,
\index{first order interpretable@@First order interpretable}%
\index{fo interpretable@@FO-interpretable. Abbreviation of ``first order interpretable''}%
if there are $\calL$-formulas as listed in (1) and (2) below
\vspace{1mm}
\begin{itemize}
\addtolength{\parskip}{-11pt}
\addtolength{\itemsep}{06pt}
\item[(1)]
A formula
$\phi_{\srfs{U}}(\vecbfix\kern2pt)$
such that $n \eqdf \rfs{lngth}(\vecbfix\kern2pt) \neq 0$, and
\item[(2)]
For every $m \stin \bbN^+$ and an $m$-place relation symbol $P$ in $\calL^*$
a formula $\phi_P(\vecbfix\fsup{1},\ldots,\vecbfix\fsup{m})$,
such that $\rfs{lngth}(\vecbfix\Fsup{j}) = n$,
\vspace{-05.7pt}
\end{itemize}
such that for every $\pair{M}{M^*} \stin R$ there is
$\surj{\tau}{\phi_{\srfs{U}}[M^n]}{\abs{M^*}}$ such that: 
for every $m \stin \bbN^+$, every $m$-place relation symbol $P$ in $\calL^*$
(including the equality symbol) and for every
$\veca\fsup{1},\ldots,\veca\fsup{m} \stin \phi_{\srfs{U}}[M^n]$:  
$\trpl{\tau(\veca\fsup{1})}{\ldots}{\tau(\veca\fsup{m})} \stin P^{M^*}$ iff 
$M \models \phi_P[\veca\fsup{1},\ldots,\veca\fsup{m}]$.

The object 
$\vecphi \eqdf \pair{\phi_{\srfs{U}}}{\setm{\phi_P}
{P \mbox{ is a relation symbol in } \calL^*}}$
is called a\break
\underline{first order interpretation (FO-interpretation)} of 
$K^*$ in $K$ with respect to $R$. The function $\tau$ is called
an \underline{interpreting mapping} of $M^*$ in $M$
with respect to $R$ and $\vecphi$.
\index{fo-interpretation@@FO-interpretation}%
\index{interpretating mapping@@Interpretating mapping}%
\end{rm}
\end{definition}

{\bf Remark } The natural situation is that the relation $R$
in the above definition will have the following additional properties:
$\rfs{Dom}(R) = K$ and $\rfs{Rng}(R) = K^*$.
This is always the case in what follows.
But there is no need to include these additional requirements in the definition.

\begin{theorem}\label{t2.2-08-14}
Suppose that $K^*$ is FO-interpretable in $K$ with respect to $R$.
Let $M_1,M_2 \stin K$ and $M_1^*,M_2^* \stin K^*$ be such that for $i = 1,2$, 
$\pair{M_i}{M_i^*} \stin R$.
Assume that $M_1\cong M_2$. Then $M_1^*\cong M_2^*$.
\end{theorem}

\begin{definition}\label{d2.3-08-14}
\begin{rm}
(a)
If $K$ is a class of structures all in the same language $\calL$, 
then $\calL$ is denoted by $\calL(K)$.
\index{N@l@@$\calL(K)$. The language of the class of structures $K$}%

(b) Let $\phi(\vecbfix,\vecbfiy\kern1.5pt)$ be an $\calL(M)$-formula,
$\rfs{lngth}(\vecbfix\kern1.5pt) = k$
and $\rfs{lngth}(\vecbfiy\kern1.5pt) = \ell$. 
Then
$$
\phi[M^k,M^{\ell}] \eqdf 
\setm{\pair{\veca}{\vecb} \stin |M|^k \sttimes |M|^{\ell}}
{M\models\phi[\veca,\vecb]}.
$$
\index{N@AAAA@@$\phi[M^k,M^{\ell}] \eqdf \setm{\pair{\veca}{\vecb} \stin |M|^k \sttimes |M|^{\ell}}{M\models\phi[\veca,\vecb]}$}%

In the sequel $K$ and $K^*$ denote classes of structures in the 
languages $\calL$ and $\calL^*$ respectively, and $R \subseteq K \sttimes K^*$.

(c)
Let $K,K^*$ and $R$ be as above.
$\trpl{K}{K^*}{R}$ is called a \underline{subuniverse system},
if for every $\pair{M}{M^*} \stin R$, $|M|\cntd |M^*|$.
\index{subuniverse system@@Subuniverse system}%

(d) Let $\trpl{K}{K^*}{R}$ be a subuniverse system. $K^*$ is
\underline{first order strongly}
\underline{interpretable (FS-interpretable)} in $K$ relative to $R$,
\index{fs-interpretable@@FS-interpretable}%
\index{first order strongly interpretable@@First order strongly interpretable}%
if there are:
\begin{itemize}
\addtolength{\parskip}{-11pt}
\addtolength{\itemsep}{06pt}
\item[(1)]
An FO-interpretation
$\vecphi = \pair{\phi_{\srfs{U}}(\vecbfix\kern2pt)}{\ldots}$
of $K^*$ in $K$,
\item[(2)]
An $\calL(K)$-formula $\phi_{\srfs{Imap}}(\vecbfix,\bfix\kern1.5pt)$ and
\item[(3)]
An $\calL(K^*)$-formula
$\phi^*_{\srfs{Imap}}(\vecbfix,\bfix\kern1.5pt)$
\vspace{-05.7pt}
\end{itemize}
such that for every $\pair{M}{M^*} \stin R$,
$\tau \eqdf \phi^*_{\srfs{Imap}}[(M^*)^{\length(\vecx)},M^*]$
is an interpreting mapping of $M^*$ in $M$ with respect to $\vecphi$,
and
{\thickmuskip=2mu \medmuskip=1mu \thinmuskip=1mu
$\phi_{\srfs{Imap}}[M^{\length(\vecx)},M] =
\setm{\pair{\veca}{b} \stin \tau}{b \stin |M|}$.
\kern-1ptWe call $\trpl{\phi^*_{\srfs{Imap}}}{\phi_{\srfs{Imap}}}{\vecphi\kern1pt}$
a \underline{first order strong interpretation} \underline{(FS-interpretation)}
of $K^*$ in $K$ relative to $R$.
}
\index{fs-interpretation@@FS-interpretation}%
\end{rm}
\end{definition}

Let $M,N$ be structure in the same language and $f$ be a function.
The notation $\iso{f}{M}{N}$ means that $f$ is an isomorphism between
$M$ and $N$.
\index{N@AAAA@@$\iso{f}{M}{N}$ means that $f$ is an isomorphism between $M$ and $N$}%

\begin{theorem}\label{t2.4-08-14}
Suppose that $K^*$ is FS-interpretable in $K$ relative to $R$.
For $i = 1,2$ let $\pair{M_i}{M_i^*} \stin R$, and let $\iso{\rho}{M_1}{M_2}$.
\underline{Then} there is a unique $\iso{\rho^*}{M_1^*}{M_2^*}$
such that $\rho\subseteq\rho^*$.
\end{theorem}

We shall need only a weak fact about the transitivity of the relation
of FS-interpretability. It is stated in Observation~\ref{new-o2.6}.

\begin{definition}\label{new-d1.5}
\begin{rm}
(a)
Let $K$ and $K^*$ be classes of structures in the languages
$\calL$ and $\calL^*$ respectively and $R \subseteq K \times K^*$.
The triple $\trpl{K}{K^*}{R}$ is called a \underline{same-universe system},
if for every $\pair{M}{M^*} \stin R$, $\abs{M} = \abs{M^*}$.
\index{same-universe-system@@Same-universe system}%

(b)
Let $\trpl{K}{K^*}{R}$ be a same-universe system.
One says that $K^*$ is
\underline{definably interpretable} in $K$ relative to $R$,
\index{definably interpretable@@Definably interpretable. $K^*$ is definably interpretable in $K$ relative to $R$}%
if for every $n \stin \bbN^+$ and every $n$-place relation symbol $P$ in
$\calL^*$ there is an $\calL$-formula $\psi_P(x_1,\ldots,x_n)$
such that for every $\pair{M}{M^*} \stin R$,
$P^{M^*} = \psi_P[M^n]$.

(c)
Let $\trpl{K}{K^*}{R}$ be a same-universe system.
One says that $K$ and $K^*$ are
\underline{bi-definably-interpretable} relative to $\dbltn{R}{R\inverse}$,
if $K^*$ is definably interpretable in $K$ relative to $R$,
and $K$ is definably interpretable in $K^*$ relative to $R\inverse$.

Note that the relation
``$K$ and $K^*$ are bi-definably-interpretable
relative to $\dbltn{R}{R\inverse}$'' is indeed symmetric.
\index{bi-definably-interpretable@@Bi-definably interpretable. $K$ and $K^*$ are bi-definably-interpretable\\\indent relative to $\dbltn{R}{R\inverse}$}%
\end{rm}
\end{definition}

\begin{observation}\label{new-o2.6}
Let $K,K_1,K_2$ be classes of structures, 
$R \subseteq K \times K_1$ and $S \subseteq K_1 \times K_2$.
Assume that $\trpl{K}{K_1}{R}$ is a substructure system,
and that $K_1$ is FS-interpretable in $K$ relative to $R$. 
Also assume that $\trpl{K_1}{K_2}{S}$ is a same-universe system,
and that $K_1$ and $K_2$
are bi-definably-interpretable with respect to $\dbltn{S}{S\inverse}$.
\underline{Then} $K_2$ is FS-interpretable with respect to $S \circ R$.
\end{observation}

\noindent
{\bf Proof }
The proof is trivial.
\hfill\hollowqed
\vspace{2mm}

\begin{definition}\label{d2.7}
\begin{rm}
(a) Let $K$ be a class of structures in the language $\calL$ and $L \subseteq K$.
One says that $L$ is a
\underline{first order definabe (FO-definable) subclass} of $K$,
if there is an $\calL$-sentence $\phi$ such that
$L = \setm{M \stin K}{M \models \phi}$.
One also says that $\phi$ defines $L$, or that $\phi$ is a defining sentence
for $L$.

(b) Let $\calL$ be a language,
and $K$ be a class of pairs of the form $\pair{M}{a}$,
where $M$ is an $\calL$-structure and $a \stin \abs{M}$.
Let $L \subseteq K$.
One says that $L$ is a
\underline{first order definabe (FO-definable) subclass} of $K$,
if there is an $\calL$-formula $\phi(\bfix\kern1pt)$ such that
$L = \setm{\pair{M}{a} \stin K}{M \models \phi[a]}$.
One also says that $\phi$ defines $L$, or that $\phi$ is a defining formula
for $L$.
\index{fo definable@@FO-definable subclass (first order definabe subclass)}%
\index{defining formula@@defining formula}%
\index{defining sentence@@defining sentence}%

(c) Let $K$ be a class of structures all in the same language~$\calL$,
and let $\phi$ be an $\calL$-sentence.
Then $\phi(K) \eqdf \setm{M \stin K}{M \models \phi}$.
\index{N@AAAA@@$\phi(K) \eqdf \setm{M \stin K}{M \models \phi}$}%
\end{rm}
\end{definition}

\begin{observation}\label{o2.8}
Let $K$ and $K^*$ be classes of structures
in the languages $\calL$ and $\calL^*$
respectively and $R \subseteq K \sttimes K^*$.
Let $\phi_1,\ldots,\phi_k$ be $\calL$-sentences.
Suppose that:
\begin{itemize}
\addtolength{\parskip}{-11pt}
\addtolength{\itemsep}{06pt}
\item[{\rm(1)}]
$K = \bigcup_{i = 1}^k \phi_i(K)$,
\item[{\rm(2)}]
For every $i = 1,\ldots,n$, $K^*$ is FS-interpretable in $\phi_i(K_i)$
with respect to $R \restriction \phi_i(K)$
\vspace{-05.7pt}
\end{itemize}
Then $K^*$ is FS-interpretable in $K$ with respect to $R$.
\end{observation}

\noindent
{\bf Proof }
The proof is trivial.
\hfill\solidqed
\vspace{2mm}

\noindent
{\bf The relativization of a formula}

The notion of ``relativization'' will be used several times. Let us
define it here and explain its role.
Let $\calL$ be a language and $\alpha(\bfix\kern1pt)$ be an $\calL$-formula
whose only free variable is $\bfix$.
For every $\calL$-formula $\phi$ we define the $\calL$-formula $\phi^{(\alpha)}$
which will be called the \underline{relativization} of $\phi$ to $\alpha$.
\index{relativization@@Relativization of a formula to another formula}%
\index{N@AAAA@@$\phi^{(\alpha)}$. The relativization of $\phi$ to $\alpha$}%
$\phi^{(\alpha)}$ is defined by induction.
\begin{itemize}
\addtolength{\parskip}{-11pt}
\addtolength{\itemsep}{06pt}
\item[(1)]
If $\phi$ is an atomic formula, then $\phi^{(\alpha)} \eqdf \phi$.
\item[(2)]
If $\star$ is a binary connective, then
$(\psi \star \theta)^{(\alpha)} \eqdf \psi^{(\alpha)} \star \theta^{(\alpha)}$.
\item[(3)]
$(\fstneg\psi)^{(\alpha)} \eqdf \fstneg\psi^{(\alpha)}$.
\item[(4)]
$(\exists \bfiy\kern2pt \psi(\bfiy\kern1pt))^{(\alpha)} \eqdf
\exists \bfiy\kern1pt(\alpha(\bfiy\kern1pt) \stwedge
\psi^{(\alpha)}(\bfiy\kern1pt))$
and
\\
$(\forall \bfiy\kern2pt \psi(\bfiy\kern1pt))^{(\alpha)} \eqdf
\forall \bfiy\kern1pt(\alpha(\bfiy\kern1pt) \rightarrow
\psi^{(\alpha)}(\bfiy\kern1pt))$
\end{itemize}

Recall that we have decided to deal only with first order languages
which do not contain function symbols.
\begin{prop}\label{p2.7}
Let $\calL$ be a language, $\alpha(\bfix\kern1pt)$ be an $\calL$-formula
whose only free variable is $\bfix$, $M$ be an $\calL$-structure,
and $\phi(\bfiy\fsub{1},\ldots,\bfiy\fsub{n})$ be an $\calL$-formula.
Suppose that $\alpha[M] \neq \stemptyset$ and let $N$ be the substructue
of $M$ whose universe is $\alpha[M]$.
Then for every $a_1,\ldots,a_n \stin \alpha[M]$:
$M \models \phi^{(\alpha)}[a_1,\ldots,a_n]$
iff\break
$N \models \phi[a_1,\ldots,a_n]$.
\end{prop}

\noindent
{\bf Proof }
The proof is trivial.
\hfill\hollowqed
\vspace{2mm}

\section{\bf A theorem of R. McKenzie, and some-trivialities}\label{s3}
Let $\rfs{Sym}(A)$ denote the symmetric group of $A$.
For a relation $R$ (and in particular for a function $R$),
let $\rfs{Dom}(R)$ and $\rfs{Rng}(R)$ denote respectively
the domain and range of $R$.
\index{N@sym@@$\rfs{Sym}(A)$. The symmetric group of $A$}%
\index{N@dom0@@$\rfs{Dom}(R)$. The domain of a relation $R$}%
\index{N@rng@@$\rfs{Rng}(R)$. The the range of a relation $R$}%

Here are two trivial observations.

\begin{observation}\label{o1.2}
Let $S \subseteq A^A$ be an FT semigroup.

{\rm(a)}
If $A$ is finite, then either $S = \rfs{Sym}(A)$,
or for some $n \leq \abs{A}$,\break
$S = {Sym}(A) \stcup \setm{f \stin A^A}{\abs{\rfs{Rng}(f)} \leq n}$.

{\rm(b)}
If $S$ contains a function $f$ such that $\rfs{Rng}(f)$ is finite
and $\rfs{Rng}(f) \neq A$, then $S$ contains a constant function.
\end{observation}

\noindent
{\bf Proof }
The very easy proofs are left to the reader.
\hfill\hollowqed
\vspace{2mm}

\begin{definition}\label{d1.3}
\begin{rm}
Let $S \subseteq A^A$ be a function semigroup.
Set
$$
\rfs{Gr}(S) = \setm{g \stin S \stcap \rfs{Sym}(A)}{g\inverse \stin S}.
$$
(Note that $\rfs{Gr}(S)$ may be empty.)
\index{N@gr@@$\rfs{Gr}(S) \eqdf \setm{g \stin S \stcap \rfs{Sym}(A)}{g\inverse \stin S}$}%
\end{rm}
\end{definition}

Clearly, if $\rfs{Gr}(S) \neq \stemptyset$,
then $\rfs{Gr}(S)$ is a subsemigroup of $S$,
and $\pair{\rfs{Gr}(S)}{\circ}$ is a group.

Let
$$
\phi_{\srfs{Id}}(\bfif\kern2pt) \ \ \eqdf\ \ %
\forall \bfig(\bfif \bfig = \bfig \bfif = \bfig\kern1pt).
$$
and
$$
\phi_{\srfs{Gr}}(\bfif\kern2pt) \ \ \eqdf\ \ %
\exists \bfig\kern1pt
(\phi_{\srfs{Id}}(\bfif \bfig\kern1pt) \stwedge
\phi_{\srfs{Id}}(\bfig \bfif\kern2.3pt)).
\index{N@fphiid@@$\phi_{\srfs{Id}}(\bfif\kern2pt) \eqdf \forall \bfig\kern1pt(\bfif \bfig = \bfig \bfif = \bfig\kern1pt)$}%
\index{N@fphigr@@$\phi_{\srfs{Gr}}(\bfif\kern2pt) \eqdf \exists \bfig\kern1pt(\phi_{\srfs{Id}}(\bfif \bfig\kern1pt) \stwedge \phi_{\srfs{Id}}(\bfig \bfif\kern2.3pt))$}%
$$


\begin{prop}\label{p1.4}
Let $S$ be an FT semigroup. Then

{\rm (a)} For every $f \stin S$, $S \models \phi_{\srfs{Id}}[f]$ iff
$f = \rfs{Id}_A$.

{\rm (b)} $\phi_{\srfs{Gr}}[S] = \rfs{Gr}(S)$.
\end{prop}

Let $h \stin \rfs{Sym}(A)$ and $f \stin A^A$.
Then $f^h$ denotes $h f h\inverse$.
\index{N@AAAA@@$f^h \eqdf h \circ f \circ h\inverse$}%

Let $S \subseteq A^A$ be an FT semigroup.
We define the structure $\rfs{Act}^{\srfs{Gr}}(S)$ as follows.
$$
\rfs{Act}^{\srfs{Gr}}(S)\ \ \eqdf\ \ %
\trpl{S \stcup A}{\circ^S}{\rfs{App}^S_{\srfs{Gr}}},
$$
\index{N@actgr@@$\rfs{Act}^{\srfs{Gr}}(S) \eqdf \trpl{S \stcup A}{\circ^S}{\rfs{App}^S_{\srfs{Gr}}}$}%
where $\fnn{\rfs{App}^S_{\srfs{Gr}}}{\rfs{Gr}(S) \sttimes A}{A}$
is the application function.
Namely, for every $f \stin \rfs{Gr}(S)$ and $a \stin A$,
$\rfs{App}^S_{\srfs{Gr}}(f,a) = f(a)$.
\vspace{2mm}
Note that even if $\rfs{Gr}(S)$ is empty, the above definition is valid.
However, this is irrelevant since for every function semigroup $S$
appearing in this work, $\rfs{Gr}(S) \neq \stemptyset$.  

We do not treat FT semigroups $S$ for which
$\abs{\bfs{Dom}(S)} \stin \threetn{1}{2}{6}$.
This is (up to isomorphism) a finite class.
Its treatment is easy, but the semigroups in this class do not fit
conviniently into the treatment of the general case.
So including them, will be a nuisance.

Because of the same reason, we shall not treat FT clones $C$ for which
$\abs{\bfs{Dom}(C)} \stin \threetn{1}{2}{6}$.

Though not needed in the proofs of Theorems A and C, it is worthwhile noting
the fact stated in the proposition below.
This fact will be needed if one would wish to incorporate into the picture
also FT semigroups or clones $X$
for which $\abs{\bfs{Dom}(X)} \stin \threetn{1}{2}{6}$.

\begin{prop}\label{1p1.5}\label{p3.4}
{\rm(a)}
There is a sentence $\phi^{\srfs{Smgr}}_{\srfs{KFT}}$ in the language of
$K^{\srfs{Alg}}_{\srfs{FT}}$
such that for every FT semigroup $S$:
$\rfs{Alg}(S) \models \phi^{\srfs{Smgr}}_{\srfs{KFT}}$ iff $S \stin K_{\srfs{FT}}$.

{\rm(b)}
There is a sentence $\phi^{\srfs{Cln}}_{\srfs{KFT}}$ in the language of
$K^{\srfs{Alg}}_{\srfs{FT-cln}}$
such that for every FT clone $C$:
$\rfs{Alg}(C) \models \phi^{\srfs{Cln}}_{\srfs{KFT}}$ iff
$S \stin K_{\srfs{FT-cln}}$.
\end{prop}

\noindent
{\bf Proof }
(a) By Proposition~\ref{p1.4}(b), for every $n \stin \bbN$ there is a sentence
$\phi_{\srfs{group}=n}$ in the language of semigroups
such that for every FT semigroup $S$:\break
$S \models \phi_{\srfs{group}=n}$ iff
$\abs{\rfs{Gr}(S)} = n$.
Let $\phi_{\srfs{KFT}}$ be the sentence which says
$\abs{\rfs{Gr}(S)} \neq 1,2,120$.
Then by Observation~\ref{o1.2}(a), $\phi_{\srfs{KFT}}$ is as required.

(b) $\phi^{\srfs{Cln}}_{\srfs{KFT}}$ is obtained from
$\phi^{\srfs{Smgr}}_{\srfs{KFT}}$, by adjusting $\phi^{\srfs{Smgr}}_{\srfs{KFT}}$
to the language of $K^{\srfs{Alg}}_{\srfs{FT-cln}}$.
\vspace{2mm}\hfill\solidqed

Let $K$ be a class of function semigroups.
Define
$$
K^{\srfs{ActGr}} \eqdf
\setm{\rfs{Act}^{\srfs{Gr}}(S)}{S \stin K}.
$$
The language of $K^{\srfs{ActGr}}$ is denoted by
$\calL_{\srfs{ActGr}}$
\index{N@actgr@@$K^{\srfs{ActGr}} \eqdf \setm{\rfs{Act}^{\srfs{Gr}}(S)}{S \stin K}$}%
\index{N@lft-actgr@@$\calL_{\srfs{ActGr}}$}%

We reformulate the theorem of McKenzie,
in such a way that it will be convenient to apply it. 
(McKenzie's Theorem is the first theorem in the introduction.)

\begin{theorem}\label{t1.5}\label{t3.5} {\rm (R. Mckenzie)}
$K_{\srfs{FT}}^{\srfs{ActGr}}$
is first order strongly interpretable in~$K^{\srfs{Alg}}_{\srfs{FT}}$.
\end{theorem}
\vspace{2mm}

\noindent
{\bf How the proof of Theorem A is interpolated.}
\\
We shall prove the following theorem.
\vspace{2mm}

\noindent
{\bf Theorem E }
$K_{\srfs{FT}}^{\srfs{Act}}$ and $K_{\srfs{FT}}^{\srfs{ActGr}}$
are bi-definably-interpretable.
\hfill\hollowqed
\vspace{2mm}

\noindent
Then, to obtain Theorem $A$, we shall combine Theorem~\ref{t1.5}, Theorem~C and
Observation~\ref{new-o2.6}.

\noindent
{\bf Remark } Proving Theorem E amounts to proving that
$\rfs{App}^S$ is definable in $\rfs{Act}^{\srfs{Gr}}(S)$.
\vspace{2mm}

\noindent
{\bf Remark } Let $K_{\srfs{No-cnst}}$ be the class of all
$S \stin K_{\srfs{FT}}$ which do not contain a constant function.
To be precise, rather than proving Theorem E,
we shall prove the following statement.

\noindent
{\bf Theorem E$^*$ }
$K_{\srfs{No-cnst}}^{\srfs{Act}}$ and $K_{\srfs{No-cnst}}^{\srfs{ActGr}}$
are bi-definably-interpretable.
\hfill\hollowqed
\\
(Theorem E$^*$ is proved in Proposition~\ref{p9.2-15-08-17}(b).)

The class of $S$'es which contain a constant function will be dealt with
separately.

As long as we deal with the proof of Theorem E$^*$,
a ``formula'' will mean a first order formula in
$\calL_{\srfs{ActGr}}$.

\section{\bf The beginning of the proof of Theorem A:
semigroups which contain a constant function}\label{s4}

We adopt the following conventions.

\noindent
{\bf Conventions } 
\begin{itemize}
\addtolength{\parskip}{-11pt}
\addtolength{\itemsep}{06pt}
\item[(1)]
$S$ always denotes a function semigroup.
\item[(2)]
$\bfs{Dom}(S)$ is always denoted by $A$.
\item[(3)]
All the formulas that we write are first order formulas in
$\calL_{\srfs{ActGr}}$.
When we say a ``formula'',
we mean an $\calL_{\srfs{ActGr}}$-formula.
\hfill\hollowqed
\vspace{5.7pt}
\end{itemize}

\begin{notation}\label{n1.1}
\begin{rm}
(a) For a class $K$ of function semigroups,
$K^{\srfs{Alg}}$ denotes the class $\setm{\rfs{Alg}(S)}{S \stin K}$,
and $K^{\srfs{Act}} \eqdf \setm{\rfs{Act}(S)}{S \stin K}$.
(The analogous notation $K^{\srfs{ActGr}}$, has been introduced before.)
\index{N@alg1@@$K^{\srfs{Alg}} \eqdf \setm{\rfs{Alg}(S)}{S \stin K}$}%
\index{N@act1@@$K^{\srfs{Act}} \eqdf \setm{\rfs{Act}(S)}{S \stin K}$}%

(b) A pair of the form $\pair{S}{f}$, where $S$ is a function semigroup
and $f \stin S$ is called an \underline{augmented function semigroup}.
If $K$ is a class of function semigroup, then
$K^{\srfs{Aug}} \eqdf \setm{\pair{S}{f}}{S \stin K \mbox{ and } f \stin S}$.
\index{augmented@@augmented function semigroup}%
\index{N@aug@@$K^{\srfs{Aug}} \eqdf \setm{\pair{S}{f}}{S \stin K \mbox{ and } f \stin S}$}%

(c)
We extend the notation introduced in Part (a) to augmented function semigroups.
Let $K$ be a class of function semigroups and
Let $J \subseteq K^{\srfs{Aug}}$.
Then $J^{\srfs{Alg}}$ denotes the class
$\setm{\pair{\rfs{Alg}(S)}{f}}{\pair{S}{f} \stin~J}$.
The classes $J^{\srfs{Act}}$ and $J^{\srfs{ActGr}}$
are defined in the same way.
\hfill\hollowqed
\index{N@alg2@@$J^{\srfs{Alg}} \eqdf \setm{\pair{\rfs{Alg}(S)}{f}}{\pair{S}{f} \stin J}$}%
\index{N@act2@@$J^{\srfs{Act}} \eqdf \setm{\pair{\rfs{Act}(S)}{f}}{\pair{S}{f} \stin J}$}%
\index{N@actgr2@@$J^{\srfs{ActGr}} \eqdf \setm{\pair{\rfs{ActGr}(S)}{f}}{\pair{S}{f} \stin J}$}%
\end{rm}
\end{notation}

Let $K_{\srfs{$\exists$-cnst}}$ be the class of all $S \stin K_{\srfs{FT}}$
such that $S$ contains a constant function.
\index{N@kexstsscnst@@$K_{\srfs{$\exists$-cnst}}$. The class of all $S \stin K_{\srfs{FT}}$ such that $S$ contains a constant\\\indent function}%
We deal with this subclass of $K_{\srfs{FT}}$ separately.
Afterwords we shall deal with the class
$K_{\srfs{FT}} \swsetminus K_{\srfs{$\exists$-cnst}}$.
We abbreviate $(K_{\srfs{$\exists$-cnst}})^{\srfs{Alg}}$
by $K_{\srfs{$\exists$-cnst}}^{\srfs{Alg}}$ and so on.

We need to prove that $K_{\srfs{$\exists$-cnst}}^{\srfs{Act}}$
is FS-interpretable in $K^{\srfs{Alg}}_{\srfs{$\exists$-cnst}}$.
However, a stronger statement is true, and its proof is pretty trivial
and well-known.
Let $K_{\srfs{$\forall$-cnst}}$ be the class of all function semigroups
$S$ such
that $\abs{\bfs{Dom}(S)} \geq~2$,
and $S$ contains all constant functions from
$\bfs{Dom}(S)$ to $\bfs{Dom}(S)$.

\begin{prop}\label{p4.1-08-15}
{\rm(a)} $K_{\srfs{$\exists$-cnst}} \subseteq K_{\srfs{$\forall$-cnst}}$.

{\rm(b)} $K_{\srfs{$\forall$-cnst}}^{\srfs{Act}}$ is FS-interpretable
in $K_{\srfs{$\forall$-cnst}}^{\srfs{Alg}}$.

{\rm(c)} $K^{\srfs{Act}}_{\srfs{$\exists$-cnst}}$
is FS-interpretable in $K^{\srfs{Alg}}_{\srfs{$\exists$-cnst}}$.
\end{prop}

\noindent
{\bf Proof }
(a) This part is trivial.

(b) The following two facts are trivial.

{\bf Fact 1 } For a function semigroup $S$,
let $\rfs{Cnst}(S)$ denote the set of constant functions belonging to $S$.
Let $\phi_{\srfs{Cnst}}(\bfif\kern2pt) \eqdf
\forall \bfig\kern2pt(\bfif \bfig = \bfif\kern2pt)$.
Then for every\break
$S \stin K_{\srfs{$\forall$-cnst}}$,
$\phi_{\srfs{Cnst}}[S] = \rfs{Cnst}(S)$.

{\bf Fact 2 }
For a set $B$ and $b \stin B$,
let $\rfs{cnst}_b^B$ denote the constant function from $B$ to $B$
whose value is $b$.
Let $S \stin K_{\srfs{$\forall$-cnst}}$.
Then for every $f \stin S$ and $a,b \stin A$,
$f(a) = b$ iff $f \circ \rfs{cnst}_a^A = \rfs{cnst}_b^A$.

The FS-interpretability of $K^{\srfs{Act}}_{\srfs{$\forall$-cnst}}$
in $K^{\srfs{Alg}}_{\srfs{$\forall$-cnst}}$
follows easily from Fact~1 and Fact~2.
The exact details are left to the reader.

(c) Part (c) is a trivial corollary of Parts (a) and (b).
\vspace{2mm}\hfill\solidqed

We shall also use the trivial fact that $K^{\srfs{Alg}}_{\srfs{$\exists$-cnst}}$
is FO-definable in $K^{\srfs{Alg}}_{\srfs{FT}}$.

\begin{observation}\label{o4.2}
Let
$\phi_{\srfs{$\exists$-cnst}} \eqdf
\exists \bfif\kern2.9pt \forall \bfig\kern1pt(\bfif \bfig = \bfif\kern2.3pt)$.
Then for every $S \stin  K_{\srfs{FT}}$,
$\rfs{Alg}(S) \models \phi_{\srfs{$\exists$-cnst}}$ iff
$S \stin K_{\srfs{$\exists$-cnst}}$.
\end{observation}
\index{N@fphiexists-cnst@@$\phi_{\srfs{$\exists$-cnst}} \eqdf \exists \bfif\kern2.9pt \forall \bfig\kern1pt(\bfif \bfig = \bfif\kern2.3pt)$}%

\noindent
{\bf Proof } The proof is trivial \hfill\solidqed
\vspace{2mm}

\section{\bf \bf An outline of the continuation of the proof}\label{s5}

Let $\fnn{f}{A}{A}$. Define
\vspace{-3mm}
$$
\rfs{Fxd}(f) \eqdf \setm{a \stin A}{f(a) = a}
\vspace{-3mm}
$$
and
\vspace{-3mm}
$$
\rfs{Fxd-img}(f) \eqdf \setm{a \stin A}{f(a) \stin \rfs{Fxd}(f)}.
$$
\index{N@fxd@@$\rfs{Fxd}(f) \eqdf \setm{a \stin \rfs{Dom}(f)}{f(a) = a}$}%
\index{N@fxdimg@@$\rfs{Fxd-img}(f) \eqdf \setm{a \stin A}{f(a) \stin \rfs{Fxd}(f)}$}%
We shall prove the following lemma.

\noindent
{\bf Lemma F }
There is a formula $\phi_{\srfs{Fxd-img}}(\bfif,\bfix\kern1.5pt)$
such that for every $S \stin K_{\srfs{FT}}$ which does not contain a constant
function, $f \stin S$ and $a \stin A$:
\vspace{-2mm}
\reqnomode
\begin{align}\mbox{
$a \stin \rfs{Fxd-img}(f)$%
\ \ iff\ \ %
$\rfs{Act}^{\srfs{Gr}}(S) \models \phi_{\srfs{Fxd-img}}[f,a]$.
\tag*{\hollowqed}
}\end{align}
\index{N@fphifxd-img@@$\phi_{\srfs{Fxd-img}}(\bfif,\bfix\kern1.5pt)$}%
\vspace{-5mm}

Once Lemma F is proved, it is a very short way to proving Theorem~A.
\vspace{2mm}

Let $K_{\srfs{No-cnst}}$ be the class of all semigroups $S \stin K_{\srfs{FT}}$
such that $S$ does not contain a constant function.
\index{N@knocnst@@$K_{\srfs{No-cnst}} \eqdf K_{\srfs{FT}} \swsetminus K_{\srfs{$\exists$-cnst}}$}%
We shall partition $K_{\srfs{No-cnst}}^{\srfs{Aug}}$
into four subclasses $J_1,J_2,J_3,J_4$.
For every $i = 1,\ldots,4$: (1) we shall show that
$(J_i)^{\srfs{ActGr}}$
is FO-definable in $(K^{\srfs{Aug}}_{\srfs{No-cnst}})^{\srfs{ActGr}}$.
(2) we shall find a formula
$\phi^{J_i}_{\srfs{Fxd-img}}(\bfif,\bfix\kern1.5pt)$
such that for every $\pair{S}{f} \stin J_i$ and $a \stin A$:
$f(a) \stin \rfs{Fxd}(f)$\ \ iff\ \ %
$\rfs{Act}^{\srfs{Gr}}(S) \models \phi^{J_i}_{\srfs{Fxd-img}}[f,a]$.

The formula $\phi_{\srfs{Fxd-img}}$ of Lemma F is obtained
from the $\phi^{J_i}_{\srfs{Fxd-img}}$'s
and the defining formulas of the $J_i$'s.
We desrcibe how this is done.
Let $\phi_{J_i}(\bfif\kern2.3pt)$ be the formula defining $(J_i)^{\srfs{ActGr}}$
in $(K^{\srfs{Aug}}_{\srfs{No-cnst}})^{\srfs{ActGr}}$.
Then
\vspace{-2mm}
$$\mbox{
$\phi_{\srfs{Fxd-img}}(\bfif,\bfix\kern1pt) \eqdf
\bigvee_{i = 1}^4
\left(\rule{0pt}{12pt}
\phi_{J_i}(\bfif\kern2.3pt) \kern1pt\stwedge\kern1pt
\phi^{J_i}_{\srfs{Fxd-img}}(\bfif,\bfix\kern1pt)\right)$.
}
\vspace{-2mm}
$$
It is a triviality that, indeed, $\phi_{\srfs{Fxd-img}}$ is as needed.
\vspace{2mm}

We turn to the definition of the $J_i$'s,
but we start with some preliminary notions.

\begin{definition}\label{d4.1.08-14}\label{d5.2}
\begin{rm}
(a)
Let $\fnn{f}{A}{A}$. The function $f$ is called a \underline{semi-constant}
\index{semiconstan@@Semi-constant function}%
function, if there is $B \subseteq A$ such that $\abs{B} \geq 2$,
$f \restriction B$ is a constant function
and $f \restriction (A \swsetminus B) = \rfs{Id}_{A \setminus B}$.

(b)
Let $K_{\srfs{$\exists$-scnst}}$ be the class of all $S \stin K_{\srfs{FT}}$
such that $S$ contains\break
a semi-constant function,
but does not contain constant functions.
\index{N@kexstsscnst@@$K_{\srfs{$\exists$-scnst}}$.
The class of all $S \stin K_{\srfs{FT}} \swsetminus K_{\srfs{$\exists$-cnst}}$ such that $S$ contains a\\\indent semi-constant function}%
Let
$K_{\srfs{No-scnst}} \eqdf
K_{\srfs{FT}} \swsetminus (K_{\srfs{$\exists$-scnst}} \stcup
K_{\srfs{$\exists$-cnst}})$.
\index{N@knoscnst@@$K_{\srfs{No-scnst}} \eqdf K_{\srfs{FT}} \swsetminus (K_{\srfs{$\exists$-scnst}} \stcup K_{\srfs{$\exists$-cnst}})$}%

(c) Let $\fnn{f}{A}{A}$. Define the equivalence relation $\sim_f$ on $A$
as follows: $a \sim_f b$, if $a,b \stin A$ and $f(a) = f(b)$.
\index{N@AAAA@@$\sim_f$. $a \sim_f b$ if $f(a) = f(b)$}%
We denote $a \kern1.5pt\slash\kern-2.5pt \sim_f$ by $[a]_f$.
\index{N@AAAA@@$[a]_f \eqdf a \kern1.5pt\slash\kern-2.5pt \sim_f$}%

(d) For $\fnn{f}{A}{A}$ define
$\bfB(f) \eqdf \setm{[a]_f}{a \stin A \mbox{ and } \abs{[a_f]} \geq 2}$.
\index{N@b@@$\bfB(f) \eqdf \setm{[a]_f}{a \stin A \mbox{ and } \abs{[a_f]} \geq 2}$}%

(e)
Let $S \stin K_{\srfs{FT}}$. Set
$\rfs{F}_{\srfs{\sbf{B}$\geq$3}}(S) \eqdf \setm{f \stin S}{\abs{\bfB(f)} \geq 3}$
\index{N@fb1@@$\rfs{F}\Nsub{\srfs{\sbf{B}$\geq$3}}(S) \eqdf \setm{f \stin S}{\abs{\bfB(f)} \geq 3}$}%
and\break
$\rfs{F}_{\srfs{\sbf{B}$\leq$2}}(S) \eqdf
\setm{f \stin S}{\abs{\bfB(f)} \leq 2}$.
\index{N@fb2@@$\rfs{F}\Nsub{\srfs{\sbf{B}$\leq$2}}(S) \eqdf \setm{f \stin S}{\abs{\bfB(f)} \leq 2}$}%
\end{rm}
\end{definition}

\index{N@j1@@$J_1 \eqdf (K_{\srfs{$\exists$-scnst}})^{\srfs{Aug}}$}

We are ready to define $J_1,\ldots,J_4$.

\begin{definition}\label{d5.3}
\begin{rm}
$J_1 \eqdf (K_{\srfs{$\exists$-scnst}})^{\srfs{Aug}}$,
\begin{itemize}
\addtolength{\parskip}{-11pt}
\addtolength{\itemsep}{06pt}
\item[]
$J_2 \eqdf \setm{\pair{S}{f}}
{S \stin K_{\srfs{No-scnst}} \mbox{ and }
f \stin \rfs{F}_{\srfs{\sbf{B}$\geq$3}}(S)}$,
\index{N@j2@@$J_2$}%
\item[]
$J_3 \eqdf \setm{\pair{S}{f}}
{S \stin K_{\srfs{No-scnst}} \text{ and }
\ f \stin \rfs{F}_{\srfs{\sbf{B}$\leq$2}}(S) \swsetminus \rfs{Gr}(S)}$
\vspace{-05.7pt}
\end{itemize}
and
\begin{itemize}
\addtolength{\parskip}{-11pt}
\addtolength{\itemsep}{06pt}
\item[]
$J_4 \eqdf \setm{\pair{S}{f}}{S \stin K_{\srfs{No-scnst}} \text{ and }
f \stin \rfs{Gr}(S)}$.
\end{itemize}
\end{rm}
\end{definition}
\index{N@j3@@$J_3$}%
\index{N@j4@@$J_4$}%
\vspace{-4mm}

\begin{observation}\label{o5.3-08-17}\label{o5.4}
$\fourtn{J_1}{J_2}{J_3}{J_4}$ is a partition of
$(K_{\srfs{No-cnst}})^{\srfs{Aug}}$.
\end{observation}

\noindent
{\bf Proof } The proof is trivial. (Indeed, there is nothing to prove.)
\hfill\solidqed
\vspace{2mm}

\section{\bf The definability of the relation
``$\boldsymbol{f(b) \stin \bfs{Fxd}(f)}$''
in FT semigroups which contain a semi-constant function}\label{s6}

According to the plan of the proof laid out in the preceeding section,
we have to show that
(i) $(J_1)^{\srfs{ActGr}}$
is FO-definable in $((K_{\srfs{No-cnst}})^{\srfs{Aug}})^{\srfs{ActGr}}$,
and that (ii) there is a formula
$\phi^{J_1}_{\srfs{Fxd-img}}(\bfif,\bfix\kern1pt)$,
which for members $\pair{S}{f}$ of $J_1$
expresses in $\rfs{Act}^{\srfs{Gr}}(S)$ the fact that
$x \stin \rfs{Fxd-img}(f)$.
More precisely, we need an $\calL_{\srfs{ActGr}}$-formula
$\phi^{J_1}_{\srfs{Fxd-img}}(\bfif,\bfix\kern1.5pt)$
such that for every $\pair{S}{f} \stin J_1$ and $x \stin \bfs{Dom}(S)$:
$\rfs{Act}^{\srfs{Gr}}(S) \models \phi^{J_1}_{\srfs{Fxd-img}}[f,x]$
iff $x \stin \rfs{Fxd-img}(f)$.

Fact (i) is trivial.
It will be proved in Proposition~\ref{p5.12-08-15}.
Fact (ii) - the expressibility of ``$x \stin \rfs{Fxd-img}(f)$'',
will be proved in Lemma~\ref{l5.11-08-16},
which is the last lemma in this section.

\begin{remark}\label{r1.1}
\begin{rm}
Let $\calL$ be a language and for $i = 1,2,3$
let $K_i$ be a class of pairs of the form $\pair{M}{a}$,
where $M$ is an $\calL$-structure and $a \stin \abs{M}$.
Suppose that $K_1 \subseteq K_2 \subseteq K_3$ and that $\phi(\bfia\kern1pt)$
defines $K_1$ in $K_3$.
Then, of course, $\phi(\bfia\kern1pt)$ defines $K_1$ in $K_2$.
So instead of proving Fact (i) we shall prove that
$\setm{\pair{\rfs{Act}^{\srfs{Gr}}(S)}{f}}{\pair{S}{f} \stin J_1}$
is FO-definable in
$\setm{\pair{\rfs{Act}^{\srfs{Gr}}(S)}{f}}{S \stin K_{\srfs{FT}}
\mbox{ and } f \stin S}$.
\end{rm}
\end{remark}

\begin{definition}\label{d3.3-08-14}
\begin{rm}
Let $A$ be a set and $a,b \stin A$ be distinct. The transposition of $A$ which
switches $a$ and $b$ is denoted by $\rvpair{a}{b}^A$.
\index{N@AAAA@@$\rvpair{a}{b}^A$. The transposition of $A$ which switches $a$ and $b$.\\\indent(Abbreviated by $\rvpair{a}{b}$)}%
Since $A$ will be fixed throughout the context, we abbreviate
$\rvpair{a}{b}^A$ by $\rvpair{a}{b}$.
\end{rm}
\end{definition}

\begin{observation}\label{o1.6}\label{o4.4.0814}
There is an an $\calL_{\srfs{ActGr}}$-formula
$\phi_{\srfs{tr}}(\bfif,\bfix,\bfiy\kern1pt)$ 
which for semigroups $S \stin K_{\srfs{FT}}$
says that that $\bfix \neq \bfiy$ and that
$\bfif = \rvpair{\bfix}{\bfiy\kern1pt}$. 
\end{observation}

\noindent
{\bf Proof }
Let
$$
\phi^{\srfs{ActGr}}_{\srfs{Dom}}(\bfix\kern1pt) \eqdf
\exists \bfif\kern2.3pt \exists \bfiy\kern1pt
(\rfs{Act}(\bfif,\bfix\kern1pt) = \bfiy\kern1pt).
$$
Then for every $S \stin K_{\srfs{FT}}$,
$$
\phi^{\srfs{ActGr}}_{\srfs{Dom}}[\rfs{Act}^{\srfs{Gr}}(S)] = \bfs{Dom}(S).
$$
Let $\phi_{\srfs{tr}}(\bfif,\bfix,\bfiy\kern1pt)$
be the $\calL_{\srfs{ActGr}}$-formula which says:
\begin{itemize}
\addtolength{\parskip}{-11pt}
\addtolength{\itemsep}{06pt}
\item[(1)]
$\bfix \neq \bfiy$,
\item[(2)]
$\rfs{Act}(\bfif,\bfix\kern1pt) = \bfiy$
and
$\rfs{Act}(\bfif,\bfiy\kern1pt) = \bfix$, and
\item[(3)]
For every $\bfiz$\kern1pt: if $\phi^{\srfs{ActGr}}_{\srfs{Dom}}(\bfiz\kern1pt)$
and  $\bfiz \neq \bfix,\bfiy$,
then $\rfs{Act}(\bfif,\bfiz\kern1pt) = \bfiz$.
\vspace{-05.7pt}
\end{itemize}
Then $\phi_{\srfs{tr}}(\bfif,\bfix,\bfiy\kern1pt)$
is as required
\rule{20mm}{0pt}\hfill\solidqed

\begin{observation}\label{o1.7}
Let $\phi_{\srfs{sim}}(\bfif,\bfix,\bfiy\kern1pt)$
\index{N@fphisim@@$\phi_{\srfs{sim}}(\bfif,\bfix,\bfiy\kern1pt)$}%
be the formula which says: Either\break
{\rm(i)} $\bfix = \bfiy$,
or {\rm(ii)} $\bfix \neq \bfiy$ and
$\bfif = f \rvpair{\bfix}{\bfiy\kern1pt}$. 
\\
\underline{Then} for every $S \stin K_{\srfs{FT}}$,
$f \stin S$ and $a,b \stin A$:
$\rfs{Act}^{\srfs{Gr}}(S) \models \phi_{\srfs{sim}}[f,a,b]$
iff $a \sim_f b$.
\end{observation}

\noindent
{\bf Proof } The proof is trivial.
\hfill\solidqed
\vspace{2mm}

A function $\fnn{f}{A}{A}$ is called a projection if
$f \restriction \rfs{Rng}(f) = \rfs{Id}_{\srfs{Rng}(f)}$.

\begin{observation}\label{o1.8}
{\rm(a)} Let $S \subseteq A^A$ be an FT semigroup and $f \stin S$.
Then $f$ is a constant function iff for every $g \stin S$, $fg = f$.
So there is an $\calL_{\srfs{ActGr}}$-formula
$\phi_{\srfs{Cnst}}(\bfif\kern2pt)$
such that for every FT semigroup $S$ and $f \stin S$,
$\rfs{Act}^{\srfs{Gr}}(S) \models \phi_{\srfs{Cnst}}[f]$
iff $f$ is a constant function.
\index{N@fphicnst@@$\phi_{\srfs{Cnst}}(\bfif\kern2pt)$. An $\calL^{\srfs{ActGr}}_{\srfs{FT}}$-formula which says that $\bfif$ is a constant function}%

{\rm(b)} Let $\fnn{f}{A}{A}$.
Then $f$ is a projection iff $f^2 = f$.
So the formula $\bfif\fsup{2} = \bfif$ defines the property
of being a projection in all function semigroups. 
So there is a $\calL_{\srfs{ActGr}}$-formula $\phi_{\srfs{Prj}}(\bfif\kern2pt)$
such that for every FT semigroup $S$ and $f \stin S$,
$\rfs{Act}^{\srfs{Gr}}(S) \models \phi_{\srfs{Prj}}[f]$
iff $f$ is a projection.
\index{N@fphiprj@@$\phi_{\srfs{Prj}}(\bfif\kern2pt)$. An $\calL^{\srfs{ActGr}}_{\srfs{FT}}$-formula which says that $\bfif$ is a projection}%

{\rm(c)} Let $\fnn{f}{A}{A}$.
Then $f$ is a semi-constant function iff
{\rm(i)} $f$ is a projection,
and \rm{(ii)} There is $a \stin A$ such that $\abs{[a]_f} \geq 2$,
and for every $b \stin A \swsetminus [a]_f$, $[b]_f = \sngltn{b}$.

So there is an $\calL^{\srfs{ActGr}}_{\srfs{FT}}$-formula
$\phi_{\srfs{scnst}}(\bfif\kern2.3pt)$
such that for every $S \stin K_{\srfs{FT}}$ and $f \stin S$,
$\rfs{Act}^{\srfs{Gr}}(S) \models \phi_{\srfs{scnst}}[f]$
iff $f$ is semi-constant.
\index{N@fphiscnst@@$\phi_{\srfs{scnst}}(\bfif\kern2pt)$. An $\calL^{\srfs{ActGr}}_{\srfs{FT}}$-formula saying that $\bfif$ is a semi-constant function}%


{\rm(d)} There is an $\calL^{\srfs{ActGr}}_{\srfs{FT}}$-sentence
$\phi_{\exists\srfs{-scnst}}$ such that for every $S \stin K_{\srfs{FT}}$,
$\rfs{Act}^{\srfs{Gr}}(S) \models \phi_{\srfs{$\exists$-scnst}}$
iff $S \stin K_{\exists\srfs{-scnst}}$.
\index{N@fphiexists-scnst@@$\phi_{\srfs{$\exists$-scnst}}$. A sentence which says that $S \stin K_{\srfs{$\exists$-scnst}}$}%
\end{observation}

\noindent
{\bf Proof } The proof is trivial.
\vspace{2mm}
\hfill\solidqed

\begin{prop}\label{p5.12-08-15}
{\rm(a)}
There is a sentence $\phi_{\srfs{$\exists$-cnst}}$
such that for every\break
$S \stin K_{\srfs{FT}}$:
$\rfs{Act}^{\srfs{Gr}}(S) \models \phi_{\srfs{$\exists$-cnst}}$ iff
$S$ contains a constant function.

{\rm(b)}
There is a sentence $\phi_{\srfs{$\exists$-scnst}}$
such that for every $S \stin K_{\srfs{FT}}$:
$\rfs{Act}^{\srfs{Gr}}(S) \models \phi_{\srfs{$\exists$-scnst}}$
iff $S \stin K_{\srfs{$\exists$-scnst}}$.
Hence
$\setm{\pair{\rfs{Act}^{\srfs{Gr}}(S)}{f}}{\pair{S}{f} \stin J_1}$
is FO-definable in
$\setm{\pair{\rfs{Act}^{\srfs{Gr}}(S)}{f}}{S \stin K_{\srfs{FT}}
\mbox{ and } f \stin S}$.
\end{prop}

\noindent
{\bf Proof }
This is a trivial corollary from Observation~\ref{o1.8}(a) and (d).
\hfill\solidqed
\vspace{2mm}

Let $f \stin A^A$ the set $\rfs{Idp}(f)$ is defined as follows.
$$
\rfs{Idp}(f) = \setm{a \stin A}{f\inverse[\sngltn{a}] = \sngltn{a}}.
$$
The members of $\rfs{Idp}(f)$ are called \underline{identity points} of $f$.
\index{N@idp@@$\rfs{Idp}(f) \eqdf \setm{a \stin A}{f\inverse[\sngltn{a}] = \sngltn{a}}$}%

Let $\fnn{f}{A}{A}$ be a semi-constant function.
Then the unique $a \stin A$ such that $\abs{f\inverse[\sngltn{a}]} \geq 2$
is called the \underline{constant value} of $f$,
and is denoted by $\rfs{cnst}(f)$,
\index{N@cnst0@@$\rfs{cnst}(f)$. The constant value of a semi-constant function $f$}%
and the set $f\inverse[\sngltn{a}]$ is called
the \underline{constant domain} of $f$, and is denoted by $\rfs{Cnst-dom}(f)$.
Note that $\rfs{cnst}(f) \stin \rfs{Cnst-dom}(f)$
and that $\rfs{Idp}(f) = A \swsetminus \rfs{Cnst-dom}(f)$.
\index{N@cnst1@@$\rfs{Cnst-dom}(f)$. The constant domain of a semi-constant function $f$}%

We shall use the following trivial fact.

\begin{observation}\label{o5.7-08-15}
Let $S$ be an FT semigroup which contains a semi-constant function which is not
constant.
Let $a,d,p \stin A$ be such that $a,d \neq p$.\break
{\thickmuskip=8mu \medmuskip=7mu \thinmuskip=6mu
Then there is a semi-constant function $g \stin S$
such that $a = \rfs{cnst}(g)$,}\break
$d \stin \rfs{Cnst-dom}(g)$ and $p \stnot\stin \rfs{Cnst-dom}(g)$.
\end{observation}

\noindent
{\bf Proof }
The proof is trivial, but we verify it nevertheless.

Let $h \stin S$ be a semi-constant function which is not constant.
Let\break
$a' = \rfs{cnst}(h)$. 
If $a' = a$, then let $h_1 = h$.
If $a' \neq a$, then let $\pi = \rvpair{a}{a'}$ and $h_1 = h^{\pi}$.
Then $h_1$ is semi-constant and not constant, and $\rfs{cnst}(h_1) = a$.

If $d = a$. Then $d \stin \rfs{Cnst-dom}(h_1)$. In this case, set $h_2 = h_1$.

Assume that $d \neq a$.
If $d \stin \rfs{Cnst-dom}(h_1)$, then let $h_2 = h_1$.
If\break
$d \stnot\stin \rfs{Cnst-dom}(h_1)$,
then let $d' \stin \rfs{Cnst-dom}(h_1) \swsetminus \sngltn{a}$,
$\sigma = \rvpair{d}{d'}$ and $h_2 = h_1^{\sigma}$.
Then $h_2$ is semi-constant and not constant, $\rfs{cnst}(h_1) = a$
and $d \stin \rfs{Cnst-dom}(h_2)$.

If $p \stnot\stin \rfs{Cnst-dom}(h_2)$, then let $g = h_2$.
If $p \stin \rfs{Cnst-dom}(h_2)$,
then let $p' \stin A \swsetminus \rfs{Cnst-dom}(h_2)$, $\eta = \rvpair{p}{p'}$
and $g = h_2^{\eta}$.
Then $g$ is semi-constant and not constant, $\rfs{cnst}(h_1) = a$,
$d \stin \rfs{Cnst-dom}(h_2)$ and $p \stnot\stin \rfs{Cnst-dom}(g)$.
\hfill\solidqed
\vspace{2mm}

\begin{prop}\label{p1.10}\label{p5.8-08-15}
Let $S \stin K_{\srfs{FT}}$, $f \stin S$ be a semi-constant function
which is not a constant function and $b \stin A$.
Then the following are equivalent:
\begin{itemize}
\addtolength{\parskip}{-11pt}
\addtolength{\itemsep}{06pt}
\item[{\rm(1)}]
$b = \rfs{cnst}(f)$
\item[{\rm(2)}]
{\rm(i)} $\abs{[b]_f} \geq 2$ and
\\
{\rm(ii)} For every $a \stin A \swsetminus [b]_f$, 
$(f \rvpair{a}{b})^3 = f \rvpair{a}{b}$. 
\vspace{-05.7pt}
\end{itemize}
\end{prop}

\noindent
{\bf Proof }
$\stneg$(1) $\Rightarrow$ $\stneg$(2):
Suppose that $b \stin A$ and $b \neq \rfs{cnst}(f)$.
If $\abs{[b]_f} = 1$. Then $b$ does not fulfill 2(i).
So assume that $\abs{[b]_f} \geq 2$.
Set $c = f(b)$. Recall that we assumed that $f$ is not a constant function.
So let $a \stin \rfs{Idp}(f)$, and let $\sigma = \rvpair{a}{b}$.
\smallskip 
\\ 
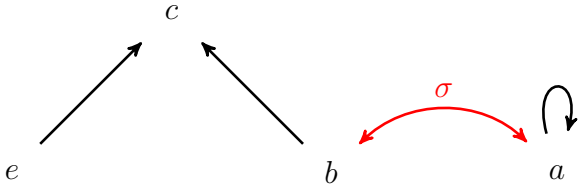
\begin{tikzpicture}[->,>=stealth',shorten >=1pt,auto,node distance=3.0cm,on grid,semithick,
                    every state/.style={fill=white,draw=none,circle,text=black}]
  \node [state] (A)                                                                               {$c$};
  \node [state] (B)  [below left= of A]                                                    {$e$};
  \node [state] (C)  [below right= of A]                                                           {$b$};
  \node [state] (D)  [right= of C]                                                    {$a$};
\path   (C) [<->, bend left=45,red,line width=1pt] edge node {$\sigma$}  (D) ;    
\path (B) [line width=1pt] edge node {} (A);
\path (C) [line width=1pt] edge node {} (A);
\path  (D) [line width=1pt]  edge [loop above]        node [swap]{}(\Cnr);
\end{tikzpicture}
\\
Then $f \sigma(b) = a$, $f \sigma(a) = c$ and $f \sigma(c) = c$.
\smallskip 
\\
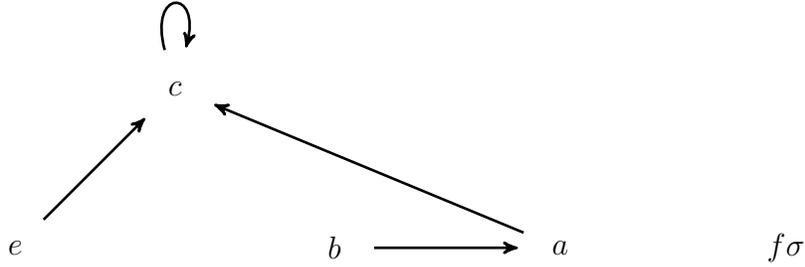
\centerline{\begin{tikzpicture}[->,>=stealth',shorten >=1pt,auto,node distance=3.0cm,on grid,semithick,
                    every state/.style={fill=white,draw=none,circle ,text=black}]
  \node [state] (A)                                                                               {$c$};
  \node [state] (B)  [below left= of A]                                                    {$e$};
  \node [state] (C)  [below right= of A]                                                           {$b$};
  \node [state] (D)  [right= of C]                                                   
{$a$};
 \node [state] (E)  [right= of D]                                                   
{$f\sigma$};
\path (B) [line width=1pt] edge node {} (A);
\path (C) [line width=1pt] edge node {} (D);
\path (D) [line width=1pt] edge node {} (A);
\path  (A) [line width=1pt]  edge [loop above]        node [swap]{}(\Cnr);
\end{tikzpicture}}
\\
It follows that $(f \sigma)^3(b) = c$ and $(f \sigma)(b) = a$. Hence
$$
(f \sigma)^3 \neq f \sigma.
$$

(1) $\Rightarrow$ (2):
Suppose that $b = \rfs{cnst}(f)$.
Then $\abs{[b]_f} \geq 2$. That is, (2)(i) holds.

Let $a \stin \rfs{Idp}(f)$ and set $\sigma = \rvpair{a}{b}$.
Then $f \sigma(b) = a$ and $f \sigma(a) = b$.
Also, for every $d \stin \rfs{Idp}(f) \swsetminus \sngltn{a}$, $f \sigma(d) = d$,
and for every $e \stin [b]_f \swsetminus \sngltn{b}$, $f \sigma(e) = b$.
\smallskip 
\\ 
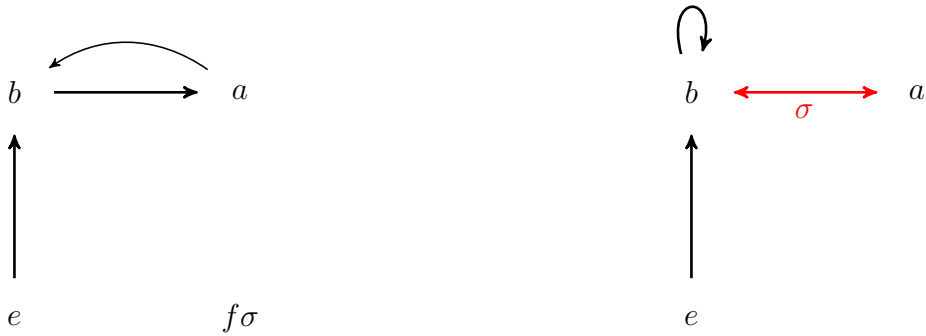
\centerline{\begin{tikzpicture}[->,>=stealth',shorten >=1pt,auto,node distance=3.0cm,on grid,semithick,
                    every state/.style={fill=white,draw=none,circle ,text=black}]
 \node [state] (O)                       {};                   
  \node [state] (A)[left= of O]         {$a$};
  \node [state] (B)  [left = of A]      {$b$};
  \node [state] (C)  [below = of B]     {$e$};
  \node [state] (G)  [right = of C]     {$f\sigma$};
  \node [state] (D)  [right= of O]      {$b$};
  \node [state] (E)  [right= of D]      {$a$};
  \node [state] (F)  [below = of D]     {$e$};
\path (C) [line width=1pt] edge node {} (B);
\path (B) [line width=1pt] edge node {} (A);
\path (A)  edge [bend right=35]  node {} (B);
\path (F) [line width=1pt] edge node {} (D);
\path  (D) [line width=1pt]  edge [loop above]        node [swap]{}(\Cnr);
\path   (E) [<->, red,line width=1pt] edge node {$\sigma$}  (D) ;  
\end{tikzpicture}}
\\
It follows that
\begin{list}{}
{\setlength{\leftmargin}{39pt}
\setlength{\labelsep}{05pt}
\setlength{\labelwidth}{25pt}
\setlength{\itemindent}{-00pt}
\addtolength{\topsep}{-04pt}
\addtolength{\parskip}{-02pt}
\addtolength{\itemsep}{-05pt}
}
\item[(i)] $(f \sigma)^3(b) = a = f \sigma(b)$
\ \ and\ \ %
$(f \sigma)^3(a) = b = f \sigma(a)$,
\item[(ii)] for every $d \stin \rfs{Idp}(f) \swsetminus \sngltn{a}$,
\ $(f \sigma)(d) = d = (f \sigma)^3(d)$,
and
\item[(iii)] for every $e \stin [b]_f \swsetminus \sngltn{b}$,
\ $(f \sigma)(e) = e = (f \sigma)^3(e)$.
\vspace{-02.0pt}
\end{list}
Hence
\reqnomode
\begin{align}
(f \sigma)^3 = f \sigma.
\tag*{\solidqed}
\end{align}
\vspace{2mm}

Our next goal is to find an $\calL^{\srfs{ActGr}}_{\srfs{FT}}$-formula
$\phi_{\srfs{Fxd-img}}^{\srfs{$\exists$-scnst}}(\bfif,\bfib\kern1pt)$
which in semigroups belonging to $K_{\srfs{FT}}$,
containing a semi-constant function but not containing a constant function,
will say that
$\bfif\kern1.5pt(\bfib\kern1pt)$ is a fixed point of $\bfif$.
Let $\phi_{\srfs{Fxd-img}}^{\srfs{$\exists$-scnst}}(\bfif,\bfib\kern1pt)$
be the formula which says:
\index{N@fphifxd-img-exists-scnst@@$\phi_{\srfs{Fxd-img}}^{\srfs{$\exists$-scnst}}(\bfif,\bfib\kern1pt)$}%
\begin{itemize}
\addtolength{\parskip}{-11pt}
\addtolength{\itemsep}{06pt}
\item
\rule{2mm}{0pt}
There is $\bfia \stin [\bfib]_{\sbfif}$ such that for every semi-constant
function $\bfig$:
\begin{list}{}
{\setlength{\leftmargin}{39pt}
\setlength{\labelsep}{05pt}
\setlength{\labelwidth}{25pt}
\setlength{\itemindent}{-00pt}
\addtolength{\topsep}{-04pt}
\addtolength{\parskip}{-02pt}
\addtolength{\itemsep}{-05pt}
}
\item[] \underline{if}
(i) $\bfia \stnot\stin \rfs{Cnst-dom}(\bfig\kern1pt)$
and (ii) $\rfs{cnst}(\bfig\kern1pt) \stnot\stin [\bfib]_{\sbfif}$,
\item[] \underline{then}
$[\bfib]_{\sbfif \bfig \bfif\kern2pt} \subseteq [\bfib]_{\sbfif\fsup{2}}$.
\vspace{-02.0pt}
\end{list}
\vspace{-05.7pt}
\end{itemize}
We have to check that (i) and (ii) and the consequent of
$\phi_{\srfs{Fxd-img}}^{\srfs{$\exists$-scnst}}(\bfif,\bfib\kern1pt)$
can be expressed by a $\calL^{\srfs{ActGr}}_{\srfs{FT}}$-formula.
Indeed, $\bfia \stnot\stin \rfs{Cnst-dom}(\bfig\kern1pt)$ is equivalent to
$[\bfia\kern0.4pt]_{\sbfig} = \sngltn{\bfia\kern0.4pt}$, so (i) is expressible.
By Proposition~\ref{p5.8-08-15},
$\rfs{cnst}(\bfig\kern1pt) \stnot\stin [\bfib]_{\sbfif}$ is expressible.
The consequent:
$[\bfib]_{\sbfif \bfig \bfif\kern2pt} \subseteq [\bfib]_{\sbfif\fsup{2}}$,
is also expressible.

\begin{prop}\label{p1.11}
{\thickmuskip=2mu \medmuskip=1mu \thinmuskip=1mu
Let $S \stin K_{\srfs{FT}}$, $f \stin S$ and $b \stin A$.
Suppose that
$f(b) \stin \rfs{Fxd}(f)$.
}
\underline{Then}
$\rfs{Act}^{\srfs{Gr}}(S) \models
\phi_{\srfs{Fxd-img}}^{\srfs{$\exists$-scnst}}[f,b]$.
\end{prop}

\noindent
{\bf Proof }
Let $a \eqdf f(b)$.
We show that $a$ serves as the required $\bfia$ in
$\phi_{\srfs{Fxd-img}}^{\srfs{$\exists$-scnst}}$.
Let $g \stin S$ be a semi-constant function fulfilling (i) and (ii) of
$\phi_{\srfs{Fxd-img}}^{\srfs{$\exists$-scnst}}$.
Let $x \stin [b]_{f g f}$ and we prove that $x \stin [b]_{f^2}$.
\smallskip 
\\ 
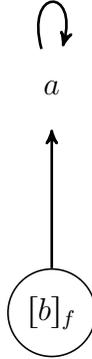
\centerline{\begin{tikzpicture}[->,>=stealth',shorten >=1pt,auto,node distance=3.0cm,on grid,semithick,
                    every state/.style={fill=white,draw=black,circle ,text=black}]
  \node [state, draw=none] (A)                                                                               {$a$};
  \node [state] (B)  [below = of A]                                                    {$[b]_f$};
\path (B) [line width=1pt] edge node {} (A);
\path  (A) [line width=1pt]  edge [loop above]        node [swap]{}(\Cnr);
\end{tikzpicture}}
$$
f g f(b) = f g(a) = f(a) = a.
$$
(The second equality holds since $g$ fulfills
(i) of $\phi_{\srfs{Fxd-img}}^{\srfs{$\exists$-scnst}}$.)
We conclude\break
that $f g f(x) = a$.
Suppose by contradiction that $g f(x) \neq f(x)$.
Then\break
$g f(x) = \rfs{cnst}(g)$.
The fact that $g$ fulfills
(ii) of $\phi_{\srfs{Fxd-img}}^{\srfs{$\exists$-scnst}}$ means that\break
$g f(x) = \rfs{cnst}(g) \stnot\stin [b]_f$. That is, $f g f(x) \neq f(b) = a$.
A contradiction.
It follows that $g f(x) = f(x)$. So $a = f g f(x) = f^2(x)$.
For $b$, too, $f^2(b) = f(a) = a$.
So $x \stin [b]_{f^2}$.
We have proved that $[b]_{f g f} \subseteq [b]_{f^2}$.
\hfill\solidqed
\vspace{2mm}

\begin{prop}\label{p1.12}
Let $S \stin K_{\srfs{FT}}$, and suppose that
$S$ contains a semi-constant function, but does not contain constant functions.
Let $f \stin S$ and $b \stin A$. Assume that $a \eqdf f(b) \stnot\stin \rfs{Fxd}(f)$.
\underline{Then}
$\rfs{Act}^{\srfs{Gr}}(S) \stnot\models
\phi_{\srfs{Fxd-img}}^{\srfs{$\exists$-scnst}}[f,b]$.
\end{prop}

\noindent
{\bf Proof }
Let $c = f(a)$. So $c \neq a$.
We have to show that a statement of the form
$$\mbox{
``There is $\bfia \stin [b]_f$ such that $\ldots$''
}$$
does not hold.
So let $b\fprime$ be an arbitrary member of $[b]_f$,
and we show that $b\fprime$ cannot serve as the needed $\bfia$.

We start by choosing a $d$ such that
$$
d \stin \rfs{Rng}(f) \swsetminus (f\inverse[\sngltn{c}] \stcup \sngltn{b\fprime}).
$$
We have to show that such a $d$ exists. Suppose otherwise.
Then $\rfs{Rng}(f) \subseteq f\inverse[\sngltn{c}] \stcup \sngltn{b\fprime}$
and so $\rfs{Rng}(f^2) \subseteq \dbltn{c}{f(b\fprime)}$.
Since $\abs{A} > 2$, $\rfs{Rng}(f^2) \neq A$. Also, $\rfs{Rng}(f^2)$ is finite.
So by Observation~\ref{o1.2}(b), $S$ contains a constant function.
A contradiction, so the required $d$ exists.

Now we choose a semi-constant function $g \stin S$ such that
\begin{itemize}
\addtolength{\parskip}{-11pt}
\addtolength{\itemsep}{06pt}
\item[(1)]
$d \stin \rfs{Cnst-dom}(g)$,
\item[(2)]
$b\fprime \stnot\stin \rfs{Cnst-dom}(g)$ and
\item[(3)]
$a = \rfs{cnst}(g)$.
\vspace{-05.7pt}
\end{itemize}
We have to show that such a $g$ exists.

Let us first verify that $b\fprime \neq a,d$.
(1)
$d$ was chosen to be different from $b\fprime$.
(2)
Since $b\fprime \stin [b]_f$, $f(b\fprime) = f(b) = a$.
However, $f(a) = c \neq a$. So $b\fprime \neq a$.

Recall that $S$ contains a semi-constant function which is not constant.
We have also shown that $b\fprime \neq a,d$.
So by Observation~\ref{o5.7-08-15}, such a $g$ exists.

Let $e$ be such that $f(e) = d$.
\smallskip 
\\ 
\centerline{\begin{tikzpicture}[->,>=stealth',shorten >=1pt,auto,node distance=3.0cm,on grid,semithick,
                    every state/.style={fill=white,draw=none,circle ,text=black}]
 \node [state] (O)                                                                                 
    {$a$};                   
  \node [state] (A)[above= of O]
   {$c$};
    \node [state] (C)  [below = of O]                                                           {$b$};
    \node [state] (B)  [left = of C]                                                    {$d$};
  \node [state] (D)  [right= of C]                                                    {$b^{\prime}$};
  \node [state] (E)  [below = of B]                                                    {$e$};
\path (O) [line width=1pt] edge node {} (A);
\path (C) [line width=1pt] edge node {} (O);
\path (D) [line width=1pt] edge node {} (O);
\path (E) [line width=1pt] edge  node {} (B);
\path (B) [red, line width=1pt] edge node {$g$} (O);
\path  (O) [red, line width=1pt]  edge [loop left ]  node [swap]{}(\Cnr);
\end{tikzpicture}}
\\
Then
$$
f g f(e) = f(g(f(e))) = f(g(d)) = f(a) = c.
$$
Also,
$$
f g f(b) = f(g(a)) = f(a) = c.
$$
So $e \stin [b]_{f g f}$.
However,
$$
f^2(e) = f(d) \neq c\ \ \ %
and\ \ \ %
f^2(b) = c.
$$
So $e \stnot\stin [b]_{f^2}$,
which implies that $[b]_{f g f} \stnot\subseteq [b]_{f^2}$.
We found a semi-constant function $g$ such that
(i) $b\fprime \stnot\stin \rfs{Cnst-dom}(g)$,
(ii) $\rfs{cnst}(g) \stnot\stin [b]_f$,
and nevertheless, $[b]_{f g f} \stnot\subseteq [b]_{f^2}$.
So $b\fprime$ does not fulfill the requirement on $\bfia$ in
$\phi_{\srfs{Fxd-img}}^{\srfs{$\exists$-scnst}}(\bfif,\bfib\kern1pt)$.
Hence $\rfs{Act}^{\srfs{Gr}}(S) \stnot\models
\phi_{\srfs{Fxd-img}}^{\srfs{$\exists$-scnst}}[f,b]$.
\hfill\solidqed
\vspace{2mm}

Combining the last two propositions one obtains the following statement.

\begin{lemma}\label{l1.13}\label{l5.11-08-16}
Let $S \stin K_{\srfs{$\exists$-scnst}}$.
Then for every $f \stin S$ and $b \stin A$:
\reqnomode
\begin{align}
\tag*{\hollowqed}
\mbox{
$\rfs{Act}^{\srfs{Gr}}(S) \models \phi_{\srfs{Fxd-img}}^{\srfs{$\exists$-scnst}}[f,b]$
\ \ \ iff\ \ \ \hspace{2pt}%
$f(b) \stin \rfs{Fxd}(f)$.
}
\end{align}
\end{lemma}

\section{\bf The definability of the relation
``$\boldsymbol{f(b) \stin \bfs{Fxd}(f)}$''
in FT semigroups which do not contain a semi-constant function,
the case that $\abs{\bfB\boldsymbol{(f)}} \geq \hbox{\bf 3}$}\label{s7}

Recall that $K_{\srfs{No-scnst}}$ is the class of all semigroups
$S \stin K_{\srfs{FT}}$ such that $S$ does not contain a semi-constant function.
The classes $J_2,J_3$ and $J_4$ were defined in Definition~\ref{d5.3},
and $\threetn{J_2}{J_3}{J_4}$ is a partition of
$\setm{\pair{S}{f}}{S \stin K_{\srfs{No-scnst}} \mbox{ and } f \stin S}$.
\vspace{2mm}

{\bf The goal for the next two sections.}
\\
We shall find an $\calL^{\srfs{ActGr}}_{\srfs{FT}}$-formula
$\phi_{\srfs{Fxd-img}}^{\srfs{No-scnst}}(\bfif,\bfib\kern1pt)$
which in semigroups belonging to $K_{\srfs{No-scnst}}$,
will say in that $\bfif\kern1pt(\bfib\kern1pt)$ is a fixed point of $\bfif$.

In this section we find such a formula which works only for
$\pair{S}{f}$'s which belong to $J_2$,
namely, $\pair{S}{f}$'s such that $\abs{\bfB(f)} \geq 3$.
In the next section we shall find the analogous formulas for $J_3$ and $J_4$.
\vspace{2mm}

We also need the facts that for $i = 2,3,4$,
$\setm{\pair{\rfs{Act}^{\srfs{Gr}}(S)}{f}}{\pair{S}{f} \stin J_i}$
is FO-definable in
$\setm{\pair{\rfs{Act}^{\srfs{Gr}}(S)}{f}}
{S \stin K_{\srfs{FT}} \mbox{ and } f \stin S}$.
These facts are rather trivial, and we verify them now.

\begin{observation}\label{o7.1}
{\rm(a)}
There is a sentence $\phi_{\srfs{No-scnst}}$
such that for every\break
$S \stin K_{\srfs{FT}}$,
$\rfs{Act}^{\srfs{Gr}}(S) \models \phi_{\srfs{No-scnst}}$
iff $S \stin K_{\srfs{No-scnst}}$.

{\rm(b)}
There is a formula $\phi_{\srfs{\sbf{B}$\geq$3}}(\bfif\kern2.3pt)$
such that for every $S \stin K_{\srfs{FT}}$ and every $f \stin S$:
$\abs{\bfB(f)} \geq 3$ iff
$\rfs{Act}^{\srfs{Gr}}(S) \models \phi_{\srfs{\sbf{B}$\geq$3}}[f]$.

{\rm(c)}
For every $i = 2,3,4$,
$\setm{\pair{\rfs{Act}^{\srfs{Gr}}(S)}{f}}{\pair{S}{f} \stin J_i}$
is FO-definable in
$\setm{\pair{\rfs{Act}^{\srfs{Gr}}(S)}{f}}
{S \stin K_{\srfs{FT}} \mbox{ and } f \stin S}$.
\end{observation}

\noindent
{\bf Proof }
(a)
Part (a) follows trivially from Proposition~\ref{p5.12-08-15}(a) and (b).

(b)
Part (b) follows from the fact that there is a formula
$\phi_{\srfs{sim}}(\bfif,\bfix,\bfiy\kern1pt)$ such that
for every $S \stin K_{\srfs{FT}}$,
$f \stin S$ and $a,b \stin A$:
$\rfs{Act}^{\srfs{Gr}}(S) \models \phi_{\srfs{sim}}[f,a,b]$
iff $a \sim_f b$. See Observation~\ref{o1.7}.
Indeed $\phi_{\srfs{\sbf{B}$\geq$3}}(\bfif\kern2pt)$ is the following formula.
\index{N@fphib@@$\phi_{\srfs{\sbf{B}$\geq$3}}(\bfif\kern2pt)$}%
{\thickmuskip=2mu \medmuskip=1mu \thinmuskip=1mu
\begin{align*}
(\exists \bfix\Fsub{1}, \bfix\Fsub{2}, \bfix\Fsub{3})
(\exists \bfiy\Fsub{1}, \bfiy\Fsub{2}, \bfiy\Fsub{3})
\left(
\bigwedge_{i = 1}^3 (\bfix\Fsub{i} \neq \bfiy\Fsub{i})
\stwedge
\bigwedge_{i = 1}^3 (\bfix\Fsub{i} \sim_{\sbfif} \bfiy\Fsub{i})
\kern3pt\stwedge 
\bigwedge_{1 \leq i < j \leq 3} (\bfix\Fsub{i} \stnot\sim_{\sbfif} \bfix\Fsub{j})
\right)
\end{align*}
}

(c)
In Part (a) we found a sentence $\phi_{\srfs{No-scnst}}$
such that for every $S \stin K_{\srfs{FT}}$,
$\rfs{Act}^{\srfs{Gr}}(S) \models \phi_{\srfs{No-scnst}}$
iff $S \stin K_{\srfs{No-scnst}}$.
In Part (b) we found a formula
$\phi_{\srfs{\sbf{B}$\geq$3}}(\bfif\kern2.3pt)$
such that for every $S \stin K_{\srfs{FT}}$ and every $f \stin S$:
$\abs{\bfB(f)} \geq 3$ iff
$\rfs{Act}^{\srfs{Gr}}(S) \models \phi_{\srfs{\sbf{B}$\geq$3}}[f]$.
In Proposition~\ref{p1.4} we found a formula
$\phi_{\srfs{Gr}}(\bfif\kern2.3pt)$
such that for every $S \stin K_{\srfs{FT}}$ and $f \stin S$:
$f \stin \rfs{Gr}(S)$ iff
$\rfs{Act}^{\srfs{Gr}}(S) \models \phi_{\srfs{Gr}}[f]$.

Now define
$$
\phi^{J_2}(\bfif\kern2.3pt) \eqdf
\phi_{\srfs{No-scnst}} \stwedge
\phi_{\srfs{\sbf{B}$\geq$3}}(\bfif\kern2.3pt),
$$
$$
\phi^{J_3}(\bfif\kern2.3pt) \eqdf
\phi_{\srfs{No-scnst}} \stwedge
\stneg\phi_{\srfs{\sbf{B}$\geq$3}}(\bfif\kern2.3pt) \stwedge
\stneg\phi_{\srfs{Gr}}(\bfif\kern2.3pt)
$$
and
$$
\phi^{J_4}(\bfif\kern2.3pt) \eqdf
\phi_{\srfs{No-scnst}} \stwedge \phi_{\srfs{Gr}}(\bfif\kern2.3pt).
$$
Clearly, $\phi^{J_2}(\bfif\kern2.3pt)$, $\phi^{J_3}(\bfif\kern2.3pt)$
and $\phi^{J_4}(\bfif\kern2.3pt)$ are as required.
\index{N@fphij2@@$\phi^{J_2}(\bfif\kern2.3pt)$}%
\index{N@fphij3@@$\phi^{J_3}(\bfif\kern2.3pt)$}%
\index{N@fphij4@@$\phi^{J_4}(\bfif\kern2.3pt)$}%
\hfill\solidqed
\vspace{2mm}

We shall need the following notation. Let $a,b,c \stin A$ be distinct.
Define
$\rvtrpl{a}{b}{c}^A \eqdf 
\threetn{\pair{a}{b}}{\pair{b}{c}}{\pair{c}{a}}
\stcup \rfs{Id}_{A \setminus \threetn{\pair{a}{b}}{\pair{b}{c}}{\pair{c}{a}}}$.
The superscript $A$ is usually omitted.
\index{N@AAAA@@$\rvtrpl{a}{b}{c}$}%

The next Observation~\ref{o1.2}(b) is an analogue of Observation~\ref{o1.2}(b).

\begin{observation}\label{o2.1}\label{o7.2}
Let $S \subseteq A^A$ be an FT semigroup.
If $S$ contains a function~$f$ such that
{\rm(i)} $f[A \swsetminus \rfs{Idp}(f)]$ is finite,
and {\rm(ii)} $f$ is not $\onetoone$,
\underline{then} $S$ contains a semi-constant function.
\end{observation}

\noindent
{\bf Proof }
For a function $f \stin A^A$ set
$\bfs{n}(f) \eqdf \abs{f[A \swsetminus \rfs{Idp}(f)]}$.
Note that $\bfs{n}(f) = 1$ iff $f$ is semi-constant.
Let $S \stin K_{\srfs{FT}}$, and $A = \bfs{Dom}(S)$. Suppose that there is
$f \stin S$ such that
$f[A \swsetminus \rfs{Idp}(f)]$ is finite, and $f$ is not 1-1.
Let
$m = \min\left(\rule{0pt}{12pt}
\setm{\bfs{n}(g)}{g \stin S, \mbox{ and $g$ is not 1-1}}
\right)$.
Suppose by contradiction that $m > 1$.
Let $g \stin S$ be such that $\bfs{n}(g) = m$.
Let $a,b \stin A$ be distinct elements such that $g(a) = g(b)$,
and let $c,d \stin g[A \swsetminus \rfs{Idp}(g)]$ be distinct.
There is $\pi \stin S$ such that (i) $\pi(a) = c$ and $\pi(b) = d$,
and\break
(ii) $\pi[\rfs{Idp}(g)] = \rfs{Idp}(g)$.
We verify that, indeed, such a $\pi$ exists.
Note that $a,b,c,d \stnot\stin \rfs{Idp}(g)$.

{\bf Case 1 } $\dbltn{a}{b} \stcap \dbltn{c}{d} = \stemptyset$.
Then $\pi \eqdf \rvpair{a}{c} \rvpair{b}{d}$ is as required.

{\bf Case 2 } $\dbltn{a}{b} \stcap \dbltn{c}{d}$ is a singleton.
We may assume that $b = c$.
Then $\pi \eqdf \rvtrpl{a}{b}{d}$ is as required.

{\bf Case 3 } $\dbltn{a}{b} = \dbltn{c}{d}$.
Then $\pi \eqdf \rfs{Id}_A$ is as required.

Clearly, $\pi \stin S$ and $g^{\pi}(c) = g^{\pi}(d)$.
Hence $h \eqdf g^{\pi} g \stin S$.
It follows that $\rfs{Idp}(h) = \rfs{Idp}(g)$.
The fact ``\kern1pt$g^{\pi}(c) = g^{\pi}(d)\kern1pt$'' implies that
\vspace{-2mm}
$$
\abs{h[A \swsetminus \rfs{Idp}(h)]} < m.
\vspace{-2mm}
$$
A contradiction. So $m = 1$. This means that $S$ contains a semi-constant
function.
\hfill\solidqed
\vspace{2mm}

We introduce two formulas $\alpha_{\srfs{Not-fxd}}(\bfif,\bfib\kern1pt)$ and
$\beta_{\srfs{Not-fxd}}(\bfif,\bfib\kern1pt)$ which will be used in
the construction of a formula expressing
the property that\break
$\bfif(\bfib\kern1pt) \stnot\stin \rfs{Fxd}(\bfif\kern2pt)$.
If $f(b) \stin \rfs{Fxd}(f)$, then both $\alpha_{\srfs{Not-fxd}}$
and $\beta_{\srfs{Not-fxd}}$ fail for $f$ and~$b$. 
If $f(b) \stnot\stin \rfs{Fxd}(f)$, then {\it in most cases} either
$\alpha_{\srfs{Not-fxd}}$
or $\beta_{\srfs{Not-fxd}}$ hold for $f$ and~$b$.
However, there are a few instances in which
$f(b) \stnot\stin \rfs{Fxd}(f)$, but\break
$\rfs{Act}^{\srfs{Gr}}(S) \stnot\models
(\alpha_{\srfs{Not-fxd}} \vee \beta_{\srfs{Not-fxd}})[f,b]$.
Dealing with these exceptional cases requires some extra work.

Let $\fnn{f}{A}{A}$, $b \stin A$ and $\pi = \rvpair{c}{d}$ be a transposition
of $A$.
We say that $\pi$ is \underline{$\pair{f}{b}$\kern2pt-\kern1ptpermissible},
if $c,d \stnot\stin [b]_f$.

Let $\alpha_{\srfs{Not-fxd}}(\bfif,\bfib\kern1pt)$ be the formula which says:
\begin{itemize}
\addtolength{\parskip}{-11pt}
\addtolength{\itemsep}{06pt}
\item
There is an $\pair{\bfif}{\bfib}$\kern2pt-\kern1ptpermissible transposition
$\pi$ such that
$[\bfib]_{\sbfif\Fsup{2}} \neq [\bfib]_{\sbfif\kern1pt \pi \sbfif}\kern1pt$.%
\vspace{-05.7pt}
\end{itemize}
Let $\beta_{\srfs{Not-fxd}}(\bfif,\bfib\kern1pt)$ be the formula which says:
\begin{itemize}
\addtolength{\parskip}{-11pt}
\addtolength{\itemsep}{06pt}
\item
There is an $\pair{\bfif}{\bfib}$\kern2pt-\kern1ptpermissible transposition
$\sigma$ such that $\alpha_{\srfs{Not-fxd}}(\bfif \sigma,\bfib\kern1pt)$.%
\vspace{-05.7pt}%
\end{itemize}

Since $\alpha_{\srfs{Not-fxd}}$ and $\beta_{\srfs{Not-fxd}}$ appear many times,
for a while, we abbreviate $\alpha_{\srfs{Not-fxd}}$ and $\beta_{\srfs{Not-fxd}}$
by respectively $\alpha$ and $\beta$.
\index{N@alphanot-fxd@@$\alpha_{\srfs{Not-fxd}}$. Abbreviated by $\alpha$}%
\index{N@betanot-fxd@@$\beta_{\srfs{Not-fxd}}$. Abbreviated by $\beta$}%

\begin{prop}\label{p2.2}\label{p7.3}
Let $S \stin K_{\srfs{FT}}$, $f \stin S$ and $b \stin A$.
Suppose that $f(b) \stin \rfs{Fxd}(f)$. Then
$\rfs{Act}^{\srfs{Gr}}(S) \stnot\models \alpha[f,b]$
and $\rfs{Act}^{\srfs{Gr}}(S) \stnot\models \beta[f,b]$.
\end{prop}

\noindent
{\bf Proof }
The proof is trivial.
\hfill\solidqed
\vspace{2mm}

\begin{definition}\label{d7.4}
\begin{rm}
Let $\fnn{f}{A}{A}$.

(a)
A triple $\trpl{b}{a}{c} \stin A^3$ is called an \underline{$f$-long triple},
if $f(b) = a$, $f(a) = c$ and $b \neq a \neq c$.
($b$ is allowed to be equal to $c$.)
\smallskip 
\\ 
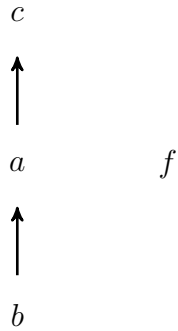
\centerline{\begin{tikzpicture}[->,>=stealth',shorten >=1pt,auto,node distance=2.0cm,on grid,semithick,
                    every state/.style={fill=white,draw=none,circle ,text=black}]
 \node [state] (O)                                                                                 
    {$a$};                   
  \node [state] (A)[above= of O]
   {$c$};
    \node [state] (C)  [below = of O]                                                           {$b$};
    \node [state] (B)  [right = of O]                                                    {$f$};
\path (O) [line width=1pt] edge node {} (A);
\path (C) [line width=1pt] edge node {} (O);
\end{tikzpicture}}
\\


(b)
A quadruple $\fourtpl{r}{p}{q}{t} \stin A^4$ is called an
\underline{$f$-long-wide quadruple},
if\break
(i) $f(p) = f(q) = t$, (ii) $p,q,t$ are pairwise distinct and (iii) $f(r) = q$.
\smallskip 
\\ 
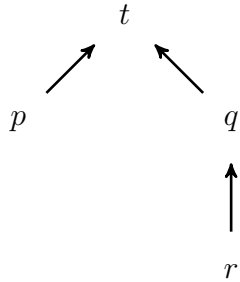
\centerline{\begin{tikzpicture}[->,>=stealth',shorten >=1pt,auto,node distance=2.0cm,on grid,semithick,
                    every state/.style={fill=white,draw=none,circle ,text=black}]
 \node [state] (O)                                                                                 
    {$t$};                   
  \node [state] (A)[below left= of O]
   {$p$};
  \node [state] (C)  [below right = of O]                                                           {$q$};
   \node [state] (B)  [below = of C]                                                           {$r$};
\path (A) [line width=1pt] edge node {} (O);
\path (C) [line width=1pt] edge node {} (O);
\path (B) [line width=1pt] edge node {} (C);
\end{tikzpicture}}
%
\end{rm}
\end{definition}

\begin{remark}\label{r7.5}
Note the following equivalence.
Let $\fnn{f}{A}{A}$ and $b \stin A$.
Then the following are equivalent:
{\rm(1)} $b \stnot\stin \rfs{Fxd-Img}(f)$;
{\rm(2)} $\trpl{b}{f(b)}{f^2(b)}$ is an $f$-long triple.
\end{remark}

\begin{prop}\label{p2.6}\label{p7.6}
Let $S \stin K_{\srfs{FT}}$,
$f \stin S$, $\trpl{b}{a}{c}$ be an $f$-long triple,
and\break
$\fourtpl{r}{q}{s}{t}$ be an $f$-long-wide quadruple such that $t \neq a,c$.
\underline{Then} $\rfs{Act}^{\srfs{Gr}}(S) \models \alpha[f,b]$.
\end{prop}

\noindent
{\bf Proof }
\smallskip 
\\ 
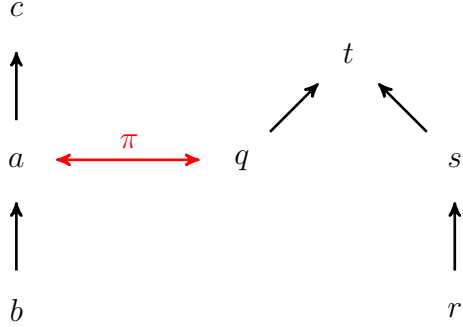
\centerline{\begin{tikzpicture}[->,>=stealth',shorten >=1pt,auto,node distance=2.0cm,on grid,semithick,
                    every state/.style={fill=white,draw=none,circle ,text=black}]
 \node [state] (G)                                                                                 
    {$a$}; 
  \node [state] (O) [node distance=3.0cm, right = of G]                                                                                                                                                    
    {$q$};                   
  \node [state] (A)[above right= of O]
  {$t$};
  \node [state] (D)[below right= of A]
   {$s$};
  \node [state] (E)[below= of D]
   {$r$};
  \node [state] (C)[above= of G]
   {$c$};
  \node [state] (B)  [below = of G]                                                    {$b$};
\path (O) [line width=1pt] edge node {} (A);
\path (B) [line width=1pt] edge node {} (G);
\path (G) [line width=1pt] edge node {} (C);
\path (D) [line width=1pt] edge node {} (A);
\path (E) [line width=1pt] edge node {} (D);
\path   (G) [<->, red,line width=1pt] edge node {$\pi$}  (O) ; 
\end{tikzpicture}}
\\
Since $f(a) = c \neq t = f(q)$, $a \neq q$. So $\pi \eqdf \rvpair{a}{q}$
is a transposition.
Now, $f(q) = t \neq a = f(b)$. So $q \stnot\stin [b]_f$.
Also $f(a) = c \neq a = f(b)$. So $a \stnot\stin [b]_f$.
It follows that
\begin{equation}\mbox{
\tag{1}
$\pi$ is $\pair{f}{b}$\kern2pt-\kern1ptpermissible.
}\end{equation}
Clearly, $f^2(r) = t \neq c = f^2(b)$.
So
\begin{equation}\mbox{
\tag{2} $r \stnot\stin [b]_{f^2}$.
}\end{equation}
{\thickmuskip=2mu \medmuskip=1mu \thinmuskip=1mu
Since $\fourtpl{r}{q}{s}{t}$ is an $f$-long-wide quadruple, $s \neq q$.
Also, $f(s) = t \neq c = f(a)$. So $s \neq a$.}
Hence $\pi(s) = s$. It follows that
\begin{equation}
\tag{3}
f \pi f(r) = f \pi(s) = f(s) = t.
\end{equation}

Since $f(s) = t \neq c = f(a)$, $s \neq a$.
Also, since $\fourtpl{r}{q}{s}{t}$ is an $f$-long-wide quadruple, $s \neq q$.
So $\pi(s) = s$. Hence
\begin{equation}
\tag{4}
f \pi f(b) = f \pi(a) = f(q) = t.
\end{equation}
By (3) and (4)
\begin{equation}
\tag{5}
r \stin [b]_{f \pi f}.
\end{equation}
By (2) and (5),
\begin{equation}
\tag{6}
[b]_{f \pi f} \neq [b]_{f^2}.
\end{equation}
By (1) and (6), $\rfs{Act}^{\srfs{Gr}}(S) \models \alpha[f,b]$.
\hfill\solidqed
\vspace{2mm}

Let $f \stin A^A$, $B \stin \bfs{B}(f)$ and $a \stin B$.
Set $\rfs{Img}_f(B) \eqdf f(a)$.
\index{N@img2@@$\rfs{Img}_f(B)$}%

\begin{lemma}\label{l7.7}
Let $\pair{S}{f} \stin J_2$
$($that is,
$S \stin K_{\srfs{No-Scnst}}$ and $\abs{\bfs{B}(f)} \geq 3$$)$.
Then for every $b \stin A \swsetminus \rfs{Fxd-img}(f)$,
$\rfs{Act}^{\srfs{Gr}}(S) \models \beta[f,b]$.
\end{lemma}

\noindent
{\bf Proof }
Set $f(b) = a$ and $f(a) = c$.
So $\trpl{b}{a}{c}$ is an $f$-long triple.
We show that there is transposition $\sigma$ such that
$\sigma$ is $\pair{f}{b}$-permissible,
and $\rfs{Act}^{\srfs{Gr}}(S) \models \alpha[f \sigma,b]$.

Since $\abs{\bfB(f)} \geq 3$, there is $B \stin \bfB(f)$
such that $ t \eqdf \rfs{img}_f(B) \neq c,a$.
Let $p,q \stin B$ be distinct. So $p \neq t$ or $q \neq t$.
We may assume that $p \neq t$.
We shall find $\sigma \stin S$ and $u,v \stin A$ such that
$\qdrpl{u}{p}{v}{t}$ is an $f \sigma$-long-wide quadruple.

We apply Observation~\ref{o7.2}.
It is given that $S$ does not contain a semi-constant function,
and that $f$ is not 1-1. So $f^2$ is not 1-1.
By Observation~\ref{o7.2}, $f^2[A \swsetminus \rfs{Idp}(f^2)]$ is infinite.
So there is $d$ such that
\vspace{2mm}
\begin{list}{}
{\setlength{\leftmargin}{39pt}
\setlength{\labelsep}{05pt}
\setlength{\labelwidth}{25pt}
\setlength{\itemindent}{-00pt}
\addtolength{\topsep}{-04pt}
\addtolength{\parskip}{-02pt}
\addtolength{\itemsep}{-05pt}
}
\item[{\bf (i)}]
\rule{03mm}{0pt}
$ d \stin f^2[A \swsetminus \rfs{Idp}(f^2)] \swsetminus
\left(\rule{0pt}{12pt}
\sixtn{a}{b}{c}{p}{q}{t} \stcup f[\sixtn{a}{b}{c}{p}{q}{t}]\right).
$
\end{list}
It follows that $d \stin \rfs{Rng}(f^2) \swsetminus \rfs{Idp}(f^2)$.
This implies that there is $u \neq d$ such that $f^2(u) = d$.
Let $v = f(u)$.
Hence
\begin{list}{}
{\setlength{\leftmargin}{39pt}
\setlength{\labelsep}{05pt}
\setlength{\labelwidth}{25pt}
\setlength{\itemindent}{-00pt}
\addtolength{\topsep}{-04pt}
\addtolength{\parskip}{-02pt}
\addtolength{\itemsep}{-05pt}
}
\item[{\bf (ii)}] $f(v) = d$.
\vspace{-02.0pt}
\end{list}
Recall that $f(q) = t$.
However, by {\bf(i)}, $d \neq t$. So
\begin{list}{}
{\setlength{\leftmargin}{39pt}
\setlength{\labelsep}{05pt}
\setlength{\labelwidth}{25pt}
\setlength{\itemindent}{-00pt}
\addtolength{\topsep}{-04pt}
\addtolength{\parskip}{-02pt}
\addtolength{\itemsep}{-05pt}
}
\item[{\bf (iii)}] $v \neq q$.
\vspace{-02.0pt}
\end{list}

The facts: $u \neq d$, $f(u) = v$ and $f(v) = d$ imply that
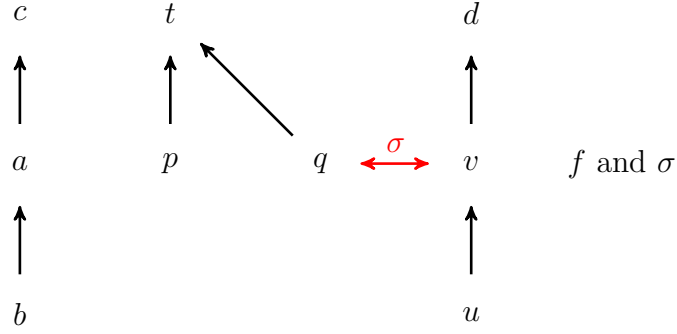
\begin{list}{}
{\setlength{\leftmargin}{39pt}
\setlength{\labelsep}{05pt}
\setlength{\labelwidth}{25pt}
\setlength{\itemindent}{-00pt}
\addtolength{\topsep}{-04pt}
\addtolength{\parskip}{-02pt}
\addtolength{\itemsep}{-05pt}
}
\item[{\bf (iv)}] $u \neq v$.
\vspace{-02.0pt}
\end{list}

\rule{0pt}{0pt}
\vspace{-4mm}
\smallskip 
\\ 
\centerline{\begin{tikzpicture}[->,>=stealth',shorten >=1pt,auto,node distance=2.0cm,on grid,semithick,
                    every state/.style={fill=white,draw=none,circle ,text=black}]
 \node [state] (C)                                                                                 
    {$c$}; 
  \node [state] (A)[below = of C]
  {$a$};
  \node [state] (B)[below = of A]
  {$b$};
  \node [state] (T)[right= of C]
   {$t$};
  \node [state] (P)[below= of T]
   {$p$};
  \node [state] (Q)[right= of P]
   {$q$};
   \node [state] (V)[right= of Q]
   {$v$};
  \node [state] (U)[below= of V]
   {$u$};
  \node [state] (D)  [above = of V] 
   {$d$};
   \node [state] (F)  [right = of V] 
   {$f \mbox{ and } \sigma$};
\path (B) [line width=1pt] edge node {} (A);
\path (A) [line width=1pt] edge node {} (C);
\path (P) [line width=1pt] edge node {} (T);
\path (Q) [line width=1pt] edge node {} (T);
\path (U) [line width=1pt] edge node {} (V);
\path (V) [line width=1pt] edge node {} (D);
\path   (Q) [<->, red,line width=1pt] edge node {$\sigma$}  (V) ; 
\end{tikzpicture}}
\\
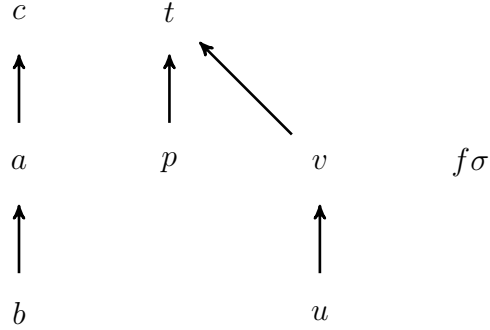
\centerline{\begin{tikzpicture}[->,>=stealth',shorten >=1pt,auto,node distance=2.0cm,on grid,semithick,
                    every state/.style={fill=white,draw=none,circle ,text=black}]
 \node [state] (C)                                                                                 
    {$c$}; 
  \node [state] (A)[below = of C]
  {$a$};
  \node [state] (B)[below = of A]
  {$b$};
  \node [state] (T)[right= of C]
   {$t$};
  \node [state] (P)[below= of T]
   {$p$};
  \node [state] (V)[right= of P]
   {$v$};
  \node [state] (U)[below= of V]
   {$u$};
  \node [state] (F)  [right = of V] 
   {$f \sigma$};
\path (B) [line width=1pt] edge node {} (A);
\path (A) [line width=1pt] edge node {} (C);
\path (P) [line width=1pt] edge node {} (T);
\path (V) [line width=1pt] edge node {} (T);
\path (U) [line width=1pt] edge node {} (V);
\end{tikzpicture}}
\\

We prove that {\bf(I)}:
$\sigma \eqdf \rvpair{v}{q}$ is well-defined,
$\sigma \stin S$, and $\sigma$ is\break
$\pair{f}{b}$-permissible.

By {\bf (iii)}, $\sigma \eqdf \rvpair{v}{q}$ is well-defined,
and since $S$ is fully-transpositional, $\sigma \stin S$.

Now we check that $\sigma$ is $\pair{f}{b}$-permissible.
Of course, $f(b) = a$.
We show that $f(v), f(q) \neq a$.
By {\bf (ii)}, $f(v) = d$, and by {\bf(i)}, $d \neq a$. So $f(v) \neq a$.

Recall that $f(q) = t = \rfs{img}_f(B) \neq a$. So $f(q) \neq a$.
It follows that for every $v,q \stnot\stin [b]_f$.
That is, $\sigma$ is $\pair{f}{b}$-permissible.
We proved {\bf(I)}.
\vspace{2mm}

Next we show that {\bf(II)}:
$\fourtpl{u}{p}{v}{t}$ is an $f \sigma$-long-wide quadruple.
We check that $u \stnot\stin \dbltn{v}{q}$.
By {\bf(iv)}, $u \neq v$.
By {\bf(i)}, $d \neq f(t)$. So $d \neq t$.
This fact together with facts: $f(v) = d$ and $f(q) = t$, imply that $v \neq q$.
So indeed, $u \stnot\stin \dbltn{v}{q}$.

Recall that $\sigma = \rvpair{v}{q}$. Hence it follows that $\sigma(u) = u$,
and so, $f \sigma(u) = f(u) = v$. That is,
\begin{list}{}
{\setlength{\leftmargin}{39pt}
\setlength{\labelsep}{05pt}
\setlength{\labelwidth}{25pt}
\setlength{\itemindent}{-00pt}
\addtolength{\topsep}{-04pt}
\addtolength{\parskip}{-02pt}
\addtolength{\itemsep}{-05pt}
}
\item[{\bf (v)}] $f \sigma(u) = v$.
\vspace{-02.0pt}
\end{list}

Now, $f \sigma(v) = f(q) = t$. That is,
\begin{list}{}
{\setlength{\leftmargin}{39pt}
\setlength{\labelsep}{05pt}
\setlength{\labelwidth}{25pt}
\setlength{\itemindent}{-00pt}
\addtolength{\topsep}{-04pt}
\addtolength{\parskip}{-02pt}
\addtolength{\itemsep}{-05pt}
}
\item[{\bf (vi)}] $f \sigma(v) = t$.
\vspace{-02.0pt}
\end{list}
Next we check that
\begin{list}{}
{\setlength{\leftmargin}{39pt}
\setlength{\labelsep}{05pt}
\setlength{\labelwidth}{25pt}
\setlength{\itemindent}{-00pt}
\addtolength{\topsep}{-04pt}
\addtolength{\parskip}{-02pt}
\addtolength{\itemsep}{-05pt}
}
\item[{\bf (vii)}]
$f \sigma(p) = t$.
\vspace{-02.0pt}
\end{list}
We have chosen $p$ and $q$ to be distinct members of $f\inverse[\sngltn{t}]$.
So $p \neq q$.
We also need to see that
\begin{list}{}
{\setlength{\leftmargin}{39pt}
\setlength{\labelsep}{05pt}
\setlength{\labelwidth}{25pt}
\setlength{\itemindent}{-00pt}
\addtolength{\topsep}{-04pt}
\addtolength{\parskip}{-02pt}
\addtolength{\itemsep}{-05pt}
}
\item[{\bf (viii)}]
$p \neq v$.
\vspace{-02.0pt}
\end{list}
Indeed, $f(p) = t$, $f(v) = d$,
and by {\bf(i)}, $t \neq d$. So $p \neq v$.

In the beginning of the proof we have assumed that
\begin{list}{}
{\setlength{\leftmargin}{39pt}
\setlength{\labelsep}{05pt}
\setlength{\labelwidth}{25pt}
\setlength{\itemindent}{-00pt}
\addtolength{\topsep}{-04pt}
\addtolength{\parskip}{-02pt}
\addtolength{\itemsep}{-05pt}
}
\item[{\bf (ix)}]
$p \neq t$,
\vspace{-02.0pt}
\end{list}
and we now verify that for $v$, too,
\begin{list}{}
{\setlength{\leftmargin}{39pt}
\setlength{\labelsep}{05pt}
\setlength{\labelwidth}{25pt}
\setlength{\itemindent}{-00pt}
\addtolength{\topsep}{-04pt}
\addtolength{\parskip}{-02pt}
\addtolength{\itemsep}{-05pt}
}
\item[{\bf (x)}]
$v \neq t$.
\vspace{-02.0pt}
\end{list}
$f(v) = d$ and by {\bf(i)}, $f(t) \neq d$. So, indeed, $v \neq t$.
By Facts {\bf(viii)} - {\bf(x)},\break
$p,v,t$ are pairwise distinct;
and Facts {\bf(v)} - {\bf(vii)} say that
$f \sigma(u) = v$, $f \sigma(v) = t$ and $f \sigma(p) = t$.
Hence $\fourtpl{u}{p}{v}{t}$ is, indeed, an $f \sigma$-long-wide quadruple.
This proves {\bf(II)}.
\vspace{2mm}

Let us now see that {\bf(III)}: $\trpl{b}{a}{c}$ is $f \sigma$-long.
We check that
\vspace{-2mm}
$$
\dbltn{b}{a} \stcap \dbltn{v}{q} = \stemptyset.
\vspace{-2mm}
$$
Also, $f(q) = t$ and $t \neq a,c$. Since $f(b) = a$ and $f(a) = c$,
it follows that $q \neq b,a$.
The facts: $f(v) = d$, $d \neq a,c$, $f(b) = a$ and $f(a) = c$, imply that
$v \neq b,a$.
Hence $\dbltn{b}{a} \stcap \dbltn{v}{q} = \stemptyset$.
We conclude that $\sigma(b) = b$ and $\sigma(a) = a$.
So $b \neq a \neq c$, $f \sigma(b) = a$ and $f \sigma(a) = c$.
That is, $\trpl{b}{a}{c}$ is an $f \sigma$-long triple.
We proved {\bf(III)}.
\vspace{2mm}

In the beginning of the proof we chose $t$ such that
{\bf (IV)} $t \neq c,a$.
\vspace{2mm}

Facts {\bf(II)}, {\bf(III)} and {\bf(IV)} together with Proposition~\ref{p7.6}
imply that\break
{\thickmuskip=2mu \medmuskip=1mu \thinmuskip=1mu
{\bf(V)}:  $\rfs{Act}^{\srfs{Gr}}(S) \models \alpha[f \sigma,b]$.
Facts {\bf(I)} and {\bf(V)} imply that
$\rfs{Act}^{\srfs{Gr}}(S) \models \beta[f,b]$.}
\rule{20mm}{0pt}\hfill\solidqed

\begin{cor}\label{c2.9}\label{c6.11-08-15}\label{c7.8}
{\rm(a)}
Let $\pair{S}{f} \stin J_2$.
That is, $S \stin K_{\srfs{No-scnst}}$,
and $f \stin S$ is such that $\abs{\bfB(f)} \geq 3$.
Then for every $b \stin A$ the following are equivalent:
\begin{itemize}
\addtolength{\parskip}{-11pt}
\addtolength{\itemsep}{06pt}
\item[{\rm(1)}]
$b \stnot\stin \rfs{Fxd-img}(f)$.
\item[{\rm(2)}]
$\rfs{Act}^{\srfs{Gr}}(S) \models \beta[f,b]$.
\vspace{-05.7pt}
\end{itemize}

{\rm(b)}
Denote $\stneg\beta$ by $\phi^{J_2}_{\srfs{Fxd-img}}$.
\index{N@fphifxd-img-j2@@$\phi^{J_2}_{\srfs{Fxd-img}} \eqdf \stneg\beta_{\srfs{Not-fxd}}$}%
Then for every $\pair{S}{f} \stin J_2$ and $b \stin A$,
the following are equivalent:
\begin{itemize}
\addtolength{\parskip}{-11pt}
\addtolength{\itemsep}{06pt}
\item[{\rm(1)}]
$b \stin \rfs{Fxd-img}(f)$.
\item[{\rm(2)}]
$\rfs{Act}^{\srfs{Gr}}(S) \models \phi^{J_2}_{\srfs{Fxd-img}}[f,b]$.
\vspace{-05.7pt}
\end{itemize}
%
\end{cor}

\noindent
{\bf Proof }
(a) Part (a) follows trivially from Proposition~\ref{p7.3}
and Lemma~\ref{l7.7}.

(b) Part (b) is a restatement of Part (a).
It is mentioned for later reference.
%
\rule{20mm}{0pt}\hfill\solidqed
\vspace{2mm}

Recall that our goal was to find a formula
$\phi_{\srfs{Fxd-img}}^{\srfs{No-scnst}}(\bfif,\bfib\kern1pt)$
which in\break
FT semigroups not containing a semi-constant function, would say that
$\bfif\kern1pt(\bfib\kern1pt)$ is a fixed point of $\bfif$.
The formula $\phi^{J_2}_{\srfs{Fxd-img}}$ does the job
for $\pair{S}{f}$'s which belong to $J_2$.

\section{\bf Functions $\boldsymbol{f}$ for which
$\boldsymbol{\abs{\bfB(f)} \leq 2}$}\label{s8}

\subsection{An overview}\label{ss8.1}

We continue to adopt the following conventions.
\begin{itemize}
\addtolength{\parskip}{-11pt}
\addtolength{\itemsep}{06pt}
\item[(1)]
$S$ denotes a function semigroup,
and $\bfs{Dom}(S)$ is denoted by $A$.
\item[(2)]
All the formulas that we write are first order formulas in
$\calL_{\srfs{ActGr}}$.
When we say a ``formula'', we mean an $\calL_{\srfs{ActGr}}$-formula.
\vspace{0pt}
\end{itemize}

So far we have found formulas $\phi(\bfif,\bfix\kern2pt)$ equivalent to
``$\bfix \stin \rfs{Fxd-img}(\bfif\kern2pt)$''\break
for two subclasses of
$\setm{\pair{\rfs{Act}^{\srfs{Gr}}(S)}{f}}{S \stin K_{\srfs{No-cnst}} \mbox{ and }
f \stin S}$.
\begin{itemize}
\addtolength{\parskip}{-11pt}
\addtolength{\itemsep}{06pt}
\item[(1)]
The class $(J_1)^{\srfs{ActGr}}$,
where $J_1$ is the class of all $\pair{S}{f}$'s such that\break
$S \stin K_{\srfs{$\exists$-scnst}}$ and $f \stin S$,
\item[(2)]
The class $(J_2)^{\srfs{ActGr}}$,
where $J_2$ is the class of all $\pair{S}{f}$'s such that\break
$S \stin K_{\srfs{No-scnst}}$ 
and $f \stin \rfs{F}_{\srfs{\sbf{B}$\geq$3}}(S)$.
\vspace{-05.7pt}
\end{itemize}
For $J_1$, $\phi$ is the formula
$\phi_{\srfs{Fxd-img}}^{\srfs{$\exists$-scnst}}$,
which was defined just before Proposition~\ref{p1.11}.
(See also Lemma~\ref{l1.13}).
We now denote $\phi_{\srfs{Fxd-img}}^{\srfs{$\exists$-scnst}}$ by
$\phi^{J_1}_{\srfs{Fxd-img}}$.
\index{N@fphifxd-img-j1@@$\phi^{J_1}_{\srfs{Fxd-img}}$}%

For $J_2$, $\phi$ is the formula $\phi^{J_2}_{\srfs{Fxd-img}}$,
which was defined in Corollary~\ref{c6.11-08-15}. (See also the definition of
$\beta_{\srfs{Not-fxd}}$ just before Proposition~\ref{p2.2}.)

We have also established the FO-definability of $(J_i)^{\srfs{ActGr}}$ in\break
$\setm{\pair{\rfs{Act}^{\srfs{Gr}}(S)}{f}}{S \stin K_{\srfs{No-cnst}}
\mbox{ and } f \stin S}$.

That $(J_1)^{\srfs{ActGr}}$ is FO-definable in
$\setm{\pair{\rfs{Act}^{\srfs{Gr}}(S)}{f}}{S \stin K_{\srfs{No-cnst}}
\mbox{ and } f \stin S}$
was proved in Observation~\ref{o1.8}(d).
The defining formula was the sentence $\phi_{\srfs{$\exists$-scnst}}$,
which we now denote by $\phi^{J_1}$.
\index{N@fphij1@@$\phi^{J_1}$}%
For $i = 2,3,4$, the FO-definability of $(J_i)^{\srfs{ActGr}}$
in $\setm{\pair{\rfs{Act}^{\srfs{Gr}}(S)}{f}}{S \stin K_{\srfs{No-cnst}}
\mbox{ and } f \stin S}$
was proved in Observation~\ref{o7.1}(c).

What remains to be shown is the existence of a formula
$\phi^{J_3}_{\srfs{Fxd-img}}(\bfif,\bfix\kern1pt)$
which for members of $(J_3)^{\srfs{ActGr}}$ expresses the fact
``$\bfif\kern2pt(\bfix\kern2pt) \stin \rfs{Fxd}(\bfif\kern2pt)$'',
and a formula $\phi^{J_4}_{\srfs{Fxd-img}}(\bfif,\bfix\kern1pt)$
which does the same job for members of $(J_4)^{\srfs{ActGr}}§$.

The case of $\phi^{J_4}_{\srfs{Fxd-img}}$ is trivial.
Recall that
$J_4 = \setm{\pair{S}{f}}{S \stin  K_{\srfs{No-scnst}} \text{ and}\break
f \stin \rfs{Gr}(S)}$. So
$\sngltn{f} \times A \subseteq \rfs{App}^{\srfs{Act}^{\trfs{Gr}}(S)}$
and $\rfs{Fxd-img}(f) = \rfs{Fxd}(f)$.
Hence it is trivial that $\phi^{J_4}_{\srfs{Fxd-img}}(\bfif,\bfix\kern1pt)$
is the formula presented in the next observation.
\vspace{2mm}

\begin{observation}\label{o8.1}
Let
$\phi^{J_4}_{\srfs{Fxd-img}}(\bfif,\bfix\kern1pt)
\ \ \eqdf\ \ %
\rfs{App}(\bfif,\bfix\kern1pt) = \bfix$.
\index{N@fphifxd-img-j4@@$\phi^{J_4}_{\srfs{Fxd-img}}(\bfif,\bfix\kern1pt)$}%
Then for every $\pair{S}{f} \stin J_4$ and $a \stin A$:
$\rfs{Act}^{\srfs{Gr}}(S) \models \phi^{J_4}_{\srfs{Fxd-img}}[f,a]$
iff $a \stin \rfs{Fxd-img}(f)$.
\end{observation}

\noindent
{\bf Proof }
The proof is trivial.
\hfill\solidqed
\vspace{2mm}

We now turn to the proof that $\phi^{J_3}_{\srfs{Fxd-img}}$ exists.
Recall that $K_{\srfs{No-scnst}}$ denotes the class of all semigroups
$S \stin K_{\srfs{FT}}$ such that $S$ does not contain a semi-constant function,
and that $J_3$ is defined to be the class of all $\pair{S}{f}$'s
such that $S \stin K_{\srfs{No-scnst}}$ and
$f \stin \rfs{F}_{\srfs{\sbf{B}$\leq$2}}(S) \swsetminus \rfs{Gr}(S)$.
\vspace{2mm}

\subsection{The definability of 1-1 and many-to-one pairs of $\boldsymbol{f}$,
and the definability of $\bfs{Rng}(\boldsymbol{f})$ in $\boldsymbol{J_3}$}

Let $\fnn{f}{A}{A}$.
The set of ``many-to-one preimages of $f$'' - $\rfs{Mo-pre}(f)$ and
the set of ``many-to-one images of $f$'' - $\rfs{Mo-img}(f)$
are defined as follows:
$$
\rfs{Mo-pre}(f) \eqdf
\setm{a \stin A}{\abs{[a]_f} \geq 2}.
$$
That is, $\rfs{Mo-pre}(f) = \bigcup \bfB(f)$.
\index{N@mo-pre@@$\rfs{Mo-pre}(f) \eqdf \bigcup \bfB(f)$}%
$$
\rfs{Mo-img}(f) = f[\rfs{Mo-pre}(f)].
$$
Hence $\rfs{Mo-img}(f) = \setm{a \stin A}{\abs{f\inverse[\sngltn{a}]} \geq 2}$.
\index{N@mo-img@@$\rfs{Mo-img}(f) \eqdf \setm{x \stin A}{\abs{f\inverse[\sngltn{x}]} \geq 2}$}%

The set of ``1-1 preimages of $f$'' and
the set of ``1-1 images of $f$'' are defined as follows:
$$
\rfs{Oo-pre}(f) \eqdf \setm{a \stin A}{[a]_f = \sngltn{a}}
$$
\index{N@oo-pre@@$\rfs{Oo-pre}(f) \eqdf \setm{a \stin A}{[a]_f = \sngltn{a}}$}%
and
$$
\rfs{Oo-img}(f) \eqdf f[\rfs{Oo-pre}(f)].
$$
\index{N@oo-img@@$\rfs{Oo-img}(f) \eqdf \setm{f(a)}{a \stin \rfs{Oo-pre}(f)}$}%

\noindent
{\bf An overview}
\\
The goal of this subsection is to find a formula
$\phi_{\srfs{Mo-img}}(\bfif,\bfix\kern1.5pt)$
such that for every $\pair{S}{f} \stin J_3$ and $a \stin A$:
$\rfs{Act}^{\srfs{Gr}}(S) \models \phi_{\srfs{Mo-img}}[f,a]$ iff
$a \stin \rfs{Mo-img}(f)$.

Note that $\rfs{Mo-img}(f) = \rfs{Rng}(f) \swsetminus \rfs{Oo-img}(f)$.
So as intermediate steps, we shall first find formulas
$\phi_{\srfs{In-rng}}(\bfif,\bfix\kern1.5pt)$
and $\phi_{\srfs{Oo-img}}(\bfif,\bfix\kern1.5pt)$,
expressing the properties ``$\bfix \stin \rfs{Rng}(\bfif\kern2.3pt)$''
and ``$\bfix \stin \rfs{Oo-img}(\bfif\kern2.3pt)$''.
Clearly, $\phi_{\srfs{Mo-img}}$ can be taken to be
$\phi_{\srfs{In-rng}} \stwedge \stneg\phi_{\srfs{Oo-img}}$.
The formulas $\phi_{\srfs{Oo-img}}$ and $\phi_{\srfs{In-rng}}$
have their intended meaning for the class of $\pair{S}{f}$'s
such that $S \stin K_{\srfs{FT}}$ and $\abs{\rfs{Oo-img}(f)} \geq 3$,
and $J_3$ is a subclass of this class.

\begin{observation}\label{o3.1}\label{o8.2}
There is a formula $\phi_{\srfs{Mo-pre}}(\bfif,\bfix\kern1.5pt)$
\index{N@fphimopre@@$\phi_{\srfs{Mo-pre}}(\bfif,\bfix\kern1.5pt)$}%
such that for every\break
$S \stin K_{\srfs{FT}}$, $f \stin S$ and $a \stin A$:
$$\mbox{
$\rfs{Act}^{\srfs{Gr}}(S) \models \phi_{\srfs{Mo-pre}}[f,a]$
\ \ \ iff\ \ \ \ %
$a \stin \rfs{Mo-pre}(f)$.
}$$
\end{observation}

\noindent
{\bf Proof }
%
%
The formula $\phi_{\srfs{sim}}(\bfif,\bfix,\bfiy\kern1pt)$
which was defined in Observation~\ref{o1.7} expresses the fact that
$\bfix \sim_{\sbfif} \bfiy$.
So 
$\phi_{\srfs{Mo-pre}}(\bfif,\bfix\kern1.5pt) \eqdf
\exists \bfiy\kern1pt((\bfiy \neq \bfix\kern1pt) \stwedge
\phi_{\srfs{sim}}(\bfif,\bfix,\bfiy\kern1pt))$ is as required.
\hfill\solidqed
\vspace{2mm}

\begin{observation}\label{o3.2}\label{o8.3}
{\rm(a)}
Let $S \stin K_{\srfs{No-scnst}}$
and  $f \stin \rfs{F}_{\srfs{\sbf{B}$\leq$2}}(S)$.
Then either {\rm(i)} $f \stin \rfs{Sym}(A)$
and $\setm{x \stin A}{f(x) \neq x}$ is finite,
or {\rm(ii)}  $A \swsetminus (\rfs{Mo-pre}(f) \stcup \rfs{Idp}(f))$ is infinite.
In the former case, $f \stin \rfs{Gr}(S)$.

{\rm(b)}
For every $\pair{S}{f} \stin J_3$, $f$ has infinitely many $f$-1-1 pairs.
\end{observation}

\noindent
{\bf Proof }
(a)
We distinguish between two cases.

{\bf Case 1:} $f$ is not $\onetoone$. Then $f$ does not fulfill (i).
Suppose by contradiction that $f$ does not fulfill (ii). That is,
$R \eqdf A \swsetminus (\rfs{Mo-pre}(f) \stcup \rfs{Idp}(f))$ is finite.
Since $\abs{\bfs{B}(f)} \leq 2$, there are $x,y \stin A$
such that $\rfs{Mo-pre}(f) \subseteq [x]_f \stcup [y]_f$.
Let $p = f(x)$ and $q = f(y)$.
Obviously, $A \swsetminus \rfs{Idp}(f) = \rfs{Mo-pre}(f) \stcup R$.
Hence $f[A \swsetminus \rfs{Idp}(f)] = \dbltn{p}{q} \stcup f[R]$.
This last set is finite.
By Observation~\ref{o7.2},\break
$S$ contains a semi-constant function,
contradicting the fact that $S \stin K_{\srfs{No-Scnst}}$.
So $A \swsetminus (\rfs{Mo-pre}(f) \stcup \rfs{Idp}(f))$ is infinite.

{\bf Case 2:} $f$ is $\onetoone$.
Assume that $f$ does not fulfill (ii), and we show that $f$ fulfills (i).
The fact that $f$ does not fulfill (ii) means that
$A \swsetminus (\rfs{Mo-pre}(f) \stcup \rfs{Idp}(f))$ is finite.
Clearly,
$A \swsetminus (\rfs{Mo-pre}(f) \stcup \rfs{Idp}(f)) =
A \swsetminus \rfs{Idp}(f)$.
So $A \swsetminus \rfs{Idp}(f)$ is finite.
This last fact together with the fact that $f$ is 1-1
implies that $f \stin \rfs{Sym}(A)$.
Also, $\setm{x \stin A}{f(x) \neq x}$ is finite.

Clearly, in Case 2, $f\inverse \stin S$. So $f \stin \rfs{Gr}(S)$.
\vspace{2mm}

(b) Part (b) follows trivially from Part (a), and is brought here just
for the sake of later reference.
\hfill\solidqed
\vspace{2mm}

\noindent
{\bf Remark }
The formulas of the form $\phi(\bfif,\ldots)$ defined in this section
have their intended meaning only for $f$'s which have at least three 1-1 pairs.
However, since for every $\pair{S}{f} \stin J_3$,
$f$ has infinitely many 1-1 pairs,
when applied to members of $J_3$,
these formulas will indeed express the property we wish them to express. 
\vspace{2mm}

{\thickmuskip=2mu \medmuskip=1mu \thinmuskip=1mu
Let $\fnn{f}{A}{A}$ and $p,q \stin A$.
We call $\pair{p}{q}$ an \underline{$f$-1-1 pair}, 
if $f\inverse[\sngltn{q}] = \sngltn{p}$.}
\index{1-1 pair@@$f$-1-1 pair. $\pair{p}{q}$ is an $f$-1-1 pair if $f\inverse[\sngltn{q}] = \sngltn{p}$}%
In this situation $p$ is called an \underline{$f$-1-1 preimage},
and $q$ is called an \underline{$f$-1-1 image}.%
\index{1-1 image@@$f$-1-1 image. If $\pair{p}{q}$ is an $f$-1-1 pair, then $q$ is called an $f$-1-1 image}%
\index{1-1 preimage@@$f$-1-1 preimage. If $\pair{p}{q}$ is an $f$-1-1 pair, then $p$ is called an $f$-1-1\\\indent preimage}%

We call $\pair{p}{q}$ a \underline{$f$-simple pair}, 
if $p \neq q$ and $f\inverse[\sngltn{q}] = \sngltn{p}$.
In this situation $p$ is called an \underline{$f$-simple preimage},
and $q$ is called an \underline{$f$-simple image}.
Note if $\pair{p}{q}$ is an $f$-1-1 non-$f$-simple pair,
then $p = q \stin \rfs{Idp}(f)$.
\index{simple pair@@$f$-simple pair. An $f$-1-1 pair $\pair{p}{q}$ such that $p \neq q$}%
\index{simple preimage@@$f$-simple preimage. The $f$-1-1 preimage of a simple pair}%
\index{simple image@@$f$-simple image. The $f$-1-1 image of a simple pair}%

\begin{prop}\label{p3.3}\label{p8.4}
{\rm(a)}
There is a formula $\phi_{\srfs{Oo-pre}}(\bfif,\bfix\kern1.5pt)$
such that for every $S \stin K_{\srfs{FT}}$,
$f \stin S$ and $a \stin A$:
$\rfs{Act}^{\srfs{Gr}}(S) \models \phi_{\srfs{Oo-pre}}[f,a]$ iff
$a$ is an $f$-1-1 preimage.

{\rm(b)}
There is a formula $\phi_{\srfs{1-1$\geq$3}}(\bfif\kern2.3pt)$
such that for every $S \stin K_{\srfs{FT}}$ and $f \stin S$:
$\rfs{Act}^{\srfs{Gr}}(S) \models \phi_{\srfs{1-1$\geq$3}}[f]$ iff
$f$ has at least three $f$-1-1 pairs.

{\rm(c)}
Let
$$
\phi^1_{\srfs{1-1}}
(\bfif,\bfix\fsub{1},\bfiu\fsub{1},\bfix\fsub{2},\bfiu\fsub{2}\kern1.5pt) \eqdf
(\bfix\fsub{1} \stnot\sim_{\sbfif} \bfix\fsub{2}\kern1.5pt) \stwedge
(\bfiu\fsub{1} \neq \bfiu\fsub{2}\kern1.5pt) \stwedge
\left(\rule{0pt}{12pt}\rvpair{\bfiu\fsub{1}}{\bfiu\fsub{2}\kern1.5pt} \bfif
\rvpair{\bfix\fsub{1}}{\bfix\fsub{2}\kern1.5pt} = \bfif\kern2.3pt\right).
$$
Then for every $S \stin K_{\srfs{FT}}$, $f \stin S$ and
$x_1,u_1,x_2,u_2 \stin A$, clauses {\rm(1)} and {\rm(2)} below are equivalent:
\begin{itemize}
\addtolength{\parskip}{-11pt}
\addtolength{\itemsep}{06pt}
\item[{\rm(1)}]
$\rfs{Act}^{\srfs{Gr}}(S) \models \phi^1_{\srfs{1-1}}[f,x_1,u_1,x_2,u_2]$.
\item[{\rm(2)}]
\begin{list}{}
{\setlength{\leftmargin}{39pt}
\setlength{\labelsep}{05pt}
\setlength{\labelwidth}{25pt}
\setlength{\itemindent}{-00pt}
\addtolength{\topsep}{03pt}
\addtolength{\parskip}{02pt}
\addtolength{\itemsep}{-01pt}
}
\item[{\rm(O1)}]
$u_1 \neq u_2$, and
\item[{\rm(O2)}]
\begin{list}{}
{\setlength{\leftmargin}{50pt}
\setlength{\labelsep}{05pt}
\setlength{\labelwidth}{45pt}
\setlength{\itemindent}{-00pt}
\addtolength{\topsep}{03pt}
\addtolength{\parskip}{02pt}
\addtolength{\itemsep}{-01pt}
}
\item[either {\rm(i)}]
$f\inverse[\sngltn{u_1}] = \sngltn{x_1}$
and $f\inverse[\sngltn{u_2}] = \sngltn{x_2}$,
\item[or {\rm(ii)}]
$f\inverse[\sngltn{u_1}] = \sngltn{x_2}$
and $f\inverse[\sngltn{u_2}] = \sngltn{x_1}$.
\end{list}
\end{list}
\vspace{-05.7pt}
\end{itemize}
\vspace{2mm}

{\rm(d)}
Let
$\phi^2_{\srfs{1-1}}
(\bfif,\bfix\fsub{1},\bfiu\fsub{1},\bfix\fsub{2},\bfiu\fsub{2},
\bfix\fsub{3},\bfiu\fsub{3}\kern1.5pt)$
be the formula
\begin{equation*}
(\bfiu\fsub{2} \neq \bfiu\fsub{3}) \kern2pt\stwedge\kern2pt
\mbox{$\bigwedge$}_{\kern2pti \stin \dbltn{2}{3}}
\kern2pt\phi^1_{\srfs{1-1}}
(\bfif,\bfix\fsup{1},\bfiu\fsup{1},\bfix\fsup{i},\bfiu\fsup{i}\kern1.5pt).
\end{equation*}
Then for every $S \stin K_{\srfs{FT}}$, $f \stin S$ and
$x_1,u_1,x_2,u_2,x_3,u_3 \stin A$, the following are equivalent:
\begin{itemize}
\addtolength{\parskip}{-11pt}
\addtolength{\itemsep}{06pt}
\item[{\rm(1)}]
$\rfs{Act}^{\srfs{Gr}}(S) \models \phi^2_{\srfs{1-1}}[f,x_1,u_1,x_2,u_2,x_3,u_3]$
{\thickmuskip=2mu \medmuskip=1mu \thinmuskip=1mu
\item[{\rm(2)}] $u_1,u_2,u_3$ are pairwise distinct,
and for every $i = 1,2,3$, $f\inverse[\sngltn{u_i}] = \sngltn{x_i}$.
}%
\vspace{-8mm}%
\end{itemize}%

{\rm(e)} 
There is a formula $\phi_{\srfs{Oo-img}}(\bfif,\bfix\kern1.5pt)$
such that for every $S \stin K_{\srfs{FT}}$,
every $f \stin S$
which has at least two $f$-1-1 pairs and every $a \stin A$:
$\rfs{Act}^{\srfs{Gr}}(S) \models \phi_{\srfs{Oo-img}}[f,a]$
iff
$a$ is an $f$-1-1 image.

{\rm(f)} 
There are formulas $\phi_{\srfs{Oo-pair}}(\bfif,\bfix,\bfiu\kern1.5pt)$,
$\phi_{\srfs{smpl-pair}}(\bfif,\bfix,\bfiu\kern1.5pt)$,
$\phi_{\srfs{Idp}}(\bfif,\bfix\kern1.5pt)$,\break
$\phi_{\srfs{Oo-pre}}(\bfif,\bfix\kern1.5pt)$,
$\phi_{\srfs{Oo-img}}(\bfif,\bfiu\kern1.5pt)$,
$\phi_{\srfs{smpl-pre}}(\bfif,\bfix\kern1.5pt)$
and $\phi_{\srfs{smpl-img}}(\bfif,\bfiu\kern1.5pt)$
such that for every $S \stin K_{\srfs{FT}}$,
every $f \stin S$ which has at least three $f$-1-1 pairs and every $a,b \stin A$:
\begin{itemize}
\addtolength{\parskip}{-11pt}
\addtolength{\itemsep}{06pt}
\item[{\rm(1)}]
$\pair{a}{b}$ is an $f$-1-1 pair iff
$\rfs{Act}^{\srfs{Gr}}(S) \models \phi_{\srfs{Oo-pair}}[f,a,b]$,
\item[{\rm(2)}]
$\pair{a}{b}$ is an $f$-simple pair iff
$\rfs{Act}^{\srfs{Gr}}(S) \models \phi_{\srfs{smpl-pair}}[f,a,b]$,
\item[{\rm(3)}]
$a \stin \rfs{Idp}(f)$ iff
$\rfs{Act}^{\srfs{Gr}}(S) \models \phi_{\srfs{Idp}}[f,a]$,
\item[{\rm(4)}]
$a \stin \rfs{Oo-pre}(f)$ iff
$\rfs{Act}^{\srfs{Gr}}(S) \models \phi_{\srfs{Oo-pre}}[f,a]$,
\item[{\rm(5)}]
$a \stin \rfs{Oo-img}(f)$ iff
$\rfs{Act}^{\srfs{Gr}}(S) \models \phi_{\srfs{Oo-img}}[f,a]$,
\item[{\rm(6)}]
$a$ is an $f$-simple preimage iff
$\rfs{Act}^{\srfs{Gr}}(S) \models \phi_{\srfs{smpl-pre}}[f,a]$, and
\item[{\rm(7)}]
$a$ is an $f$-simple image iff
$\rfs{Act}^{\srfs{Gr}}(S) \models \phi_{\srfs{smpl-img}}[f,a]$.
\index{N@fphi11pair@@$\phi_{\srfs{Oo-pair}}(\bfif,\bfix,\bfiu\kern1.5pt)$}%
\index{N@fphismpl-pair@@$\phi_{\srfs{smpl-pair}}(\bfif,\bfix,\bfiu\kern1.5pt)$}%
\vspace{-05.7pt}
\end{itemize}

%
\vspace{3mm}
{\rm(g)}

\kern-08.4mm

\begin{list}{}
{\setlength{\leftmargin}{65pt}
\setlength{\labelsep}{05pt}
\setlength{\labelwidth}{51pt}
\setlength{\itemindent}{-00pt}
\addtolength{\topsep}{-04pt}
\addtolength{\parskip}{-02pt}
\addtolength{\itemsep}{-08pt}
}
\item[{\rm(1)}]
Let $\fnn{f}{A}{A}$ and $x \stin \rfs{Rng}(f)$.
Then for every $y \stin A \swsetminus \sngltn{x}$,\break
$\rvpair{x}{y} f \neq~f$.
\vspace{0.5mm}
\item[{\rm(2)}] Let $\fnn{f}{A}{A}$ and $x,y \stin A \swsetminus \rfs{Rng}(f)$.
Then $\rvpair{x}{y} f = f$.
\vspace{10.0pt}
\end{list}

{\rm(h)}

\kern-08.7mm

\begin{list}{}
{\setlength{\leftmargin}{65pt}
\setlength{\labelsep}{05pt}
\setlength{\labelwidth}{51pt}
\setlength{\itemindent}{-00pt}
\addtolength{\topsep}{-04pt}
\addtolength{\parskip}{-02pt}
\addtolength{\itemsep}{-08pt}
}
\item[{\rm(1)}]
Let $\phi_{\srfs{notrng$\geq$2}}(\bfif\kern2pt)$ be the formula which says:
\vspace{-1.5mm}
\begin{list}{}
{\setlength{\leftmargin}{39pt}
\setlength{\labelsep}{05pt}
\setlength{\labelwidth}{25pt}
\setlength{\itemindent}{-00pt}
\addtolength{\topsep}{-04pt}
\addtolength{\parskip}{-02pt}
\addtolength{\itemsep}{-10pt}
}
\item[$\bullet$]
There is a transposition $\pi$ such that $\pi f = f$.
\vspace{-05.0pt}
\end{list}
Then for every $S \stin K_{\srfs{FT}}$ and $f \stin S$,
\vspace{-4mm}
$$\mbox{
$\rfs{Act}^{\srfs{Gr}}(S) \models \phi_{\srfs{notrng$\geq$2}}[f]$\ \ \ \ %
iff\ \ \ \ %
$\abs{A \swsetminus \rfs{Rng}(A)} \geq 2$.
}
\vspace{-2mm}
$$
\item[{\rm(2)}]
Let $\phi_{\srfs{In-rng}}^1(\bfif,\bfix\kern1pt)$ be the formula which says:
\begin{list}{}
{\setlength{\leftmargin}{39pt}
\setlength{\labelsep}{05pt}
\setlength{\labelwidth}{25pt}
\setlength{\itemindent}{-00pt}
\addtolength{\topsep}{-04pt}
\addtolength{\parskip}{-02pt}
\addtolength{\itemsep}{-05pt}
}
\item[]
For every $\bfiy \neq \bfix$, $\rvpair{\bfix}{\bfiy\kern1pt} \bfif \neq \bfif$.
\index{N@fphiin-rng@@$\phi_{\srfs{In-rng}}^1(\bfif,\bfix\kern1pt)$}%
\vspace{-05.0pt}
\end{list}
Then for every $S \stin K_{\srfs{FT}}$ and every $f \stin S$ such that\\
$\abs{A \swsetminus \rfs{Rng}(f)} \geq 2$,
the following holds:
For every $a \stin A$:\\
$\rfs{Act}^{\srfs{Gr}}(S) \models \phi_{\srfs{In-rng}}^1[f,a]$%
\ \ \ iff\ \ \ %
 $a \stin \rfs{Rng}(f)$.
\vspace{-02.0pt}
\end{list}
\end{prop}

\noindent
{\bf Proof }
(a)
Clearly, the formula
$\phi_{\srfs{Oo-pre}}(\bfif,\bfix\kern1.5pt) \eqdf
\forall \bfiy\kern1pt(\phi_{\srfs{sim}}(\bfif,\bfix,\bfiy\kern1pt) \rightarrow
(\bfix = \bfiy\kern1pt))$
is as required.
\vspace{2mm}

(b)
Clearly, the formula
\vspace{-2mm}
$$\mbox{
$\phi_{\srfs{1-1$\geq$3}}(\bfif\kern2.3pt) \eqdf
\exists \bfix\fsub{1} \bfix\fsub{2} \bfix\fsub{3}
\left(\bigwedge_{i = 1}^3
\phi_{\srfs{Oo-pre}}(\bfif,\bfix\fsub{i}\kern1.5pt) \stwedge
\bigwedge_{1 \leq i < j \leq 3} (\bfix\fsub{i} \neq \bfix\fsub{j})
\right)$
}
\vspace{-2mm}
$$
is as required.
\vspace{2mm}

(c)
{\bf Proof of (2) $\bm{\Rightarrow}$ (1) }
Let $S \stin K_{\srfs{FT}}$, $f \stin S$ and
$x_1,u_1,x_2,u_2 \stin A$.\break
Suppose that $u_1 \neq u_2$,
$f\inverse[\sngltn{u_1}] = \sngltn{x_1}$
and $f\inverse[\sngltn{u_2}] = \sngltn{x_2}$.
We show that
$\rfs{Act}^{\srfs{Gr}}(S) \models \phi^1_{\srfs{1-1}}[f,x_1,u_1,x_2,u_2]$.
Clearly, $u_1 \neq u_2$ and $x_1 \stnot\sim_f x_2$.
It remains to verify that
$\rvpair{u\fsub{1}}{u\fsub{2}\kern1.5pt} f
\rvpair{x\fsub{1}}{x\fsub{2}\kern1.5pt} = f\kern1pt$.
Let $c \stin A$. Assume first that $c \neq x_1,x_2$.
Then
$$
\rvpair{u\fsub{1}}{u\fsub{2}\kern1.5pt} f
\rvpair{x\fsub{1}}{x\fsub{2}\kern1.5pt}(c) =
\rvpair{u\fsub{1}}{u\fsub{2}\kern1.5pt} f(c).
$$
For every $i = 1,2$, $x_i$ is the unique element $x$ such that $f(x) = u_i$.
Hence $f(c) \neq u_1,u_2$.
So
$$
\rvpair{u\fsub{1}}{u\fsub{2}\kern1.5pt} f(c) = f(c).
$$
Suppose that $c = x_1$.
Then
\vspace{-2mm}
$$
\rvpair{u\fsub{1}}{u\fsub{2}\kern1.5pt} f
\rvpair{x\fsub{1}}{x\fsub{2}\kern1.5pt}(x_1) =
\rvpair{u\fsub{1}}{u\fsub{2}\kern1.5pt} f(x_2) =
\rvpair{u\fsub{1}}{u\fsub{2}\kern1.5pt}(u_2) = u_1 = f(x_1).
\vspace{-2mm}
$$
The same computation holds for $x_2$.
So
$\rvpair{u\fsub{1}}{u\fsub{2}\kern1.5pt} f
\rvpair{x\fsub{1}}{x\fsub{2}\kern1.5pt} = f\kern1pt$,
and hence
$\rfs{Act}^{\srfs{Gr}}(S) \models \phi^1_{\srfs{1-1}}[f,x_1,u_1,x_2,u_2]$.

The same argument applies to the case that
$f\inverse[\sngltn{u_1}] = \sngltn{x_2}$
and\break
$f\inverse[\sngltn{u_2}] = \sngltn{x_1}$.
So if $f,x_1,x_2,u_1,u_2$ fulfill (O1) and (O2),
then $\rfs{Act}^{\srfs{Gr}}(S) \models \phi^1_{\srfs{1-1}}[f,x_1,u_1,x_2,u_2]$.
\vspace{2mm}

{\bf Proof of (1) $\bm{\Rightarrow}$ (2) }
Assume now that
$\rfs{Act}^{\srfs{Gr}}(S) \models \phi^1_{\srfs{1-1}}[f,x_1,u_1,x_2,u_2]$.
We show that $f,x_1, u_1,x_2,u_2$ fulfill Clause (2) of Part (c).

Clause (2) of Part (c) is equivalent to the following statement.
\vspace{-3mm}
\begin{list}{}
{\setlength{\leftmargin}{39pt}
\setlength{\labelsep}{05pt}
\setlength{\labelwidth}{25pt}
\setlength{\itemindent}{-00pt}
\addtolength{\topsep}{03pt}
\addtolength{\parskip}{02pt}
\addtolength{\itemsep}{-01pt}
}
\item[{\rm(P1)}]
$u_1 \neq u_2$,
\item[{\rm(P2)}]
\begin{list}{}
{\setlength{\leftmargin}{50pt}
\setlength{\labelsep}{05pt}
\setlength{\labelwidth}{45pt}
\setlength{\itemindent}{-00pt}
\addtolength{\topsep}{03pt}
\addtolength{\parskip}{02pt}
\addtolength{\itemsep}{-01pt}
}
\item[either {\rm(i)}]
$f(x_1) = u_1$ and $f(x_2) = u_2$
\item[or {\rm(ii)}]
$f(x_1) = u_2$ and $f(x_2) = u_1$,
\end{list}
and
\item[{\rm(P3)}]
$f\inverse[\dbltn{u_1}{u_2}] = \dbltn{x_1}{x_2}$.
\end{list}
Clearly, $u_1 \neq u_2$ because this is part of what $\phi^1_{\srfs{1-1}}$ says.
Hence (P1) holds.

Now we show that $f(x_1) \stin \dbltn{u_1}{u_2}$.
Suppose by contradiction that\break
$f(x_1) \stnot\stin \dbltn{u_1}{u_2}$.
Then
$$
f(x_2) = 
\rvpair{u\fsub{1}}{u\fsub{2}\kern1.5pt} f
\rvpair{x\fsub{1}}{x\fsub{2}\kern1.5pt}(x_2) =
\rvpair{u\fsub{1}}{u\fsub{2}\kern1.5pt} f(x_1) = f(x_1).
$$
Hence $x_1 \sim_f x_2$. However, this contradicts the fact that
$\rfs{Act}^{\srfs{Gr}}(S) \models\break
\phi^1_{\srfs{1-1}}[f,x_1,u_1,x_2,u_2]$.
So indeed, $f(x_1) \stin \dbltn{u_1}{u_2}$.

Similarly, $f(x_2) \stin \dbltn{u_1}{u_2}$.
Since by $\phi^1_{\srfs{1-1}}[f,x_1,u_1,x_2,u_2]$,
$x_1 \stnot\sim_f x_2$,
either (i) $f(x_1) = u_1$ and  $f(x_2) = u_2$,
or (ii) $f(x_1) = u_2$ and  $f(x_2) = u_1$.
We have proved (P2).

Let $c \neq x_1,x_2$, and we show that $f(c) \neq u_1,u_2$.
Suppose by contradiction that $f(c) = u_1$. Then
$$
u_1 = f(c) =
\rvpair{u\fsub{1}}{u\fsub{2}\kern1.5pt} f
\rvpair{x\fsub{1}}{x\fsub{2}\kern1.5pt}(c) =
\rvpair{u\fsub{1}}{u\fsub{2}\kern1.5pt} f(c) =
\rvpair{u\fsub{1}}{u\fsub{2}\kern1.5pt}(u_1) = u_2.
$$
But the fact that $\phi^1_{\srfs{1-1}}(f,x_1,u_1,x_2,u_2)$ holds,
implies that $u_1 \neq u_2$.
A contradiction. So $f(c) \neq u_1$.

Similarly, $f(c) \neq u_2$.
This shows that $f\inverse[\dbltn{u_1}{u_2}] = \dbltn{x_1}{x_2}$.
That is, (P3) holds.
\vspace{2mm}

(d)
Part (d) follows quite trivially from Part (c).
That Clause (2) implies Clause (1) is immediate.
(It follows, of course, from Part (c).)

{\thickmuskip=0mu \medmuskip=0mu \thinmuskip=0mu
We prove that $(1) \Rightarrow (2)$.
Assume that
$\rfs{Act}^{\srfs{Gr}}\kern-2pt(S) \models
\phi^2_{\srfs{1-1}}[f,x_1,u_1,x_2,u_2,x_3,u_3]$.
}%
Since $u_1,u_2,u_3$ satisfy the first conjunct of $\phi^2_{\srfs{1-1}}$,
$u_2 \neq u_3$, and the fact that $u_1,u_2,u_3$ satisfy
$\phi^1_{\srfs{1-1}}
(\bfif,\bfix\fsup{1},\bfiu\fsup{1},\bfix\fsup{i},\bfiu\fsup{i}\kern1.5pt)$,
implies that for $i = 2,3$, $u_1 \neq u_i$.
Hence $u_1,u_2,u_3$ are pairwise distinct.

By Part (c), since
$\rfs{Act}^{\srfs{Gr}}\kern-2pt(S) \models
\phi^1_{\srfs{1-1}}[f,x_1,u_1,x_2,u_2]$,
$f(x_1) \stin \dbltn{u_1}{u_2}$,
and since
$\rfs{Act}^{\srfs{Gr}}\kern-2pt(S) \models
\phi^1_{\srfs{1-1}}[f,x_1,u_1,x_3,u_3]$,
$f(x_1) \stin \dbltn{u_1}{u_3}$.
It follows that\break
$(*)$ $f(x_1) = u_1$.
Part (c) together with Fact $(*)$ and the fact that
$\rfs{Act}^{\srfs{Gr}}\kern-2pt(S) \models
\phi^1_{\srfs{1-1}}[f,x_1,u_1,x_i,u_i]$
imply that for $i = 2,3$, $f(x_i) = u_i$.
We have shown that 
$f,x_1,u_1,x_2,u_2$ fulfill Clause (2) of Part (d).
\vspace{2mm}

(e)
Let $\phi_{\srfs{Oo-img}}(\bfif,\bfiu\kern1.5pt)$ be the formula which says: 
\begin{itemize}
\addtolength{\parskip}{-11pt}
\addtolength{\itemsep}{06pt}
\item
There are $\bfix,\bfiy$ and $\bfiv$
such that $\phi^1_{\srfs{1-1}}(\bfif,\bfix,\bfiu,\bfiy,\bfiv\kern1.5pt)$.
\vspace{-05.7pt}
\end{itemize}
\index{N@fphioo-img@@$\phi_{\srfs{Oo-img}}(\bfif,\bfix\kern1.5pt)$}%
We prove that for every $S \stin K_{\srfs{FT}}$,
every $f \stin S$ which has at least two $f$-1-1 pairs and every $a \stin A$:
$\rfs{Act}^{\srfs{Gr}}(S) \models \phi_{\srfs{Oo-img}}[f,a]$ iff
$a$ is an $f$-1-1 image.

{\bf Proof of $\bm{\Leftarrow}$.}
Let $u$ be an $f$-1-1 image.
Let $x$ be such that $f(x) = u$.
Also, let $\pair{y}{v}$ be an $f$-1-1 pair different from $\pair{x}{u}$.
So $u \neq v$.
By Part (c) of this proposition,
$\rfs{Act}^{\srfs{Gr}}(S) \models \phi^1_{\srfs{1-1}}[f,x,u,y,v]$.

{\bf Proof of $\bm{\Rightarrow}$.}
Suppose that $\rfs{Act}^{\srfs{Gr}}(S) \models \phi_{\srfs{Oo-img}}[f,u]$,
and let $x,y,v$ be the elements whose existence is assured by
$\phi_{\srfs{Oo-img}}$.
So by Part~(c) of this proposition (O2) holds.
Namely, 
either (i) $f\inverse[\sngltn{u}] = \sngltn{x}$
and $f\inverse[\sngltn{v}] = \sngltn{y}$,
or\break
(ii) $f\inverse[\sngltn{u}] = \sngltn{y}$
and $f\inverse[\sngltn{v}] = \sngltn{x}$.
In either case, $u$ is an $f$-1-1 image.
\vspace{2mm}

(f)
Let
\begin{equation}
\tag*{\rule{4mm}{0pt}(1)}
\phi_{\srfs{Oo-pair}}(\bfif,\bfix,\bfiu\kern1.5pt) \eqdf
(\exists \bfiy,\bfiv,\bfiz,\bfiw\kern1pt)\kern2pt
\phi^2_{\srfs{1-1}}(\bfif,\bfix,\bfiu,\bfiy,\bfiv,\bfiz,\bfiw\kern1.5pt),
\vspace{-5.5mm}
\end{equation}
\begin{equation}
\tag*{\rule{4mm}{0pt}(2)}
\phi_{\srfs{smpl-pair}}(\bfif,\bfix,\bfiu\kern1.5pt) \eqdf
\phi_{\srfs{Oo-pair}}(\bfif,\bfix,\bfiu\kern1.5pt) \stwedge
(\bfiu \neq\bfix\kern1pt),
\vspace{-5.5mm}
\end{equation}
\begin{equation}
\tag*{\rule{4mm}{0pt}(3)}
\phi_{\srfs{Idp}}(\bfif,\bfix\kern1.5pt) \eqdf
\phi_{\srfs{Oo-pair}}(\bfif,\bfix,\bfix\kern1.5pt),
\vspace{-5.5mm}
\end{equation}
\begin{equation}
\tag*{\rule{4mm}{0pt}(4)}
\phi_{\srfs{Oo-pre}}(\bfif,\bfix\kern1.5pt) \eqdf
\exists \bfiu\kern1pt\phi_{\srfs{Oo-pair}}(\bfif,\bfix,\bfiu\kern1.5pt),
\vspace{-5.5mm}
\end{equation}
\begin{equation}
\tag*{\rule{4mm}{0pt}(5)}
\phi_{\srfs{Oo-img}}(\bfif,\bfiu\kern1.5pt) \eqdf
\exists \bfix \kern1pt\phi_{\srfs{Oo-pair}}(\bfif,\bfix,\bfiu\kern1.5pt),
\vspace{-5.5mm}
\end{equation}
\begin{equation}
\tag*{\rule{4mm}{0pt}(6)}
\phi_{\srfs{smpl-pre}}(\bfif,\bfix\kern1.5pt) \eqdf
\exists \bfiu\kern2pt \phi_{\srfs{smpl-pair}}(\bfif,\bfix,\bfiu\kern1.5pt)
\vspace{-2mm}
\end{equation}
and
\vspace{-2mm}
\begin{equation}
\tag*{\rule{4mm}{0pt}(7)}
\phi_{\srfs{smpl-img}}(\bfif,\bfiu\kern1.5pt) \eqdf
\exists \bfix\kern2.5pt \phi_{\srfs{smpl-pair}}(\bfif,\bfix,\bfiu\kern1.5pt).
\end{equation}
That the formula in (1) is as required follows from Part (d).
The proof that formulas (2) - (7) are as required is trivial.
\vspace{2mm}

(g) The proof of Part (g) is trivial.
\vspace{2mm}

(h) Part (h) is a trivial corollary of Part (g).
\vspace{6mm}
\hfill\solidqed

We now turn to finding the formula $\phi{\srfs{In-rng}}$.
Let $\phi_{\srfs{Not-in-rng}}(\bfif,\bfix\kern2pt)$
be the formula which says:
\begin{list}{}
{\setlength{\leftmargin}{39pt}
\setlength{\labelsep}{05pt}
\setlength{\labelwidth}{25pt}
\setlength{\itemindent}{-00pt}
\addtolength{\topsep}{-04pt}
\addtolength{\parskip}{-02pt}
\addtolength{\itemsep}{-05pt}
}
\item[] \underline{either}
\item[(N1)]
There is $\bfiy$
such that
(I) $\phi_{\srfs{smpl-pair}}(\bfif,\bfix,\bfiy\kern1pt)$,
\underline{and}
(II) $\rvpair{\bfix}{\bfiy\kern1pt} \bfif\Fsup{2} = \bfif\Fsup{2}$,
\item[] \underline{or}
\item[(N2)]
\begin{list}{}
{\setlength{\leftmargin}{26pt}
\setlength{\labelsep}{05pt}
\setlength{\labelwidth}{25pt}
\setlength{\itemindent}{-00pt}
\addtolength{\topsep}{03pt}
\addtolength{\parskip}{02pt}
\addtolength{\itemsep}{-01pt}
}
\item[(III)]
There is no $\bfiy$
such that $\phi_{\srfs{smpl-pair}}(\bfif,\bfix,\bfiy\kern1pt)$,
\underline{and}
\item[(IV)]
There are $\bfiy,\bfiz \neq \bfix$
such that $\phi_{\srfs{Oo-pair}}(\bfif,\bfiy,\bfiz\kern1pt)$
and\\
$\rvpair{\bfix}{\bfiz\kern1pt} (\bfif \rvpair{\bfix}{\bfiy\kern1pt}\kern-2pt)^2
= (\bfif \rvpair{\bfix}{\bfiy\kern1pt}\kern-2pt)^2$,
\end{list}
\end{list}
\index{N@fphinot-in-rng@@$\phi_{\srfs{Not-in-rng}}(\bfif,\bfix\kern2pt)$. For $f$'s with at least three 1-1 pairs it says that\\\indent$\bfix \stnot\stin \rfs{Rng}(\bfif\kern2.6pt)$}%

\begin{prop}\label{p3.5}\label{p8.5}
{\rm(a)}
Let $S \stin K_{\srfs{FT}}$, $f \stin S$ be such that $f$ has at least three
$f$-1-1 pairs,
and $b \stin A$. Then the following are equivalent:
\begin{list}{}
{\setlength{\leftmargin}{39pt}
\setlength{\labelsep}{05pt}
\setlength{\labelwidth}{25pt}
\setlength{\itemindent}{-00pt}
\addtolength{\topsep}{-04pt}
\addtolength{\parskip}{-02pt}
\addtolength{\itemsep}{-05pt}
}
\item[{\rm(R1)}]
$b \stnot\stin \rfs{Rng}(f)$.
\item[{\rm(R2)}]
$\rfs{Act}^{\srfs{Gr}}(S) \models \phi_{\srfs{Not-in-rng}}[f,b]$.
\vspace{-02.0pt}
\end{list}
Remark: The existence of at least three $f$-1-1 pairs is needed
in order to ensure:
{\rm(i)} the equivalence between the formula $\phi_{\srfs{smpl-pair}}$
and the $f$-simplicity of a pair; and
{\rm(ii)} the equivalence between the formula $\phi_{\srfs{Oo-pair}}$
and the $f$-1-1-ness of a pair.
See Proposition~\ref{p3.3}(f)(2).

{\rm(b)}
Let
$\phi_{\srfs{In-rng}}(\bfif,\bfix\kern1.5pt) \eqdf
\stneg\phi_{\srfs{Not-in-rng}}(\bfif,\bfix\kern1.5pt)$.
Then for every $S \stin K_{\srfs{FT}}$,
every $f \stin S$ which has at least three $f$-1-1 pairs and every $b \stin A$:
$\rfs{Act}^{\srfs{Gr}}(S) \models \phi_{\srfs{In-rng}}[f,b]$
iff $b \stin \rfs{Rng}(f)$.
\index{N@fphiin-rng@@$\phi_{\srfs{In-rng}}(\bfif,\bfix\kern1.5pt)$. For $f$'s with at least three 1-1 pairs it says that\\\indent$\bfix \stin \rfs{Rng}(\bfif\kern2.6pt)$}%
\end{prop}

\noindent
{\bf Proof:}
(a)
{\bf Proof of $\bfs{(R1)} \boldsymbol{\Rightarrow} \bfs{(R2)}$ }
Suppose that $b \stnot\stin \rfs{Rng}(f)$.
So $b \neq f(b)$. Hence either
(i) $\pair{b}{f(b)}$ is $f$-simple,
or (ii) $\abs{[b]_f} \geq 2$.
We show that if (i) happens,
then $f,b$ fulfill (N1) of $\phi_{\srfs{Not-in-rng}}$,
and if (ii) happens,
then $f,b$ fulfill (N2) of $\phi_{\srfs{Not-in-rng}}$.

(i) Assume that $\pair{b}{f(b)}$ is $f$-simple.
It is trivial to check that $f(b) \stnot\stin \rfs{Rng}(f^2)$.
Obviously, $b \stnot\stin \rfs{Rng}(f^2)$.
This implies that $\rvpair{b}{f(b)} f^2 = f^2$.
Hence choosing $f(b)$ to replace $\bfy$, shows that
(N1) holds for $f$ and $b$.

(ii) Suppose that $\abs{[b]_f} \geq 2$.
So trivially, (III) in (N2) is fulfilled.

Let $a \stin A$ be such that $[a]_f = \sngltn{a}$.
Then $b \neq a$. Set $\pi = \rvpair{b}{a}$ and $c = f(a)$.
(In (IV) of (N2), $a$ will substitute $\bfiy$, and $c$ will substitute $\bfiz$.)
So
\begin{list}{}
{\setlength{\leftmargin}{39pt}
\setlength{\labelsep}{05pt}
\setlength{\labelwidth}{25pt}
\setlength{\itemindent}{-00pt}
\addtolength{\topsep}{-04pt}
\addtolength{\parskip}{-02pt}
\addtolength{\itemsep}{-05pt}
}
\item[$\pmb{(\dagger)}$]
$\rfs{Act}^{\srfs{Gr}}(S) \models
\phi_{\srfs{Oo-pair}}(\bfif,\bfiy,\bfiz\kern1pt)[f,a,c]$.
\vspace{-02.0pt}
\end{list}
Then $f \pi(b) = c$.
Moreover, $(*)$: $(f \pi)\inverse[\sngltn{c}] = \sngltn{b}$.
Here is why. Suppose that $(f \pi)(x) = c$. That is, $f(\pi(x)) = c$.
Hence $\pi(x) = a$. So $x = b$.
We have thus checked that $\pair{b}{c}$ is $f \pi$-simple.
\vspace{2mm}
\\ 
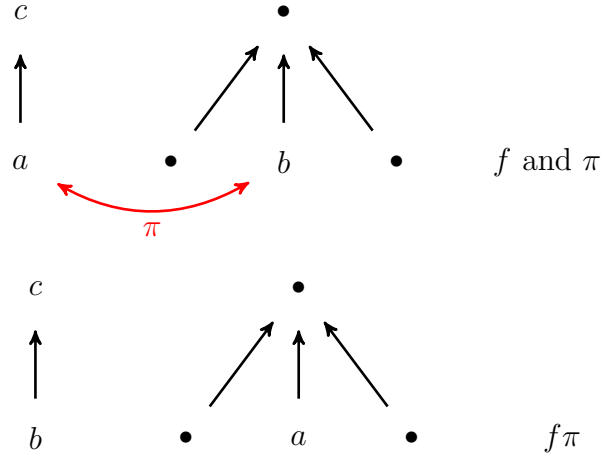
\centerline{\begin{tikzpicture}[->,>=stealth',shorten >=1pt,auto,node distance=2.0cm,on grid,semithick,
                    every state/.style={fill=white,draw=none,circle ,text=black}]
 \node [state] (D)
   {$c$};
 \node [state] (F)[below = of D]
   {$a$}; 
 \node [state] (G) [node distance=3.5cm, right = of D]                                                                                 
    {$\bullet$}; 
 \node [state] (B) [below = of G]                                                                                 
    {$b$};
 \node [state] (O) [node distance=1.5cm, right = of B]                                                                                                                                                    
    {$\bullet$};                   
 \node [state] (A)[node distance=1.5cm, left= of B]
  {$\bullet$};
 \node [state] (C) [right = of O]                                                                                 
    {$f \mbox{ and } \pi$};
\path (O) [line width=1pt] edge node {} (G);
\path (B) [line width=1pt] edge node {} (G);
\path (A) [line width=1pt] edge node {} (G);
\path (F) [line width=1pt] edge node {} (D);
\path   (B) [<->, red, bend left=30,line width=1pt] edge node {$\pi $} (F); 
\end{tikzpicture}}
\centerline{\begin{tikzpicture}[->,>=stealth',shorten >=1pt,auto,node distance=2.0cm,on grid,semithick,
                    every state/.style={fill=white,draw=none,circle ,text=black}]
 \node [state] (D)
   {$c$};
 \node [state] (F)[below = of D]
   {$b$}; 
 \node [state] (G) [node distance=3.5cm, right = of D]                                                                                 
    {$\bullet$}; 
 \node [state] (B) [below = of G]                                                                                 
    {$a$};
 \node [state] (O) [node distance=1.5cm, right = of B]                                                                                                                                                    
    {$\bullet$};                   
 \node [state] (A)[node distance=1.5cm, left= of B]
  {$\bullet$};
 \node [state] (C) [right = of O]                                                                                 
    {$f\pi$};
\path (O) [line width=1pt] edge node {} (G);
\path (B) [line width=1pt] edge node {} (G);
\path (A) [line width=1pt] edge node {} (G);
\path (F) [line width=1pt] edge node {} (D);
\end{tikzpicture}}

The fact that $b \stnot\stin \rfs{Rng}(f)$ implies that
$(**)$: $b \stnot\stin \rfs{Rng}(f \pi)$.
Facts $(*)$ and $(**)$ imply that $c \stnot\stin \rfs{Rng}((f \pi)^2)$.
Also $b \stnot\stin \rfs{Rng}((f \pi)^2)$.
So
\begin{list}{}
{\setlength{\leftmargin}{45pt}
\setlength{\labelsep}{05pt}
\setlength{\labelwidth}{25pt}
\setlength{\itemindent}{-00pt}
\addtolength{\topsep}{-04pt}
\addtolength{\parskip}{-02pt}
\addtolength{\itemsep}{-05pt}
}
\item[\pmb{($\ddagger$)}]
$\rvpair{b}{c} (f \pi)^2 = (f \pi)^2$.
\vspace{-02.0pt}
\end{list}
Substitute $a$ for $\bfiy$, and $c$ for $\bfiz$ in (IV) of (N2).
Then $\pmb{(\dagger)}$ and $\pmb{(\ddagger)}$
show that $f$ and $b$ fulfill (IV) of (N2).
\vspace{2mm}

{\bf Proof of\rule{3pt}{0pt}\ %
$\boldsymbol{\stneg}\bfs{(R1)} \Rightarrow \boldsymbol{\stneg}\bfs{(R2)}$ }
{\thickmuskip=2mu \medmuskip=1mu \thinmuskip=1mu
Suppose that $b \stin \rfs{Rng}(f)$,
and we show that
$\rfs{Act}^{\srfs{Gr}}(S) \stnot\models \phi_{\srfs{Not-in-rng}}[f,b]$.
}%

{\bf Case 1 } $\abs{[b]_f} \geq 2$ or $f(b) = b$.
In this case, there is no $y \stin A$ such that $\pair{b}{y}$ is $f$-simple.
So $f,b$ do not fulfill (N1).

We show that $f,b$ do not fulfill (IV) of (N2).
Suppose by contradiction they do.
Let $a$ be the required $\bfiy$ and $c$ be the required $\bfiz$.
So $a \neq b$. Set $\pi = \rvpair{b}{a}$.
Then $f \pi(b) = f(a)$.
Since $\pair{a}{c}$ is an $f$-1-1 pair, $f(a) = c$.
We show that {\bf ($\boldsymbol{\star}$)}: $c \stin \rfs{Rng}((f \pi)^2)$.
Suppose first that $f(b) \neq b$.
Let $d$ be such that $f(d) = b$.
Then $(f \pi)^2(d) = f \pi f(d) = f \pi(b) = f(a) = c$.
Hence $c \stin \rfs{Rng}((f \pi)^2)$.

Suppose that $f(b) = b$.
Then $(f \pi)^2(a) = f \pi f(b) = f \pi(b) = f(a) = c$.
So again, $c \stin \rfs{Rng}((f \pi)^2)$.
We have shown {\bf ($\boldsymbol{\star}$)}.

Now, {\bf ($\boldsymbol{\star}$)} implies that
$\rvpair{b}{c} (f \pi)^2 \neq (f \pi)^2$.
Hence $\bfiy = a$ and $\bfiz = c$ are not as needed in (IV) of (N2).
A contradition. Hence $f,b$ do not fulfill (IV) of (N2).
So $f,b$ do not fulfill (N2).

{\thickmuskip=6mu \medmuskip=5mu \thinmuskip=4mu
Since $f,b$ do not fulfill (N1) and do not fulfill (N2),
$\rfs{Act}^{\srfs{Gr}}(S) \stnot\models$}\break
$\phi_{\srfs{Not-in-rng}}[f,b]$.

{\bf Case 2 } $\abs{[b]_f} = 1$ and $f(b) \neq b$.
These facts imply that $\pair{b}{f(b)}$ is $f$-simple.
So (III) in (N2) does not hold, and hence $f,b$ do not fulfill
(N2) of $\phi_{\srfs{Not-in-rng}}$.

Let $y$ be such that $\pair{b}{y}$ is $f$-simple. Then $y = f(b)$.
Let $d$ be such that $f(d) = b$. Then $f^2(d) = f(b) = y$.
Hence $y \stin \rfs{Rng}(f^2)$.
By Observation~\ref{p3.3}(g)(1), $\rvpair{b}{y} f^2 \neq f^2$.
So $f,b$ do not fulfill (N1) of $\phi_{\srfs{Not-in-rng}}$.
Hence $\rfs{Act}^{\srfs{Gr}}(S) \stnot\models \phi_{\srfs{Not-in-rng}}[f,b]$.
\vspace{2mm}

(b) Part (b) is a restatement of Part (a).
We stated it just that we can refer to it later.
\rule{20mm}{0pt}\hfill\solidqed
\vspace{2mm}

\begin{cor}\label{c3.7}
There is a formula $\phi_{\srfs{Mo-img}}(\bfif,\bfix\kern1.5pt)$
such that for every\break
$S \stin K_{\srfs{FT}}$,
every $f \stin S$ which has at least three $f$-1-1 images and every $a \stin A$,
$$\mbox{
$\rfs{Act}^{\srfs{Gr}}(S) \models \phi_{\srfs{Mo-img}}[f,a]$
\ \ \ iff\ \ \ \ %
$a \stin \rfs{Mo-img}(f)$.
}$$
\end{cor}

\noindent
{\bf Proof }
$\rfs{Mo-img}(f) = \rfs{Rng}(f) \swsetminus \rfs{Oo-img}(f)$.
So let
$$
\phi_{\srfs{Mo-img}}(\bfif,\bfix\kern1.5pt)
\ \ \eqdf\ \ %
\phi_{\srfs{In-rng}}(\bfif,\bfix\kern1.5pt) \stwedge
\stneg\phi_{\srfs{Oo-img}}(\bfif,\bfix\kern1.5pt).
$$
\index{N@fphimo-img@@$\phi_{\srfs{Mo-img}}(\bfif,\bfix\kern1.5pt)$}%
Let $S \stin K_{\srfs{FT}}$,
and $f \stin S$ be a function with at least three $f$-1-1 pairs.
Let $a \stin A$.
Then by Proposition~\ref{p3.5}(b),
$a \stin \rfs{Rng}(f)$ iff $\rfs{Act}^{\srfs{Gr}}(S) \models
\phi_{\srfs{In-rng}}[f,a]$,
and by Proposition~\ref{p3.3}(e), 
$a \stnot\stin \rfs{Oo-img}(f)$ iff
$\rfs{Act}^{\srfs{Gr}}(S) \models \stneg\phi_{\srfs{Oo-img}}[f,a]$.
So $a \stin \rfs{Mo-img}(f)$ iff
$\rfs{Act}^{\srfs{Gr}}(S) \models \phi_{\srfs{Mo-img}}[f,a]$.
\rule{20mm}{0pt}\hfill\solidqed
\vspace{2mm}

\begin{cor}\label{c3.8}
For every $S \stin K_{\srfs{No-scnst}}$,
$f \stin \rfs{F}_{\srfs{\sbf{B}$\leq$2}}(S)$:
either $f \stin \rfs{Gr}(S)$,
or for every $a \stin A$:
\vspace{-2mm}
$$\mbox{
$\rfs{Act}^{\srfs{Gr}}(S) \models \phi_{\srfs{Mo-img}}[f,a]$
\ \ \ iff\ \ \ \ %
$a \stin \rfs{Mo-img}(f)$.
}$$
\end{cor}

\noindent
{\bf Remark }
It is probably not true that for every $S \stin K_{\srfs{No-scnst}}$,
$f \stin \rfs{F}_{\srfs{\sbf{B}$\leq$2}}(S)$ and $a \stin A$:
\vspace{-4mm}
$$\mbox{
$\rfs{Act}^{\srfs{Gr}}(S) \models \phi_{\srfs{Mo-img}}[f,a]$
\ \ \ iff\ \ \ \ %
$a \stin \rfs{Mo-img}(f)$.
}$$
\vspace{-9mm}

\noindent
The reason is as follows.
For the above equivalence to hold, $f$ needs to have at least three 1-1 pairs.
However, if $\abs{\bfs{Dom}(S)} \stin \threetn{3}{4}{5}$,
then the fact that $S \stin K_{\srfs{No-scnst}}$
and $f \stin \rfs{F}_{\srfs{\sbf{B}$\leq$2}}(S)$ does not imply that
$f$ has three 1-1 pairs.
\vspace{2mm}

\noindent
{\bf Proof }
Let $S \stin K_{\srfs{No-scnst}}$
and $f \stin \rfs{F}_{\srfs{\sbf{B}$\leq$2}}(S)$.
By Observation~\ref{o3.2}, either\break
(i) $f \stin \rfs{Gr}(S)$,
or (ii) there are infinitely many $f$-1-1 pairs.
If (ii) happens, then in particular, $f$ has at least three $f$-1-1 pairs.
So by Corollary~\ref{c3.7}, for every $a \stin A$:
$a \stin \rfs{Mo-img}(f)$ iff
$\rfs{Act}^{\srfs{Gr}}(S) \models \phi_{\srfs{Mo-img}}[f,a]$.
\rule{20mm}{0pt}\hfill\solidqed
\vspace{4mm}

\noindent
\subsection{The definability of the relation
``$\boldsymbol{x \stin \bfs{Fxd-img}(f})$'' 
in $\boldsymbol{J_3}$}

\noindent
{\bf An Overview}
\\
We explain how the formula $\phi_{\srfs{Mo-img}}$ will be used
in finding a formula expressing the property $a \stin \rfs{Fxd-img}(f)$.
Let $\pair{S}{f} \stin J_3$. Hence $\abs{\bfs{B}(f)} \leq 2$.
The cases that $\abs{\bfs{B}(f)} < 2$ are easier.
So let us deal with case that $\abs{\bfs{B}(f)} =~2$.
$\dbltn{\rfs{Oo-pre}(f)}{\rfs{Mo-pre}(f)}$ is a definable partition of $A$.
For members $a \stin\break
\rfs{Oo-pre}(f)$,
it is (by now) trivial to find a formula expressing the property
``$a \stin \rfs{Fxd-img}(f)$.
So it remains to find such a formula for $a$'s belonging to $\rfs{Mo-pre}(f)$.

Since $\abs{\bfs{B}(f)} = 2$, there are $B,C \stin \bfs{B}(f)$
such that $\varPi \eqdf \dbltn{B}{C}$ is a partition of $\rfs{Mo-pre}(f)$.
Clearly, $\varPi$ is definable from $f$.
At this point of the proof we use the formula $\phi_{\srfs{Mo-img}}$.
The set $\rfs{Mo-img}(f)$ has two elements, say, $d$ and~$e$,
and $\phi_{\srfs{Mo-img}}$ defines the set $\dbltn{d}{e}$ from $f$.
It remains to find a formula $\phi(\bfif,\bfix,\bfiy\kern1pt)$ that for
$\bfix \stin B \stcup C$ and for $\bfiy\kern1pt \stin \dbltn{d}{e}$
will say that $\bfif\kern1pt(\bfix\kern1pt) = \bfiy$
and $\bfiy\kern1pt \stin \rfs{Fxd}(\bfif\kern2.3pt)$.

\begin{observation}\label{o8.8}
Let $\pair{S}{f} \stin K_{\srfs{FT}}$.
Then
\vspace{-3mm}
$$
\rfs{Fxd-img}(f) \stcap \rfs{Oo-pre}(f) = \rfs{Idp}(f).
$$
\end{observation}

\noindent
{\bf Proof }
Clearly, $\rfs{Idp}(f) \subseteq \rfs{Oo-pre}(f)$.
Let $b \stin \rfs{Fxd-img}(f) \stcap \rfs{Oo-pre}(f)$.
This means that $f(f(b)) = f(b)$.
Since $b$ is an $f$-1-1 preimage, and both $b$ and $f(b)$ are sent
by $f$ to $f(b)$, $b = f(b)$. The facts ``$b$ is an $f$-1-1 preimage''
and ``$f(b) = b$'' imply that $b \stin \rfs{Idp}(f)$.
\hfill\solidqed
\vspace{3mm}

Since for $\pair{S}{f}$'s belonging to $J_3$,
$f$ has infinitely many $f$-1-1 pairs,
the informal formulation of the formulas $\phi^{i}_{\srfs{Fxd-img}}$
defined below can be replaced by first order formulas in
$\calL_{\srfs{ActGr}}$.
(Remember that many of the subformulas appearing in $\phi^{i}_{\srfs{Fxd-img}}$
have their indended meaning only under the assumption that $f$ has at least
three 1-1 pairs.)
\vspace{3mm}

\begin{prop}\label{p8.3.9}
Let
$\phi^{\srfs{Oo}}_{\srfs{Fxd-img}}(\bfif,\bfix\kern2pt) \eqdf
\phi_{\srfs{Idp}}(\bfif,\bfix\kern2pt)$.
Then for every $\pair{S}{f} \stin~J_3$
and $x \stin \rfs{Oo-pre}(f)$,
$$\mbox{
$x \stin \rfs{Fxd-img}(f)$
\ \ \ iff\ \ \ %
$\rfs{Act}^{\srfs{Gr}}(S) \models \phi^{\srfs{Oo}}_{\srfs{Fxd-img}}[f,x]$.
}$$
\end{prop}

\noindent
{\bf Proof }
The proposition is a trivial corollary from Observation~\ref{o8.8}.
\hfill\solidqed
\vspace{2mm}

\begin{prop}\label{p8.3.10}
Let
$\phi^{\srfs{\sbf{B}=1,1}}_{\srfs{Fxd-img}}(\bfif,\bfix\kern1pt)
\ \ \eqdf\ \ %
\sngltn{\bfix\kern1pt} =
\rfs{Mo-pre}(\bfif\kern2.3pt) \stcap \rfs{Mo-img}(\bfif\kern2.3pt)$
and
$\phi^{\srfs{\sbf{B}=1}}_{\srfs{Fxd-img}}(\bfif,\bfix\kern1pt) \eqdf
(\exists \bfiy \stin [\bfix\kern1pt]_{\sbfif}\kern1.8pt)
\kern2pt\phi^{\srfs{\sbf{B}=1,1}}_{\srfs{Fxd-img}}(\bfif,\bfiy\kern1pt)$.
\underline{Then} 
for every $\pair{S}{f} \stin~J_3$ such that $\abs{\bfB(f)} = 1$,
and for every $x \stin \rfs{Mo-pre}(f)$:
$$\mbox{
$x \stin \rfs{Fxd-img}(f)$
\ \ \ iff\ \ \ %
$\rfs{Act}^{\srfs{Gr}}(S) \models
\phi^{\srfs{\sbf{B}=1}}_{\srfs{Fxd-img}}[f,x]$.
}$$
\end{prop}

\noindent
{\bf Proof }
Let $\pair{S}{f} \stin~J_3$ be such that $\abs{\bfB(f)} = 1$ and $x \stin A$.

{\bf Claim 1 } The following are equivalent:
\begin{itemize}
\addtolength{\parskip}{-11pt}
\addtolength{\itemsep}{06pt}
\item[(1)]
$\rfs{Act}^{\srfs{Gr}}(S) \models
\phi^{\srfs{\sbf{B}=1,1}}_{\srfs{Fxd-img}}[f,x]$.
\item[(2)]
$f(x) = x$.
\vspace{-05.7pt}
\end{itemize}

{\bf Proof }
Suppose that (1):
$\rfs{Act}^{\srfs{Gr}}(S) \models
\phi^{\srfs{\sbf{B}=1,1}}_{\srfs{Fxd-img}}[f,x]$.
That is, $\sngltn{x} = \rfs{Mo-pre}(f) \stcap \rfs{Mo-img}(f)$.
Since $\abs{\bfB(f)} = 1$, $\sngltn{x} = \rfs{Mo-img}(f)$.
Since\break
$x \stin \rfs{Mo-pre}(f)$, $f(x) = x$. Hence (2) holds.

{\thickmuskip=7mu \medmuskip=6mu \thinmuskip=5mu
Suppose that (2): $f(x) = x$.
Then $x \stin \rfs{Mo-pre}(f) \stcap \rfs{Mo-img}(f)$.}\break
Since $\abs{\bfB(f)} = 1$,
$\sngltn{x} = \rfs{Mo-img}(f) \stcap \rfs{Mo-img}(f)$.
That is, 
$\rfs{Act}^{\srfs{Gr}}(S) \models\break
\phi^{\srfs{\sbf{B}=1,1}}_{\srfs{Fxd-img}}[f,x]$.
This proves Claim 1.
\hfill\hollowqed
\vspace{2mm}

{\bf Proof of $\boldsymbol{\Rightarrow}$:}
Suppose that $x \stin \rfs{Fxd-img}(f)$. That is, $f(f(x)) = f(x)$.\break
By Claim 1,
$\rfs{Act}^{\srfs{Gr}}(S) \models
\phi^{\srfs{\sbf{B}=1,1}}_{\srfs{Fxd-img}}[f,x]$.
Also, $f(f(x)) = f(x)$ means that $f(x) \stin [x]_f$.
So $y \eqdf f(x)$ demonstrates that
$\rfs{Act}^{\srfs{Gr}}(S) \models
\phi^{\srfs{\sbf{B}=1}}_{\srfs{Fxd-img}}[f,x]$.

{\bf Proof of $\boldsymbol{\Leftarrow}$:}
Suppose that
$\rfs{Act}^{\srfs{Gr}}(S) \models
\phi^{\srfs{\sbf{B}=1}}_{\srfs{Fxd-img}}[f,x]$.
Let $y$ be such that $y \stin [x]_f$ and
$\rfs{Act}^{\srfs{Gr}}(S) \models
\phi^{\srfs{\sbf{B}=1,1}}_{\srfs{Fxd-img}}[f,y]$.
By Claim 1, $f(y) = y$. So $f(x) = f(y) = y$.
That is, $x \stin \rfs{Fxd-img}(f)$.
\hfill\solidqed
\vspace{2mm}

\begin{prop}\label{p8.3.11}
Let
$$
\phi^{\srfs{\sbf{B}=2,1}}_{\srfs{Fxd-img}}(\bfif,\bfix\kern1pt)
\ \ \eqdf\ \ %
\rfs{Mo-img}(\bfif\kern2.3pt) \subseteq [x_f].
$$
\underline{Then} 
for every $\pair{S}{f} \stin~J_3$ such that
\begin{itemize}
\addtolength{\parskip}{-11pt}
\addtolength{\itemsep}{06pt}
\item[{\rm(i)}]
$\abs{\bfB(f)} = 2$,
\item[{\rm(ii)}]
there is $t$ such that $\rfs{Mo-img}(f) \subseteq [t]_f$,
\vspace{-05.7pt}
\end{itemize}
and for every $x \stin \rfs{Mo-pre}(f)$:
$$\mbox{
$x \stin \rfs{Fxd-img}(f)$
\ \ \ iff\ \ \ %
$\rfs{Act}^{\srfs{Gr}}(S) \models
\phi^{\srfs{\sbf{B}=1}}_{\srfs{Fxd-img}}[f,x]$.
}$$
\end{prop}

\noindent
{\bf Proof }
Let $\pair{S}{f}$ fulfill (i) and (ii).

{\bf Proof of $\boldsymbol{\Rightarrow}$:}
Suppose that $x \stin \rfs{Fxd-img}(f)$.
Let $t$ be as assured by Clause (ii).
Obviously, $y \eqdf f(x) \stin \rfs{Mo-img}(f)$.
So $y \stin [t]_f$. Clearly, $f(y) = f(f(x)) = f(x)$. That is, $y \stin [x]_f$.
Let $z \stin \rfs{Fxd-img}(f) \swsetminus \sngltn{y}$.
So by Clause (ii), $z \stin [t]_f$,
and since $y,z \stin [t]_f$, $z \stin [x]_f$.
Since $\abs{\bfB(f)} = 2$, $\rfs{Mo-img}(f) = \dbltn{y}{z}$.
So $\rfs{Mo-img}(f) \subseteq [x]_f$.
That is,
$\rfs{Act}^{\srfs{Gr}}(S) \models
\phi^{\srfs{\sbf{B}=2,1}}_{\srfs{Fxd-img}}[f,x]$.

{\bf Proof of $\boldsymbol{\Leftarrow}$:}
Suppose that
$\rfs{Act}^{\srfs{Gr}}(S) \models
\phi^{\srfs{\sbf{B}=2,1}}_{\srfs{Fxd-img}}[f,x]$.
Let $\rfs{Mo-img}(f) = \dbltn{y}{z}$.
Then $f(x) = y$ or $f(x) = z$. We may assume that $f(x) = y$.
Since $y \stin [x]_f$, $f(y) = f(x) = y$.
That is, $f(x) = y \stin \rfs{Fxd}(f)$.
Hence $x \stin \rfs{Fxd-img}(f)$.
\rule{20mm}{0pt}\hfill\solidqed
\vspace{2mm}

\begin{prop}\label{p8.3.12}
Let
$\phi^{\srfs{\sbf{B}=2,2}}_{\srfs{Fxd-img}}(\bfif,\bfix\kern1pt)$
be the formula which says
\begin{itemize}
\addtolength{\parskip}{-11pt}
\addtolength{\itemsep}{06pt}
\item[]
For every $\bfiy\kern1pt \stin [x]_{\sbfif} \wsetminus
\left(\rule{0pt}{10pt}\sngltn{\bfix\kern2pt} \stcup
\rfs{Mo-img}(\bfif\kern2.3pt)\right)$
and $\bfiu,\bfiv\kern2pt$:
if $\pair{\bfiu}{\bfiv}\kern1.2pt$ is an $f$-simple pair
and $\bfiv \neq \bfix\kern1pt, \bfiy$,
then
$\bfiu \stin
[\bfix\kern2pt]_{\left(\bfif\kern2pt \srvpair{\bfiv}{\bfiy}\right)^2}$.
\vspace{-05.7pt}
\end{itemize}
\underline{Then} for every $\pair{S}{f} \stin J_3$
and $x \stin \rfs{Mo-pre}(f)$
such that: {\rm(i)} $\abs{[x]_{f}} \geq 3$,
and\break
{\rm(ii)} $\abs{[x]_{f} \stcap \rfs{Mo-img}(f)} = 1$:
\vspace{-3mm}

$$\mbox{
{\rm(1)} $x \stin \rfs{Fxd-img}(f)$
\ \ \ iff\ \ \ %
{\rm(2)} $\rfs{Act}^{\srfs{Gr}}(S) \models
\phi^{\srfs{\sbf{B}=2,2}}_{\srfs{Fxd-img}}[f,x]$.
}$$
\end{prop}

\noindent
{\bf Proof }
Let $\pair{S}{f} \stin J_3$, $x \stin \rfs{Mo-pre}(f)$,
and $f,x$ fulfill Clauses (i) and (ii).

{\bf Proof of (1) $\boldsymbol{\Rightarrow}$ (2):}
Suppose that $x \stin \rfs{Fxd-img}(f)$.
Let
$y \stin\break
[x]_{f} \swsetminus
\left(\rule{0pt}{10pt}\sngltn{x} \stcup \rfs{Mo-img}(f)\right)$
and $\pair{u}{v}$ be an $f$-simple pair such that and\break
{\thickmuskip=2mu \medmuskip=1mu \thinmuskip=1mu
$v \neq x,y$.
For every $z \stin [x]_f$, $f(z) = f(x)$. So since $\abs{[x]_f} \geq 2$,
$f(x) \stin \rfs{Mo-pre}(f)$.} It follows that $y \neq f(x)$.
Also $u$ is an $f$-1-1 preimage and $y$ is not an $f$-1-1 preimage.
So $u \neq y$.
Set $\pi = \rvpair{v}{y}$.
\smallskip 
\\ 
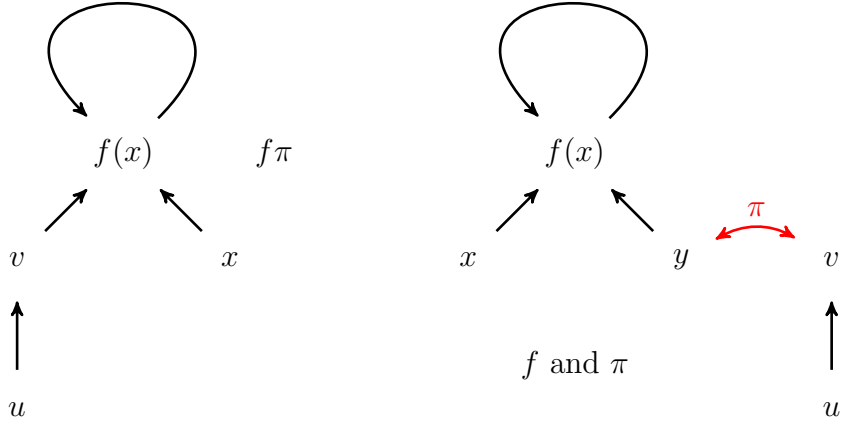
\centerline{\begin{tikzpicture}[->,>=stealth',shorten >=1pt,auto,node distance=2.0cm,on grid,semithick,
                    every state/.style={fill=white,draw=none,circle ,text=black}]
\node [state] (T)
   {$f(x)$}; 
\node [state] (H)[right = of T]
   {$f \pi$}; 
\node [state] (J)[right = of H]
   {}; 
\node [state] (F)[right = of J]
   {$f(x)$};  
 \node [state] (Y)[below right = of F]
   {$y$}; 
   \node [state] (X)[below left = of F]
   {$x$};
 \node [state] (B) [below left = of T]                                                                                 
    {$v$};
 \node [state] (G)[below right= of T]
   {$x$};
 \node [state] (U)[below = of B]
   {$u$};
 \node [state] (O) [right = of Y]                                                                                                                                                    
    {$v$};                   
 \node [state] (C) [below = of O]                                                                                 
    {$u$};
 \node [state] (A) [below left = of Y]                                                                                 
    {$f \mbox{ and } \pi$};
\path (B) [line width=1pt] edge node {} (T);
\path (U) [line width=1pt] edge node {} (B);
\path (G) [line width=1pt] edge node {} (T);
\path (Y) [line width=1pt] edge node {} (F);
\path (X) [line width=1pt] edge node {} (F);
\path  (F) [line width=1pt]  edge [loop]  node [swap] {}(\Cnr);
\path  (T) [line width=1pt]  edge [loop]  node [swap] {}(\Cnr);
\path (C) [line width=1pt] edge node {} (O);
\path   (Y) [<->, red, bend left=30,line width=1pt] edge node {$\pi $} (O); 
\end{tikzpicture}}
\\
We now compute: (i) $(f \pi)^2(u)$, and (ii) $(f \pi)^2(x)$. 
We start with (i).
\\
(1) Since $u \neq y,v$, \ $f \pi(u) = f(u) = v$.
\\
(2) $f \pi(v) = f(y) = f(x)$.
\\
So
\begin{equation}
\tag{I}
(f \pi)^2(u) = f(x).
\end{equation}
(3) It is given that $x \neq y$ and that $x \neq v$.
Hence $f \pi(x) = f(x)$.
\\
(4) $f(x) \stin \rfs{Mo-img}(f)$ and $v,y \stnot\stin \rfs{Mo-img}(f)$,
hence $v,y \neq f(x)$.
So
\vspace{-3mm}
\begin{equation}
\tag{II}
(f \pi)^2(x) = f \pi((f(x)) = f(f(x)).
\vspace{-2mm}
\end{equation}
{\bf Remark } Note that till now, we have not used the fact that
$x \stin \rfs{Fxd-img}(f)$. We shall use this fact in the next step.
\vspace{1mm}

From (II) and the fact that $x \stin \rfs{Fxd-img}(f)$, it follows that
\vspace{-2mm}
\begin{equation}
\tag{III}
(f \pi)^2(x) = f(x)
\vspace{-2mm}
\end{equation}
{\thickmuskip=1mu \medmuskip=1mu \thinmuskip=1mu
(I) and (III) imply that $u \stin [x]_{(f \pi)^2}$.
This shows that
$\rfs{Act}^{\srfs{Gr}}(S) \models
\phi^{\srfs{\sbf{B}=2,2}}_{\srfs{Fxd-img}}[f,x]$.%
}%
\vspace{2mm}

{\bf Proof of $\stneg$(1) $\boldsymbol{\Rightarrow}$ $\stneg$(2):}
Suppose that $x \stnot\stin \rfs{Fxd-img}(f)$.
Since $\abs{[x]_f} \geq 3$ and $\abs{[x]_f \stcap \rfs{Mo-img}(f)} = 1$,
there is $y$ such that 
$y \stin [x]_{f} \swsetminus
\left(\rule{0pt}{10pt}\sngltn{x} \stcup \rfs{Mo-img}(f)\right)$.
By Observation~\ref{o8.3}(b), $f$ has infinitely many $f$-simple pairs.
So there is an $f$-simple pair $\pair{u}{v}$ such that $v \neq x,y$.
We have found $y,u,v$ which fulfill the antecedent of
$\phi^{\srfs{\sbf{B}=2,2}}$.
By Remark 1, we may use facts (I) and (II)
in the proof of (1) $\Rightarrow$ (2).
So $(f \pi)^2(u) = f(x)$ and $(f \pi)^2(x) = f(f(x))$.
Since $x \stnot\stin \rfs{Fxd-img}(f)$, \ $f(f(x)) \neq f(x)$.
So $u \stnot\stin [x]_{(f \pi)^2}$.
This shows that
$\rfs{Act}^{\srfs{Gr}}(S) \stnot\models
\phi^{\srfs{\sbf{B}=2,2}}_{\srfs{Fxd-img}}[f,x]$.%
\hfill\solidqed

\begin{prop}\label{p8.3.13}
$\phi^{\srfs{\sbf{B}=2,3}}_{\srfs{Fxd-img}}(\bfif,\bfix\kern1pt)$
be the formula which says
\begin{itemize}
\addtolength{\parskip}{-11pt}
\addtolength{\itemsep}{06pt}
\item[]
For every
$\bfiy \stin [\bfix\kern1pt]_{\sbfif\kern1pt} \swsetminus
\rfs{Mo-img}(\bfif\kern2.3pt)$,
and an $\bfif$-simple pair $\pair{\bfiu}{\bfiv\kern1pt}$
such that $\bfiv \stnot\stin [\bfix\kern1pt]_{\sbfif\kern1pt}$,
there is
$\bfiz \stin \rfs{Mo-pre}(\bfif\kern2.3pt) \swsetminus
[\bfix]_{\sbfif\kern1pt}$, \ %
such that
$\bfiz \stnot\stin
[\bfiu\kern1pt]_{(\srvpair{\bfiv}{\bfiy\kern1pt} \bfif\kern2pt)^2}$.
\vspace{-05.7pt}
\end{itemize}
\underline{Then} for every $\pair{S}{f} \stin J_3$
and $x \stin \rfs{Mo-pre}(f)$
such that: {\rm(i)} $\abs{[x]_{f}} = 2$,
and\break
{\rm(ii)} $\abs{[x]_{f} \stcap \rfs{Mo-img}(f)} = 1$:
$$\mbox{
{\rm(1)} $x \stin \rfs{Fxd-img}(f)$
\ \ \ iff\ \ \ %
{\rm(2)} $\rfs{Act}^{\srfs{Gr}}(S) \models
\phi^{\srfs{\sbf{B}=2,3}}_{\srfs{Fxd-img}}[f,x]$.
}$$
\end{prop}

\noindent
{\bf Remark }
The ``For every $\bfiy \stin [\bfix\kern1pt]_{\sbfif\kern1pt} \swsetminus
\rfs{Mo-img}(\bfif\kern2.3pt)$'' in the definition of
$\phi^{\srfs{\sbf{B}=2,3}}_{\srfs{Fxd-img}}$ means
``From the two members of $[\bfix\kern1pt]_{\sbfif\kern1pt}$
pick the one which does not belong to $\rfs{Mo-img}(\bfif\kern2.3pt)$''.

\noindent
{\bf Proof }
Let $\pair{S}{f} \stin J_3$, $x \stin \rfs{Mo-pre}(f)$,
and $f,x$ fulfill Clauses (i) and (ii).

{\bf Proof of (1) $\boldsymbol{\Rightarrow}$ (2):}
Suppose that $x \stin \rfs{Fxd-img}(f)$.
Let $y$ be such that
$\sngltn{y} = [x]_f \swsetminus \rfs{Mo-img}(f)$,
and let $\pair{u}{v}$ be an $f$-simple pair such that $v \stnot\stin [x]_f$.
Clearly, $f(y) = f(x) = f(f(x))$. So $f(y) \stin \rfs{Mo-img}(f)$,
and hence $y \neq f(y) = f(f(y))$.
Obviously $v \neq y$. So set $\pi = \rvpair{y}{v}$.
We shall show that for every
$z \stin \rfs{Mo-pre}(f) \swsetminus [x]_f$,
$(\pi f)^2(u) \neq (\pi f)^2(z)$.
\smallskip 
\\ 
\centerline{\begin{tikzpicture}[->,>=stealth',shorten >=1pt,auto,node distance=2.0cm,on grid,semithick,
                    every state/.style={fill=white,draw=none,circle ,text=black}]
 \node [state] (D)
   {$x\in  Fxd-img(f)$};
\node [state] (T)[below left = of D]
   {$f(z)$}; 
\node [state] (H)[right = of T]
   {$f(z)\ne y,v$}; 
   \node [state] (F)[right = of H]
   {$f(x)$};  
 \node [state] (Y)[below = of F]
   {$y$}; 
 \node [state] (B) [below left = of T]                                                                                 
    {$\bullet$};
 \node [state] (G)[below right= of T]
   {$z$};
 \node [state] (O) [right = of Y]                                                                                                                                                    
    {$v$};                   
 \node [state] (C) [below = of O]                                                                                 
    {$u$};
 \node [state] (A) [right = of O]                                                                                 
    {$f \mbox{ and } \pi$};
\path (B) [line width=1pt] edge node {} (T);
\path (G) [line width=1pt] edge node {} (T);
\path (Y) [line width=1pt] edge node {} (F);
\path  (F) [line width=1pt]  edge [loop right]  node [swap] {}(\Cnr);
\path (C) [line width=1pt] edge node {} (O);
\path   (Y) [<->, red, bend left=30,line width=1pt] edge node {$\pi $} (O); 
\end{tikzpicture}}
\centerline{\begin{tikzpicture}[->,>=stealth',shorten >=1pt,auto,node distance=2.0cm,on grid,semithick,
                    every state/.style={fill=white,draw=none,circle ,text=black}]
 \node [state] (D)
   {$f^2(z)$};
 \node [state] (T)[below = of D]
   {$f(z)$}; 
\node [state] (H)[right = of D]
   {$f^2(z)\ne f(x)$}; 
   \node [state] (F)[right = of H]
   {$f(x)$};  
 \node [state] (Y)[below = of F]
   {$y$}; 
 \node [state] (B) [below left = of T]                                                                                 
    {$\bullet$};
 \node [state] (G)[below right= of T]
   {$z$};                  
 \node [state] (C) [below = of Y]                                                                                 
    {$u$};
 \node [state] (A) [right = of O]                                                                                 
    {$\pi f$};
\path (T) [line width=1pt] edge node {} (D);
\path (B) [line width=1pt] edge node {} (T);
\path (G) [line width=1pt] edge node {} (T);
\path (Y) [line width=1pt] edge node {} (F);
\path  (F) [line width=1pt]  edge [loop right]  node [swap] {}(\Cnr);
\path (C) [line width=1pt] edge node {} (Y);
\end{tikzpicture}}
\\

(1)
$(\pi f)^2(u) = \pi f \pi f(u) = \pi f \pi(v) = \pi f(y) =
\pi(f(y)) = f(y) = f(x)$.

Let $z \stin \rfs{Mo-pre}(f) \swsetminus [x]_f$.
Since $v$ is an $f$-1-1 image and $f(z)$ is not an $f$-1-1 image,
$f(z) \neq v$.
Recall that $y \stnot\stin \rfs{Mo-img}(f)$.
Since $z \stin \rfs{Mo-pre}(f)$, it follows that $f(z)\stin \rfs{Mo-img}(f)$.
Hence $f(z) \neq y$.
We now show that $f(f(z)) \neq f(x)$. 
Suppose by contradiction that $f(f(z)) = f(x)$.
Recall that $f(f(x)) = f(x)$.
This means that $f(x) \stin [x]_f \stcap \rfs{Mo-img}(f)$.
Since $\abs{[x]_f \stcap \rfs{Mo-img}(f)} = 1$,
it follows that $\sngltn{f(x)} = [x]_f \stcap \rfs{Mo-img}(f)$.

By the contradiction assumption, $f(z) \stin [x_f]$.
Since $f(z) \stin \rfs{Mo-img}(f)$, we also have that
$f(z) \stin [x]_f \stcap \rfs{Mo-img}(f)$.
This implies that $f(z) = f(x)$. Hence $z \stin [x]_f$.
But this contradicts the fact that
$z \stin \rfs{Mo-pre}(f) \swsetminus [x]_f$.
We conclude that
\\
(2) $f(f(z)) \neq f(x)$. 
\\
We obtain that
\\
(3) $(\pi f)^2(z) = \pi f \pi f(z) = \pi(f(f(z)))$.
\\
By (2) and the fact that $\pi$ is $\onetoone$,
\\
(4) $\pi(f(f(z))) \neq \pi(f(x))$.
\\
Now, $y,v \neq f(x)$, since both $y$ and $v$ do not belong to $\rfs{Mo-img}(f)$,
and $f(x)$ does. Hence
\\
(5) $\pi(f(x)) = f(x)$.
\\
Facts (3), (4) and (5) imply that
\\
(6) $(\pi f)^2(z) \neq f(x)$,
\\
and (1) and (6) imply that $(\pi f)^2(u) \neq (\pi f)^2(z)$.
That is, $z \stnot\stin [u]_{(\pi f)^2}$.

We have thus shown that
$\rfs{Act}^{\srfs{Gr}}(S) \models
\phi^{\srfs{\sbf{B}=2,3}}_{\srfs{Fxd-img}}[f,x]$.%

{\bf Proof of $\stneg$(1) $\boldsymbol{\Rightarrow}$ $\stneg$(2):}
Suppose that $x \stnot\stin \rfs{Fxd-img}(f)$.
Clearly,\break
$f(x) \stnot\stin [x]_f$,
because $f(x) \stin [x]_f$ would mean that $f(f(x)) = f(x)$,
which is equivalent to saying that $x \stin \rfs{Fxd-img}(f)$.
Let $w \stin [x]_f \stcap \rfs{Mo-img}(f)$.
Clearly, $w \neq f(x)$.
This is so since $w \stin [x]_f$ and $f(x) \stnot\stin [x]_f$.

Let $z$ be such that $f(z) = w$.
Since $w \neq f(x)$, it follows that $z \stnot\stin [x]_f$.
Since $w \stin \rfs{Mo-img}(f)$, \ $z \stin \rfs{Mo-pre}(f)$.
Hence
\\
\indent(1) $z \stin \rfs{Mo-pre}(f) \swsetminus [x]_f$.
\\
Since $[x]_f$ is finite, and there are infinitely many $f$-simple pairs,
there is an $f$-simple pair $\pair{u}{v}$ such that
$v \stnot\stin [x]_f$. So
\\
\indent(2) $\pair{u}{v}$ is an $f$-simple pair and $v \stnot\stin [x]_f$.
\\
Let $y \stin [x]_f \swsetminus \sngltn{w}$.
Clearly, $y \neq v$. Set $\pi \eqdf \rvpair{y}{v}$.
We shall show that $(\pi f)^2(z) = (\pi f)^2(u)$.
\smallskip 
\\ 
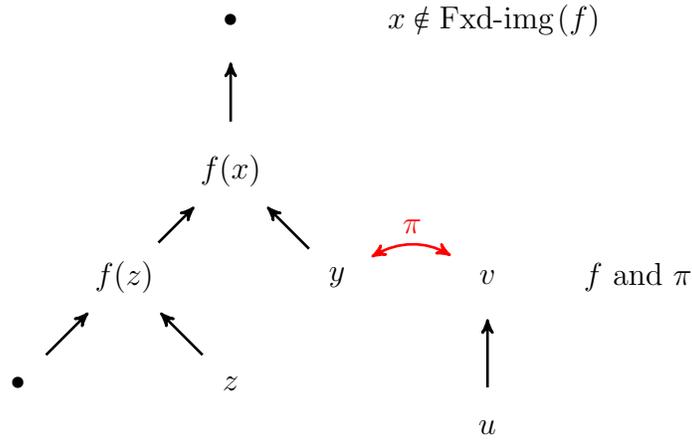
\centerline{\begin{tikzpicture}[->,>=stealth',shorten >=1pt,auto,node distance=2.0cm,on grid,semithick,
                    every state/.style={fill=white,draw=none,circle ,text=black}]
 \node [state] (D)
   {$\bullet$};
 \node [state] (F)[below = of D]
   {$f(x)$}; 
 \node [state] (T)[below left = of F]
   {$f(z)$}; 
 \node [state] (Y)[below right= of F]
   {$y$}; 
 \node [state] (B) [below left = of T]                                                                                 
    {$\bullet$};
 \node [state] (G)[below right= of T]
   {$z$};
 \node [state] (O) [right = of Y]                                                                                                                                                    
    {$v$};                   
 \node [state] (C) [below = of O]                                                                                 
    {$u$};
 \node [state] (A) [right = of O]                                                                                 
    {$f \mbox{ and } \pi$};
\node [state] (E) [node distance=3.5cm, right = of D]                                                                                 
    {$x\notin \rfs{Fxd-img}(f)$};
\path (B) [line width=1pt] edge node {} (T);
\path (G) [line width=1pt] edge node {} (T);
\path (T) [line width=1pt] edge node {} (F);
\path (Y) [line width=1pt] edge node {} (F);
\path (F) [line width=1pt] edge node {} (D);
\path (C) [line width=1pt] edge node {} (O);
\path   (Y) [<->, red, bend left=30,line width=1pt] edge node {$\pi $} (O); 
\end{tikzpicture}}
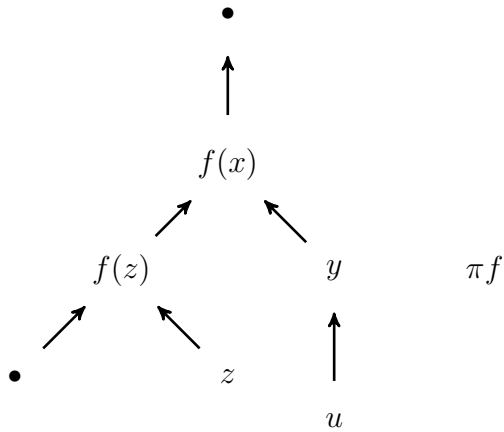
\centerline{\begin{tikzpicture}[->,>=stealth',shorten >=1pt,auto,node distance=2.0cm,on grid,semithick,
                    every state/.style={fill=white,draw=none,circle ,text=black}]
 \node [state] (D)
   {$\bullet$};
 \node [state] (F)[below = of D]
   {$f(x)$}; 
 \node [state] (T)[below left = of F]
   {$f(z)$}; 
 \node [state] (Y)[below right= of F]
   {$y$}; 
 \node [state] (B) [below left = of T]                                                                                 
    {$\bullet$};
 \node [state] (G)[below right= of T]
   {$z$};            
 \node [state] (C) [below = of Y]                                                                                 
    {$u$};
 \node [state] (A) [right = of Y]                                                                                 
    {$\pi f$};
\path (B) [line width=1pt] edge node {} (T);
\path (G) [line width=1pt] edge node {} (T);
\path (T) [line width=1pt] edge node {} (F);
\path (Y) [line width=1pt] edge node {} (F);
\path (F) [line width=1pt] edge node {} (D);
\path (C) [line width=1pt] edge node {} (Y);
\end{tikzpicture}}
\\

From the facts $f(x) \stnot\stin [x]_f$, and $y \stin [x]_f$,
one concludes that $f(x) \neq y$.
Since $v$ is an $f$-1-1 image, and $f(x)$ is not an $f$-1-1 image, \ %
$f(x) \neq v$. So
\\
\indent(3) $\pi(f(x)) = f(x)$.
\\
From (3) we conclude that
\\
\indent(4)
$(\pi f)^2(u) = \pi f \pi f(u) = \pi f \pi(v) = \pi f(y) = \pi f(x) = f(x)$.
\\
Since $y \stin [x]_f \swsetminus \sngltn{w}$, it follows that $w \neq y$.
Since $w \stin [x]_f$ and $v \stnot\stin [x]_f$, we also have that $w \neq v$.
So
\\
\indent(5) $\pi(w) = w$.
\\
Relying on (3) and (5), we now compute $(\pi f)^2(z)$.
\\
\indent(6)
$(\pi f)^2(z) = \pi f \pi f(z) = \pi f \pi(w) = \pi f(w) = \pi f(x) = f(x)$.
\\
Facts (1), (2), (4) and (6) imply that
$\rfs{Act}^{\srfs{Gr}}(S) \stnot\models
\phi^{\srfs{\sbf{B}=2,3}}_{\srfs{Fxd-img}}[f,x]$.%
\hfill\solidqed
\vspace{2mm}

Recall that our goal is to find a formula
$\phi^{J_3}_{\srfs{Fxd-img}}(\bfif,\bfix\kern1pt)$
which for $\pair{S}{f}$'s belonging to $J_3$ expresses the relation
``$\bfix \stin \rfs{Fxd-img}(\bfif\kern2.3pt)$.
To obtain $\phi^{J_3}_{\srfs{Fxd-img}}$,
we shall use the formulas defined in Propositions~\ref{p8.3.9} - \ref{p8.3.13}.
These formulas are:
(1) $\phi^{\srfs{Oo}}_{\srfs{Fxd-img}}$,
(2) $\phi^{\srfs{\sbf{B}=1}}_{\srfs{Fxd-img}}$,
(3) $\phi^{\srfs{\sbf{B}=2,1}}_{\srfs{Fxd-img}}$,
(4) $\phi^{\srfs{\sbf{B}=2,2}}_{\srfs{Fxd-img}}$ and
(5) $\phi^{\srfs{\sbf{B}=2,3}}_{\srfs{Fxd-img}}$.
\\
Denote (1) $\phi^{\srfs{Oo}}_{\srfs{Fxd-img}}$ by
$\phi\fsup{1}_{\srfs{Fxd-img}}$,
(2) $\phi^{\srfs{\sbf{B}=1}}_{\srfs{Fxd-img}}$ by
$\phi\fsup{2}_{\srfs{Fxd-img}}$, and so on.
\vspace{2mm}

We need the following observation.

\begin{observation}\label{o8.3.14}
Let $\pair{S}{f} \stin J_3$ and $x \stin A$.
Assume that $x \stin \rfs{Fxd-img}(\bfif\kern2.3pt)$.
Then one of the following situations happens.
\begin{itemize}
\addtolength{\parskip}{-11pt}
\addtolength{\itemsep}{06pt}
\item[{\rm(1)}]
$x \stin \rfs{Oo-pre}(f)$.
(This is the assumption appearing in Proposition~\ref{p8.3.9}.)
\item[{\rm(2)}]
$x \stin \rfs{Mo-pre}(f)$ and $\abs{\bfB(f)} = 1$.
(This is the assumption appearing in Proposition~\ref{p8.3.10}.)
\item[{\rm(3)}]
$x \stin \rfs{Mo-pre}(f)$, and
\begin{list}{}
{\setlength{\leftmargin}{39pt}
\setlength{\labelsep}{05pt}
\setlength{\labelwidth}{25pt}
\setlength{\itemindent}{-00pt}
\addtolength{\topsep}{-04pt}
\addtolength{\parskip}{-02pt}
\addtolength{\itemsep}{-05pt}
}
\item[{\rm(i)}]
$\abs{\bfB(f)} = 2$, and
\item[{\rm(ii)}]
There is $t$ such that $\rfs{Mo-img}(f) \subseteq [t]_f$.
\vspace{05.0pt}
\end{list}
(This is the assumption appearing in Proposition~\ref{p8.3.11}.)
\item[{\rm(4)}]
$x \stin \rfs{Mo-pre}(f)$, and
\begin{list}{}
{\setlength{\leftmargin}{39pt}
\setlength{\labelsep}{05pt}
\setlength{\labelwidth}{25pt}
\setlength{\itemindent}{-00pt}
\addtolength{\topsep}{03pt}
\addtolength{\parskip}{02pt}
\addtolength{\itemsep}{-01pt}
}
\item[{\rm(i)}] $\abs{[x]_{f}} \geq 3$, and
\item[{\rm(ii)}] $\abs{[x]_{f} \stcap \rfs{Mo-img}(f)} = 1$.
\end{list}
(This is the assumption appearing in Proposition~\ref{p8.3.12}.)
\item[{\rm(5)}]
$x \stin \rfs{Mo-pre}(f)$, and
\begin{list}{}
{\setlength{\leftmargin}{39pt}
\setlength{\labelsep}{05pt}
\setlength{\labelwidth}{25pt}
\setlength{\itemindent}{-00pt}
\addtolength{\topsep}{03pt}
\addtolength{\parskip}{02pt}
\addtolength{\itemsep}{-01pt}
}
\item[{\rm(i)}] $\abs{[x]_{f}} = 2$, and
\item[{\rm(ii)}] $\abs{[x]_{f} \stcap \rfs{Mo-img}(f)} = 1$.
\end{list}
(This is the assumption appearing in Proposition~\ref{p8.3.13}.)
\vspace{-05.7pt}
\end{itemize}
\end{observation}

Note that each of the above five situations can be expressed by an\break
$\calL_{\srfs{ActGr}}$-formula.
Denote the formula which expresses situation~(i) by\break
$\theta_{\srfs{Fxd-img}}\fsup{i}(\bfif,\bfix\kern1pt)$.
\vspace{2mm}

\noindent
{\bf Proof }
Let $\pair{S}{f} \stin J_3$ and $x \stin \rfs{Fxd-img}(f)$.

If $x \stin \rfs{Oo-pre}(f)$, then
$\rfs{Act}^{\srfs{Gr}}(S) \models \theta_{\srfs{Fxd-img}}\fsup{1}[f,x]$.
\vspace{2mm}

Suppose that $x \stnot\stin \rfs{Oo-pre}(f)$.
Then $x \stin \rfs{Mo-pre}(f)$. Hence $\abs{\bfs{B}(f)} \geq 1$.

If $\abs{\bfs{B}(f)} = 1$, then
$\rfs{Act}^{\srfs{Gr}}(S) \models \theta_{\srfs{Fxd-img}}\fsup{2}[f,x]$.
\vspace{2mm}

Suppose that $\abs{\bfs{B}(f)} \neq 1$, then $\abs{\bfs{B}(f)} = 2$.
Hence $\abs{\rfs{Mo-img}(f)} = 2$. Clearly, $f(x) \stin [x]_f$.
{\thickmuskip=2mu \medmuskip=2mu \thinmuskip=1mu
Hence $\abs{[x]_f \stcap \rfs{Mo-img}(f)} \geq 1$.
If $\abs{[x]_f \stcap \rfs{Mo-img}(f)} =~2$,} then
$\rfs{Act}^{\srfs{Gr}}(S) \models \theta_{\srfs{Fxd-img}}\fsup{3}[f,x]$.
\vspace{2mm}

Suppose that $\abs{[x]_f \stcap \rfs{Mo-img}(f)} \neq 2$.
Then $\abs{[x]_f \stcap \rfs{Mo-img}(f)} = 1$.

If $\abs{[x]_{f}} \geq 3$, then
$\rfs{Act}^{\srfs{Gr}}(S) \models \theta_{\srfs{Fxd-img}}\fsup{4}[f,x]$.
\vspace{2mm}

If $\abs{[x]_{f}} \stnot\geq 3$, then $\abs{[x]_{f}} = 2$. So
$\rfs{Act}^{\srfs{Gr}}(S) \models \theta_{\srfs{Fxd-img}}\fsup{5}[f,x]$.
\hfill\solidqed
\vspace{2mm}

\begin{cor}\label{c8.3.15}
Let
$\phi^{J_3}_{\srfs{Fxd-img}}(\bfif,\bfix\kern1pt) \eqdf
\bigvee_{i = 1}^{5}
\left(\rule{0pt}{10pt}\theta\fsup{i}_{\srfs{Fxd-img}}(\bfif,\bfix\kern1pt)
\stwedge
\phi\fsup{i}_{\srfs{Fxd-img}}(\bfif,\bfix\kern1pt)\right)$.
\index{N@fphifxd-img-j3@@$\phi^{J_3}_{\srfs{Fxd-img}}(\bfif,\bfix\kern1pt)$}%
Then for every $\pair{S}{f} \stin J_3$ and $x \stin A$:
\vspace{-3mm}
$$\mbox{
$x \stin \rfs{Fxd-img}(f)$
\ \ \ iff\ \ \ %
$\rfs{Act}^{\srfs{Gr}}(S) \models \phi^{J_3}_{\srfs{Fxd-img}}[f,x]$.
}$$
\end{cor}

\noindent
{\bf Proof }
The corollary follows trivially from Propositions~\ref{p8.3.9} - \ref{p8.3.13}
and Observation~\ref{o8.3.14}.
\hfill\solidqed
\vspace{2mm}

\section{\bf The conclusion of the proof of Theorem A}\label{9}

\begin{cor}\label{c8.1-08-16}
There is a formula $\phi^{\srfs{No-cnst}}_{\srfs{Fxd-img}}(\bfif,\bfix\kern2pt)$
such for every $S \stin K_{\srfs{No-cnst}}$, $f \stin S$ and $a \stin A$:
$a \stin \rfs{Fxd-img}(f)$ iff
$\rfs{Act}^{\srfs{Gr}}(S) \models \phi^{\srfs{No-cnst}}_{\srfs{Fxd-img}}[f,a]$.
\end{cor}

\noindent
{\bf Proof }
The corollary follows easily from nine facts that have been already proved.
These facts are collected into the next three statements.
\vspace{2mm}

{\bf Fact 1 } $\fsetn{J_1}{J_4}$ is a partition of
$K_{\srfs{No-cnst}}^{\srfs{Aug}}$.
\vspace{2mm}

{\bf Fact 2 } For every $i = 1,\ldots,4$,
$\setm{\pair{\rfs{Act}^{\srfs{Gr}}(S)}{f}}{\pair{S}{f} \stin J_i}$
is FO-definable in 
$\setm{\pair{\rfs{Act}^{\srfs{Gr}}(S)}{f}}
{S \stin K_{\srfs{No-cnst}} \mbox{ and } f \stin S}$.
\vspace{2mm}

{\bf Fact \kern-3pt3 } For every $i = 1,\ldots,4$,
there is a formula an $\calL_{\srfs{ActGr}}$-formula\break
$\phi^{J_i}_{\srfs{Fxd-img}}(\bfif,\bfix\kern1pt)$
such that for every $\pair{S}{f} \stin J_i$ and $a \stin A$:
$a \stin \rfs{Fxd-img}(f)$ iff
$\rfs{Act}^{\srfs{Gr}}(S) \models \phi^{J_i}_{\srfs{Fxd-img}}[f,a]$.
\vspace{2mm}

Denote by
$\phi_{J_i}(\bfif\kern2.3pt)$ the formula which defines
$\setm{\pair{\rfs{Act}^{\srfs{Gr}}(S)}{f}}{\pair{S}{f} \stin J_i}$\break
in
$\setm{\pair{\rfs{Act}^{\srfs{Gr}}(S)}{f}}
{S \stin K_{\srfs{No-cnst}} \mbox{ and } f \stin S}$
(see Fact 2).
We skip the trivial proof that
$\phi^{\srfs{No-cnst}}_{\srfs{Fxd-img}}(\bfif,\bfix\kern2pt)$
can be taken to be the following formula.
\vspace{-2mm}
\begin{align*}\mbox{
$\bigvee_{i = 1}^4 \left(\phi_{J_i}(\bfif\kern2.3pt) \kern1pt\stwedge\kern1pt
\phi^{J_i}_{\srfs{Fxd-img}}(\bfif,\bfix\kern1pt)\rule{0pt}{12pt}\right)$.
}
\vspace{-2mm}
\end{align*}
However, we give a references to where the nine facts were proved.
Fact 1 is Observation~\ref{o5.3-08-17}.
By Corollary~\ref{p5.12-08-15}(b),
$(J_1)^{\srfs{ActGr}}$ is an FO-definable subclass of
$((K_{\srfs{FT}})^{\srfs{Aug}})^{\srfs{ActGr}}$,
and by Lemma~\ref{l5.11-08-16},
for every $\pair{S}{f} \stin J_1$,
$\rfs{Fxd-img}(f) =
\setm{a \stin A}{\rfs{Act}^{\srfs{Gr}}(S) \models
\phi_{\srfs{Fxd-img}}^{\srfs{$\exists$-scnst}}[f,a]}$.
(In Subsection~\ref{ss8.1} we denoted
$\phi_{\srfs{Fxd-img}}^{\srfs{$\exists$-scnst}}$ by
$\phi^{J_1}_{\srfs{Fxd-img}}$.)

The same situation happens for $J_2, J_3$ and $J_4$.
The defining formula for $(J_2)^{\srfs{ActGr}}$
is $\phi_{\srfs{\sbf{B}$\geq$3}}(\bfif\kern2pt)$
(See Observation-\ref{o7.1}(b)),
and for every $\pair{S}{f} \stin J_2$,
$\rfs{Fxd-img}(f) =
\setm{a \stin A}{\rfs{Act}^{\srfs{Gr}}(S) \models
\phi^{J_2}_{\srfs{Fxd-img}}[f,b]}$.
See Corollary~\ref{c6.11-08-15}(b).

The defining formula for $(J_3)^{\srfs{ActGr}}$ is $\phi^{J_3}$.
See Observation~\ref{o3.1}.
By Corollary~\ref{c8.3.15},
{\thickmuskip=2mu \medmuskip=1mu \thinmuskip=1mu
for every $\pair{S}{f} \stin J_3$,
$\rfs{Fxd-img}(f) = \setm{a \stin A}
{\rfs{Act}^{\srfs{Gr}}(S) \models \phi^{J_3}_{\srfs{Fxd-img}}[f,a]}$
}%

In Observation~\ref{o7.1}(c) we found the formula
$\phi^{J_4}(\bfif\kern2.3pt)$
which defines\break
$(J_4)^{\srfs{ActGr}}$
in $((K_{\srfs{FT}})^{\srfs{Aug}})^{\srfs{ActGr}}$,
and in Observation~\ref{o8.1}, we found the formula
$\phi^{J_4}_{\srfs{Fxd-img}}(\bfif,\bfix\kern1pt)$
which defines $\rfs{Fxd-img}(f)$ for members of $J_4$.
\hfill\solidqed

\begin{prop}\label{p9.2-15-08-17}
There is an $\calL_{\srfs{ActGr}}$-formula
$\phi_{\srfs{App}}^{\srfs{No-cnst}}(\bfif,\bfix,\bfiy\kern2pt)$
such that for every $S \stin K_{\srfs{No-cnst}}$, $f \stin S$ and $a,b \stin A$:
\vspace{-2mm}
$$\mbox{
$f(a) = b$
\ \ iff\ \ \ %
$\rfs{Act}^{\srfs{Gr}}(S) \models \phi_{\srfs{App}}^{\srfs{No-cnst}}[f,a,b]$.
}$$
\end{prop}

\noindent
{\bf Proof }
Let $\phi_{\srfs{App}}^{\srfs{No-cnst}}(\bfif,\bfix,\bfiy\kern2pt)$
be the formula which says:
\begin{itemize}
\addtolength{\parskip}{-11pt}
\addtolength{\itemsep}{06pt}
\item[(1)]
$\phi^{\srfs{No-cnst}}_{\srfs{Fxd-img}}(\bfif,\bfix\kern2pt)$,
\ $\bfiy \stin [\bfix]_{\sbfif}$,
and there is $\bfiz \stnot\stin [\bfix]_{\sbfif\kern2pt}$ such that
\\
$\stneg\phi^{\srfs{No-cnst}}_{\srfs{Fxd-img}}(\rvpair{y}{z} \bfif,\bfix\kern2pt)$
\item[] or
\item[(2)]
$\stneg\phi^{\srfs{No-cnst}}_{\srfs{Fxd-img}}(\bfif,\bfix\kern2pt)$,
$\bfiy \neq \bfix$ and
$\phi^{\srfs{No-cnst}}_{\srfs{Fxd-img}}
(\rvpair{\bfix}{\bfiy\kern2pt} \bfif,\bfix\kern2pt)$.
\vspace{-05.7pt}
\end{itemize}

Let 
$S \stin K_{\srfs{No-cnst}}$, $f \stin S$ and $a \stin A$. Set $b = f(a)$.

{\bf Case 1 } $a \stin \rfs{Fxd-img}(f)$.
\\
So
$\rfs{Act}^{\srfs{Gr}}(S) \models \phi^{\srfs{No-cnst}}_{\srfs{Fxd-img}}[f,a]$.
The fact ``$a \stin \rfs{Fxd-img}(f)$'' means that $f(b) = b$.
Since $S \stin K_{\srfs{No-cnst}}$ and $f \stin S$,
$f$ is not a constant function.
So there is $z \stin A$ such that $f(z) \neq b$.
This means that $z \stnot\stin [a]_f$.
Set $\sigma = \rvpair{b}{z}$.
Then $\sigma f(a) = \sigma(b) = z$.
The only element $x$ such that $\sigma(x) = z$ is $b$.
So since $f(z) \neq b$,\break
$\sigma f(z) \neq z$.
That is, $\sigma f(\sigma f(a)) \neq \sigma f(a)$.
This implies that
$\rfs{Act}^{\srfs{Gr}}(S) \models\break
\stneg\phi^{\srfs{No-cnst}}_{\srfs{Fxd-img}}[\sigma f,a]$.
It follows that $f$ and $a$ and $b$ satisfy the first disjunct of
$\phi^{\srfs{No-cnst}}_{\srfs{App}}$.
Hence
$\rfs{Act}^{\srfs{Gr}}(S) \models \phi^{\srfs{No-cnst}}_{\srfs{App}}[f,a,b]$.

Let $c \neq b$, and we show that
$\rfs{Act}^{\srfs{Gr}}(S) \stnot\models
\phi^{\srfs{No-cnst}}_{\srfs{App}}[f,a,c]$.
Since $a \stin \rfs{Fxd-img}(f)$,
$\rfs{Act}^{\srfs{Gr}}(S) \models
\phi^{\srfs{No-cnst}}_{\srfs{Fxd-img}}[f,a]$.
So $f,a,c$ do not satisfy the second disjunct of
$\phi^{\srfs{No-cnst}}_{\srfs{App}}$.

If $c \stnot\stin [a]_f$,
then $f,a,c$ do not satisfy the first disjunct of
$\phi^{\srfs{No-cnst}}_{\srfs{App}}$,
and hence
$\rfs{Act}^{\srfs{Gr}}(S) \stnot\models
\phi^{\srfs{No-cnst}}_{\srfs{App}}[f,a,c]$.

So assume $c \stin [a]_f$.
We show that for every $z \stin A \swsetminus [a]_f$,\break
$a \stin \rfs{Fxd-img}(\rvpair{c}{z}f)$.
Let $z \stin A \swsetminus [a]_f$, and set $\sigma = \rvpair{c}{z}$.
Then $\sigma f(a) = \sigma(b)$. We assumed that $b \neq c$.
Since $b \stin [a]_f$ and $z \stnot\stin [a]_f$, $b \neq z$. So $\sigma(b) = b$.
Hence $\sigma f(a) = b$.
Since $a \stin \rfs{Fxd-img}(f)$, and $b = f(a)$, $f(b) = b$.
Since $b \neq c,z$, $\sigma(b) = b$. Hence $\sigma f(b) = b$.
That is, $b \stin \rfs{Fxd}(\sigma f)$.
Recall that $b = \sigma f(a)$. Hence $\sigma f(a) \stin \rfs{Fxd}(\sigma f)$.
That is, $a \stin \rfs{Fxd-img}(\sigma f)$.
We have shown that if in $\phi^{\srfs{No-cnst}}_{\srfs{App}}$,
$\bfiy$ is taken to be $c$, then there is no $z \stin A$
which fulfills the existential requirement
in the first disjunct of $\phi^{\srfs{No-cnst}}_{\srfs{App}}$.
So
$\rfs{Act}^{\srfs{Gr}}(S) \stnot\models \phi^{\srfs{No-cnst}}_{\srfs{App}}[a,c]$.

{\bf Case 2 } $a \stnot\stin \rfs{Fxd-img}(f)$.
\\
Then
$\rfs{Act}^{\srfs{Gr}}(S) \models
\stneg\phi^{\srfs{No-cnst}}_{\srfs{Fxd-img}}[f,a]$.
Clearly, $\rvpair{a}{b} f(a) = a$.
Hence $a \stin \rfs{Fxd-img}(\sigma f)$
and so $S \models \phi_{\srfs{Fxd-img}}[\sigma f,a]$.
It follows that $f,a,b$ satisfy the second disjunct of
$\phi^{\srfs{No-cnst}}_{\srfs{App}}$,
and hence
$\rfs{Act}^{\srfs{Gr}}(S) \models \phi^{\srfs{No-cnst}}_{\srfs{App}}[a,b]$.

Let $c \neq b$.
Clearly, $f,a,c$ do not fulfill the first disjunct of
$\phi^{\srfs{No-cnst}}_{\srfs{App}}$.
We show that $f,a,c$ do not satisfy the second disjunct of
$\phi^{\srfs{No-cnst}}_{\srfs{App}}$.
Set $\sigma = \rvpair{a}{c}$.
Then $\sigma f(a) = \sigma(b)$.
Since $a \stnot\stin \rfs{Fxd-img}(f)$, $b \neq a$. We assumed that $b \neq c$.
So $\sigma(b) = b$.
Hence $\sigma f(a) = b$.
Since $a \stnot\stin \rfs{Fxd-img}(f)$, $f(b) \neq b$.
Clearly, $\sigma f(b)$ is either $a$ or $c$ or $f(b)$.
We have already seen that $a,c,f(b) \neq b$.
Hence $\sigma f(b) \neq b$.
So
\vspace{-2mm}
$$
\sigma f(\sigma f(a)) = \sigma f(b) \neq b = \sigma f(a).
\vspace{-2mm}
$$
That is, $a \stnot\stin \rfs{Fxd-img}(\sigma f)$.
Hence $f, a, c$ do not satisfy the second disjunct of
$\phi^{\srfs{No-cnst}}_{\srfs{App}}$.
So
$\rfs{Act}^{\srfs{Gr}}(S) \stnot\models \phi^{\srfs{No-cnst}}_{\srfs{App}}[a,c]$.
\hfill\solidqed
\vspace{2mm}
%


Part (b) of the next corollary is Theorem E$^*$ mentioned at the end
of Section 3.

\begin{cor}\label{c9.3-08-18}
{\rm(a)}
$K_{\srfs{No-cnst}}^{\srfs{Act}}$ is definably interpretable in
$K_{\srfs{No-cnst}}^{\srfs{ActGr}}$.

{\rm(b)}
$K_{\srfs{No-cnst}}^{\srfs{Act}}$ and $K_{\srfs{No-cnst}}^{\srfs{ActGr}}$
are bi-definably-interpretable.

{\rm(c)}
$K_{\srfs{No-cnst}}^{\srfs{Act}}$ is FS-interpretable in
$K^{\srfs{Alg}}_{\srfs{No-cnst}}$.

{\rm(d)}
$K_{\srfs{FT}}^{\srfs{Act}}$ is FS-interpretable in
$K_{\srfs{FT}}^{\srfs{Alg}}$.
\end{cor}

{\bf Remark }
Notions of interpretability require a specification
of a relation $R \subseteq K \sttimes K^*$.
But in the above corollary the mention of the $R$'s was omitted.
Of course, in Parts (a), (b) and (c) the $R$'s have the form
$R = \setm{\pair{\bfX(S)}{\bfY(S)}}{S \stin K_{\srfs{\bfZ}}}$,
where $\bfX$, $\bfY$ and~$\bfZ$ are chosen appropriately.%
\vspace{2mm}

\vspace{2mm}

\noindent
{\bf Proof }
(a)
For every $\pair{M}{M^*} \stin R$, $\abs{M^*} = \abs{M}$.
So $\phi_{\srfs{U}}$ can be taken to be $\bfix = \bfix$.
The existence of the formula
$\phi^{\srfs{No-cnst}}_{\srfs{App}}$ (Proposition~\ref{p9.2-15-08-17})
means that, indeed, $K_{\srfs{No-cnst}}^{\srfs{Act}}$
is definably interpretable in $K_{\srfs{No-cnst}}^{\srfs{ActGr}}$.

(b)
Recall that
\vspace{-2mm}
$$\mbox{
$\rfs{Act}(S) = \pair{S \stcup A}{\rfs{App}^S}$%
\ \ and\ \ %
$\rfs{Act}^{\srfs{Gr}}(S) =
\trpl{S \stcup A}{\circ}{\rfs{App}^S_{\srfs{Gr}}}.
$
\vspace{-2mm}
}$$
Since $\rfs{Gr}(S)$ is definable in $\rfs{Act}(S)$,
$\rfs{App}^S_{\srfs{Gr}}$ too is definable in $\rfs{Act}(S)$.
This implies the definable interpretability of
$K_{\srfs{No-cnst}}^{\srfs{ActGr}}$ in $K_{\srfs{No-cnst}}^{\srfs{Act}}$.

(c)
By the theorem of McKenzie (Theorem \ref{t1.5}),
$K_{\srfs{No-cnst}}^{\srfs{ActGr}}$ is FS-interp\-retable in
$K^{\srfs{Alg}}_{\srfs{No-cnst}}$.
(This is because $K_{\srfs{No-cnst}} \subseteq K_{\srfs{FT}}$.)
By Part (b),
$K_{\srfs{No-cnst}}^{\srfs{Act}}$ and $K_{\srfs{No-cnst}}^{\srfs{ActGr}}$
are bi-definably-interpretable.
So by the transitivity statement - Observation~\ref{new-o2.6},
$K_{\srfs{No-cnst}}^{\srfs{Act}}$ is FS-interpretable in
$K^{\srfs{Alg}}_{\srfs{No-cnst}}$.

(d) $K_{\srfs{FT}} = K_{\srfs{$\exists$-cnst}} \stcup K_{\srfs{No-cnst}}$.
Let $\phi = \exists \bfif\kern3pt \forall \bfig\kern2pt
(\bfif \bfig = \bfif\kern2pt)$.
Then $\phi(K^{\srfs{Alg}}_{\srfs{FT}}) = K^{\srfs{Alg}}_{\srfs{$\exists$-cnst}}$
and $\stneg\phi(K^{\srfs{Alg}}_{\srfs{FT}}) = K^{\srfs{Alg}}_{\srfs{No-cnst}}$.

By part (b) $K_{\srfs{No-cnst}}^{\srfs{Act}}$ is FS-interpretable
in $K_{\srfs{No-cnst}}$,
and by Proposition~\ref{p4.1-08-15}(c),
$K_{\srfs{$\exists$-cnst}}^{\srfs{Act}}$ is FS-interpretable
in $K^{\srfs{Alg}}_{\srfs{$\exists$-cnst}}$.
Then by Observation~\ref{o2.8},
$K_{\srfs{FT}}^{\srfs{Act}}$ is FS-interpretable in $K^{\srfs{Alg}}_{\srfs{FT}}$.
\hfill\solidqed
\vspace{2mm}

\noindent
{\bf Proof of Theorem A in the introduction }
Corollary~\ref{c9.3-08-18}(d) is a restatement of Theorem A.
\hfill\solidqed
\vspace{2mm}

\noindent
{\bf Proof of Corollary B in the introduction }
The corollary follows from Theorem A and Theorem~\ref{t2.4-08-14}.
\hfill\solidqed

\section{\bf A transfer theorem from function semigroups to clones,
and a corollary for FT clones}
\label{s10}

\subsection{Reconstructing the action of a semi-transitive clone
from the action of its semigroup}
\label{ss10.1}

Let $S$ be a function semigroup on $A$.
One says that $S$ is a \underline{transitive semigroup},\break
if for every $a,b \stin A$ there is $f \stin S$
such that $f(a) = b$.
\index{transitive semigroup@@Transitive semigroup}%

In fact, in what follows a weaker assumption which we call
``semi-tran-\break
sitivity'' suffices.
Let $S$ be a function semigroup on $A$. One says that $S$ is a
\underline{semi-transitive function semigroup},
if for every $a,b \stin A$ there are $c \stin A$ and
$f,g \stin S$ such that $f(c) = a$ and $g(c) = b$.
\index{semi transitive function semigroup@@Semi-transitive semigroup}%
A clone $C$ on $A$ is called a \underline{semi-transitive clone},
if $C \stcap A^A$ is a semi-transitive semigroup.
\index{semi transitive clone@@Semi-transitive clone}%
\vspace{2mm}

The algebraic structure and the action structure of a clone were defined in the
introduction.
$$
\rfs{Alg}(C) =
\pair{C}{\sngltn{\rfs{Cmp}^C_{n,k}}_{n,k \stin \bbN^+}}
\vspace{-5mm}
$$
and
$$
\rfs{Act}(C) \eqdf
\pair{C \stcup A}{\sngltn{\rfs{App}_n}_{n \stin \bbN^+}}.
$$

We make the notations precise.
We denote the language of $\rfs{Alg}(C)$ by $\calL_{\srfs{Alg-cln}}$.
\index{N@lcln@@$\calL_{\srfs{Alg-cln}}$. The language of $\rfs{Alg}(C)$}%
Hence $\calL_{\srfs{Alg-cln}} = \setm{\rfs{Cmp}_{n,k}}{n,k \stin \bbN^+}$,
where $\rfs{Cmp}_{n,k}$ is an\break
$(n + 2)$-place relation symbol.
Of course, $(\rfs{Cmp}_{n,k})^{\srfs{Alg}(C)} = \rfs{Cmp}^C_{n,k}$.\break
Similarly,
the language of $\rfs{Act}(C)$ is $\setm{\rfs{App}_n}{n \stin \bbN^+}$,
where $\rfs{App}_n$ is an $(n + 2)$-place relation symbol.
It is denoted by $\calL_{\srfs{Act-cln}}$.
Also, $(\rfs{Act}_n)^{\srfs{Act}(C)} = \rfs{Act}^C_n$.
Recall that even though $\rfs{Cmp}_{n,k}$ and $\rfs{App}_n$ are
defined to be relation symbols, they are written as functions.
That is,
$\rfs{Cmp}^C_{n,k}(f,g_1,\ldots,g_n) = h$ means that
$\fivetpl{f}{g_1}{\ldots}{g_n}{h} \stin \rfs{Cmp}^C_{n,k}$,
and
$\rfs{App}^C_n(f,a_1,\ldots,a_n) = b$ means that
$\fivetpl{f}{a_1}{\ldots}{a_n}{b} \stin \rfs{App}^C_n$.

For a clone $C$ we now define the structures
$\rfs{Act}^+(C)$ and $\rfs{Act}^1(C)$.
The language $\calL^{\srfs{Act}\kern-1pt+}_{\srfs{Cln}}$ of $\rfs{Act}^+(C)$ is
$\setm{\rfs{App}_n}{n \stin \bbN^+} \stcup \setm{\rfs{Cmp}_{n,k}}{n,k \stin \bbN^+}$
and
\begin{equation*}
\rfs{Act}^+(C) \eqdf
\trpl{C \stcup A}{\sngltn{\rfs{App}^C_n}_{n \stin \bbN^+}}
{\sngltn{\rfs{Cmp}^C_{n,k}}_{n.k \stin \bbN^+}}.
\end{equation*}
\vspace{-1mm}
\index{N@act2@@$\rfs{Act}^+(C) \eqdf \trpl{C \stcup A}{\sngltn{\rfs{App}^C_n}_{n \stin \bbN^+}}{\sngltn{\rfs{Cmp}^C_n}_{n \stin \bbN^+}}$}%
The structure $\rfs{Act}^1(C)$ is the following reduct of $\rfs{Act}^+(C)$.
$$
\rfs{Act}^1(C) \eqdf
\trpl{C \stcup A}{\rfs{App}^C_1}{\sngltn{\rfs{Cmp}^C_n}_{n \stin \bbN^+}}.
\vspace{-6mm}
$$
\index{N@act3@@$\rfs{Act}^1(C) \eqdf \trpl{C \stcup A}{\rfs{App}^C_1}{\sngltn{\rfs{Cmp}^C_n}_{n \stin \bbN^+}}$}%

Let $K_{\srfs{Sm-trnsv}}$ be the class of all semi-transitive clones,
\index{N@ksemtrns@@$K_{\srfs{Sm-trnsv}}$}%
$$
K^{\srfs{Act}\nsup{\kern-0.6pt+}}_{\srfs{Sm-trnsv}} \eqdf
\setm{\rfs{Act}(C)}
{C \stin K_{\srfs{Sm-trnsv}}}
\vspace{-4mm}
$$
\index{N@kactsmtrnsv@@$K^{\srfs{Act}}_{\srfs{Sm-trnsv}}$}%
and
$$
K^{\srfs{Act}\nsup{1}}_{\srfs{Sm-trnsv}} \eqdf
\setm{\rfs{Act}^1(C)}{C \stin K_{\srfs{Sm-trnsv}}}.
\vspace{-4mm}
$$
\index{N@kact1@@$K^{\srfs{Act}\nsup{1}}_{\srfs{Sm-trnsv}}$}%

This subsection is devoted to the proof of the following theorem.

\begin{theorem}\label{t10.1}
$K^{\srfs{Act}\nsup{\kern-0.6pt+}}_{\srfs{Sm-trnsv}}$
is definably interpretable in $K^{\srfs{Act}\nsup{1}}_{\srfs{Sm-trnsv}}$.
\end{theorem}

\begin{observation}\label{was-o8.2}\label{o10.2}
Let $S$ be a semi-transitive function semigroup on $A$.
Then for every $n \geq 2$ and $a_1,\ldots,a_n \stin A$ there are
$c \stin A$ and $g_1,\ldots,g_n \stin S$
such that for every $i = 1,\ldots,n$, $g_i(c) = a_i$.
\end{observation}

\noindent
{\bf Proof }
By definition, the claim of the observation is true for $n = 2$.
Suppose that the claim of the observation is true for $n$,
and let $a_1,\ldots,a_{n + 1} \stin A$.
By the induction hypothesis there are $d \stin A$ and
$h_1,\ldots,h_n \stin S$ such that for every
$i = 1,\ldots,n$, $h_i(d) = a_i$.
Let $c \stin A$ and $k_1,k_2 \stin S$ be such that $k_1(c) = d$
and $k_2(c) = a_{n + 1}$.
For $i = 1,\ldots,n$ let $g_i = h_i k_1$ and let $g_{n + 1} = k_2$.
Then $c,g_1,\ldots,g_{n + 1}$ are as required.
\vspace{2mm}
\hfill\solidqed

\noindent
{\bf Proof of Theorem~\ref{t10.1}}
Since $\abs{\rfs{Act}^+(C)} = \abs{\rfs{Act}^1(C)}$,
the formula $\phi_{\srfs{U}}$ which defines the universe of
$\abs{\rfs{Act}^1(C)}$ in $\rfs{Act}^+(C)$ can be taken to be
$\bfix = \bfix$.
The relation symbols which belong to
$\calL(\rfs{Act}^+(C)) \swsetminus \calL(\rfs{Act}^1(C))$,
are $\rfs{App}_n$ for $n \geq 2$.
So it remains to show that for every $n \geq 2$,
$\rfs{App}^C_n$ is definable in $\rfs{Act}^1(C)$.
To be precise, we have to prove that:
\begin{list}{}
{\setlength{\leftmargin}{39pt}
\setlength{\labelsep}{05pt}
\setlength{\labelwidth}{25pt}
\setlength{\itemindent}{-00pt}
\addtolength{\topsep}{-04pt}
\addtolength{\parskip}{-02pt}
\addtolength{\itemsep}{-05pt}
}
{\thickmuskip=2mu \medmuskip=1mu \thinmuskip=1mu
\item[$(\star)$]
for every $n \geq 2$ there is an $\calL(\rfs{Act}^1(C))$-formula
$\phi^{\srfs{App}}_n(\bfif,\bfix\Fsub{1},\ldots,\bfix\Fsub{n},\bfy)$}
such that for every $C \stin K_{\srfs{Sm-trnsv}}$,
$f \stin \rfs{Fnc}_n(C)$ and $a_1,\ldots,a_n,b \stin A$:%
\vspace{-09.0pt}
\begin{align*}
\mbox{(I) } f(a_1,\ldots,a_n) = b\ \ \ \mbox{iff}\ \ \ %
\mbox{(II) } \rfs{Act}^1(C) \models \phi^{\srfs{App}}_n[f,a_1,\ldots,a_n,b].
\end{align*}
\vspace{-25.0pt}
\end{list}

There are $\calL(\rfs{Act}^1(C))$-formulas
$\psi_{\pmb{\kern1.3pt\rm D}}(\bfix\kern1pt)$
and for every $n \stin \bbN^+$ a formula
$\psi^n_{\pmb{\kern1pt\rm F}}(\bfix\kern1pt)$ 
such that for every clone $C$ and $p \stin C$,
$\rfs{Act}^1(C) \models \psi_{\pmb{\kern1.3pt\rm D}}[p]$ iff\break
$p \stin \bfs{Dom}(C)$,
and
$\rfs{Act}^1(C) \models
\psi^n_{\pmb{\kern1pt\rm F}}[p]$ iff $p \stin \rfs{Fnc}_n(C)$.
We skip the verification of these trivialities.

For $n \geq 2$ let
$\phi^{\srfs{App}}_n(\bfif,\bfix\Fsub{1},\ldots,\bfix\Fsub{n},\bfy)$
be the $\calL(\rfs{Act}^1(C))$-formula which says:
\begin{list}{}
{\setlength{\leftmargin}{24pt}
\setlength{\labelsep}{05pt}
\setlength{\labelwidth}{10pt}
\setlength{\itemindent}{-00pt}
\addtolength{\topsep}{-04pt}
\addtolength{\parskip}{-02pt}
\addtolength{\itemsep}{-05pt}
}
\item[]
$\psi^n_{\pmb{\kern1pt\rm F}}(\bfif\kern2.3pt)$, and
there are $\bfiz$ and $\bfig\fsub{1},\ldots,\bfig\fsub{n}$ such that\ \ %
$\psi_{\pmb{\kern1.3pt\rm D}}(\bfiz\kern1pt) \kern2pt\stwedge\kern2pt
\bigwedge_{i = 1}^n\psi^1_{\pmb{\kern1pt\rm F}}(\bfig\fsub{i})$\rule{15mm}{0pt}
and
\begin{list}{}
{\setlength{\leftmargin}{39pt}
\setlength{\labelsep}{05pt}
\setlength{\labelwidth}{25pt}
\setlength{\itemindent}{-00pt}
\addtolength{\topsep}{-04pt}
\addtolength{\parskip}{-02pt}
\addtolength{\itemsep}{-05pt}
}
\item[(C1)]
$\bigwedge_{i = 1}^n
\rfs{App}_1(\bfig\fsub{i},\bfiz\kern2pt) = \bfix\fsub{i}$, \ and
\vspace{2mm}
\item[(C2)]
$\rfs{App}_1
\left(\rule{0pt}{12pt}
\rfs{Cmp}_{n,1}(\bfif,\bfig\fsub{1},\ldots,\bfig\fsub{n}),\bfiz\kern2pt\right) =
\bfiy$.
\end{list}
\end{list}

(II) $\rightarrow$ (I): (This implication is true for any clone, not just for
semi-transitive ones.)
Let $C$ be a clone, $f \stin C$ and $b, a_1,\ldots,a_n \stin \bfs{Dom}(C)$
be such that $\rfs{Act}^1(S) \models \phi^{\srfs{App}}_n[f,a_1,\ldots,a_n,b]$.
So $f \stin \rfs{Fnc}_n(C)$,
and there are $c \stin \bfs{Dom}(C)$ and $g_1,\ldots,g_n \stin \rfs{Fnc}_1(C)$
such that
\begin{equation}
\tag{1}
\mbox{
For every $i = 1,\ldots,n$,\ \ %
$\rfs{App}^C_1(g_i,c) = a_i$
}
\end{equation}
and
\begin{equation}
\tag{2}
\rfs{App}^C_1
\left(\rule{0pt}{12pt}\rfs{Cmp}^C_{n,1}(f,g_1,\ldots,g_n),c\right) = b.
\end{equation}
(1) means that
\begin{equation}
\tag{3}
\mbox{
For every $i = 1,\ldots,n$,\ \ %
$g_i(c) = a_i$,
}
\end{equation}
and (2) means that
\begin{equation}
\tag{4}
\rfs{Cmp}^C_{n,1}(f,g_1,\ldots,g_n)(c) = b.
\end{equation}
That is,
\begin{equation}
\tag{5}
f(g_1(c),\ldots,g_n(c)) = b.
\end{equation}
By (3) and (5),
\begin{equation*}
f(a_1,\ldots,a_n) = b.
\end{equation*}

(I) $\rightarrow$ (II): (Here we shall use the semi-transitivity of $C$.)
Suppose that $f(a_1,\ldots,a_n) = b$. By the semi-transitivity of $C$
and Observation~\ref{o10.2},
there are $c \stin A$ and $g_1,\ldots,g_n \stin \rfs{Fnc}_1(C)$
such that for every $i = 1,\ldots,n$, $g_i(c) = a_i$.
Clearly,
$\rfs{Act}^1(C) \models \psi^n_{\pmb{\kern1pt\rm F}}[f]$,
$\rfs{Act}^1(C) \models \psi_{\pmb{\kern1.3pt\rm D}}[c]$,
and for every $i = 1,\ldots,n$,
$\rfs{Act}^1(C) \models \psi^1_{\pmb{\kern1pt\rm F}}[g_i]$.
The fact ``$g_i(c) = a_i$'' can be written as
$g_i(c) = \rfs{App}^C_1(g_i,c)$.
So
\begin{equation}
\tag{6}
\rfs{App}^C_1(g_i,c) = a_i.
\end{equation}
That is,
$(\rfs{App}_1)^{\srfs{Act}^1(C)}(g_i,c) = a_i$.
So
$$
\rfs{Act}^1(C) \models
\left(\rule{0pt}{12pt}\rfs{App}_1(\bfig\fsub{i},\bfiz\kern1pt) =
\bfix\fsub{i}\right)
[g_i,c,a_i].
$$
Hence $a_1,\ldots,a_n,g_1,\ldots,g_n,c$ fulfill (C1) of $\phi^{\srfs{App}}_n$.

Obviously,
\begin{equation}
\tag{7}
f(g_1(c),\ldots,g_n(c)) = \rfs{Cmp}^C_{n,1}(f,g_1,\ldots,g_n)(c).
\end{equation}
So by (6) and (7),
\begin{equation*}
\rfs{App}^C_1
\left(\rule{0pt}{12pt}
\rfs{Cmp}^C_{n,1}(f,g_1,\ldots,g_n),c\right) = b.
\end{equation*}
Recall that (C2) is the formula
$$\rfs{App}_1
\left(\rule{0pt}{12pt}
\rfs{Cmp}_{n,1}(\bfif,\bfig\fsub{1},\ldots,\bfig\fsub{n}),\bfiz\kern2pt\right) =
\bfiy.$$
Hence $f,g_1,\ldots,g_n,c,b$ satisfy (C2).

It follows that $C \models \phi^{\srfs{App}}_n[f,a_1,\ldots,a_n,b]$.
\hfill\solidqed
\vspace{2mm}

\subsection{Corollary: A reconstruction theorem for FT clones}\label{ss10.2}

In this subsection we shall prove Theorem C in the introduction.
It says that
$K^{\srfs{Act}}_{\srfs{FT-cln}}$
is FS-interpretable in $K^{\srfs{Alg}}_{\srfs{FT-cln}}$.

\begin{cor}\label{c10.3} $($From Theorem~\ref{t10.1}.$)$
\\
$\setm{\rfs{Act}^1(C)}{C \stin K_{\srfs{FT-cln}}}$
and
$\setm{\rfs{Act}^+(C)}{C \stin K_{\srfs{FT-cln}}}$
are bi-definably-interp\-retable.
\end{cor}

\noindent
{\bf Proof }
$\rfs{Act}^1(C)$ is a reduct of $\rfs{Act}^+(C)$.
So trivially,
$\setm{\rfs{Act}^1(C)}{C \stin K_{\srfs{FT-cln}}}$
is definably-interpretable in
$\setm{\rfs{Act}^+(C)}{C \stin K_{\srfs{FT-cln}}}$.
The other\break
definable-interpretability follows from Theorem~\ref{t10.1}.
This is so, since\break
$K_{\srfs{FT-cln}} \subseteq K_{\srfs{Sm-trnsv}}$.
\rule{30mm}{0pt}\hfill\solidqed
\vspace{2mm}

\begin{observation}\label{o10.4}
$\setm{\rfs{Act}(C)}{C \mbox{ is a clone}}$ and
$\setm{\rfs{Act}^+(C)}{C \mbox{ is a clone}}$ are bi-definably-interpretable.
\end{observation}

\noindent
{\bf Proof }
Since $\rfs{Act}(C)$ is a reduct of $\rfs{Act}^+(C)$,
$\setm{\rfs{Act}(C)}{C \mbox{ is a clone}}$ is definable in
$\setm{\rfs{Act}^+(C)}{C \mbox{ is a clone}}$.

Next, we define the $\calL_{\srfs{Act-cln}}$-formula $\psi_{\srfs{Cmp}_{n,k}}$.
Let
$$
\psi_{\bm{D}}(\bfix\kern1pt) \eqdf
\exists \bfiy\kern1pt \exists \bfiz\kern1pt
\left(\rule{0pt}{10pt}\rfs{Cmp}_{1,1}(\bfiy,\bfix\kern1pt) =
\bfiz\kern1pt\right).
$$
Then for every clone $C$, $\psi_{\bm{D}}[\rfs{Act}(C)] = \bfs{Dom}(C)$.
Let $\vecbfia$ denote $\bfia_1,\ldots,\bfia_k$.
Define
\begin{align*}
&
\psi_{\srfs{Cmp}_{n,k}}(\bfif,\bfig\fsub{1},\ldots,\bfig\fsub{n},\bfih\kern1pt)
\eqdf
\forall \bfia_1,\ldots,\bfia_k
\left(\rule{0pt}{16pt}
\mbox{$\bigwedge_{i = 1}^k$}\psi_{\bm{D}}(\bfia_i) \rightarrow
\right.
\\&
\left(\rule{0pt}{12pt}
\rfs{App}_k(\bfih,\vecbfia\kern1pt) = 
\rfs{App}_n\left(\rule{0pt}{10pt}\bfif,\rfs{App}_k(\bfig\fsub{1},\vecbfia\kern1pt),
\ldots, \rfs{App}_k(\bfig\fsub{n},\vecbfia\kern1pt)\right)
\kern-2pt\rule{0pt}{12pt}\right)
\left.\kern-3pt\rule{0pt}{16pt}\right).
\end{align*}
It is obvious that for every clone $C$,
$\psi_{\srfs{Cmp}_{n,k}}$ defines $\rfs{Cmp}^C_{n,k}$ in $\rfs{Act}(C)$.
\hfill\solidqed
\vspace{2mm}

We now use the reconstruction theorem for FT semigroups (Theorem A).

\begin{cor}\label{c10.5}
$($From Theorem~A in the introduction.$)$
\\
$\setm{\rfs{Act}^1(C)}{C \stin K_{\srfs{FT-cln}}}$ is FS-interp\-retable in
$\setm{\rfs{Alg}(C)}{C \stin K_{\srfs{FT-cln}}}$.
\end{cor}

\noindent
{\bf Proof }
The proof is trivial, but writing its details is long.
For a clone $C$ let $\rfs{Smgr}(C) \eqdf C \stcap A^A$.
Also, let
$K \eqdf \setm{\rfs{Smgr}(C)}{C \stin K_{\srfs{FT-cln}}}$.
Then $K \subseteq K_{\srfs{FT}}$.
(In fact,
$K = \setm{S \stin K_{\srfs{FT}}}{\rfs{Id}_A \stin S}$,
but this fact is not needed.)

Recall that the required FO-interpretation has to consist of the following
$\calL_{\srfs{Alg-Cln}}$-formulas:
$\psi_{\srfs{U}}, \psi_{=}$ and $\psi_{\srfs{App}\nsub{1}}$
and $\setm{\psi_{\srfs{Cmp}\nsub{n,k}}}{n,k \stin \bbN^+}$.
To obtain an FS-interpretation, we need two additional formulas:
$\psi^*_{\srfs{Imap}}$ in $\calL(\rfs{Act}^1(C))$
and $\psi_{\srfs{Imap}}$ in $\calL_{\srfs{Alg-Cln}}$.

Since $K \subseteq K_{\srfs{FT}}$,
$K^{\srfs{Act}}$ is FS-interpretable in $K^{\srfs{Alg}}$.
This follows from Theorem A.

Let $C \stin K_{\srfs{FT-cln}}$ and $S = \rfs{Smgr}(C)$.
$\rfs{Alg}(S)$ was defined to be $\pair{S}{\circ^S}$,
where $\circ^S$ is regarded as 3-place relation (see Definition~\ref{d1.1}(b)).
So in fact, $\circ^S = \rfs{Cmp}_{1,1}^C$.


Let $\theta_{\srfs{U}}(\vecbfif\kern2.6pt)$,
$\theta_{=}(\vecbfif,\vecbfig\kern1.5pt)$,
$\theta_{\srfs{App}}(\vecbfif,\vecbfig,\vecbfih\kern1.5pt)$,
$\theta_{\srfs{Imap}}(\vecbfif,\bfig\kern1pt)$
and 
$\theta^*_{\srfs{Imap}}(\vecbfif,\bfig\kern1pt)$
be the formulas which constitute the FS-interpretation of 
$K^{\srfs{Act}}$ in $K^{\srfs{Alg}}$.
Assume that $\vecbfif = \fseqn{\bfif\Fsub{1}}{\bfif\Fsub{k}}$,
$\vecbfig = \fseqn{\bfig\fsub{1}}{\bfig\fsub{k}}$
and $\vecbfih = \fseqn{\bfih\fsub{1}}{\bfih\fsub{k}}$.
Let $\phi_{\srfs{Smgr}}(\bfif\kern2.3pt)$ be an $\calL_{\srfs{Alg-Cln}}$-formula,
such that for every clone $C$,
$\phi_{\srfs{Smgr}}[\rfs{Alg}(C)]  = \rfs{Smgr}(C)$.
We first define $\psi_{\srfs{U}}, \psi_{=}$ and $\psi_{\srfs{App}\nsub{1}}$.
\begin{align*}
\mbox{
$\psi_{\srfs{U}}(\vecbfif\kern2.6pt) \eqdf$}
\mbox{$\left(
\bigwedge_{i = 1}^k\phi_{\srfs{Smgr}}(\bfif\Fsub{i}) \stwedge
(\theta_{\srfs{U}})^{(\phi_{\trfs{Smgr}})}(\vecbfif\kern2.6pt)\right)
\stvee$\kern3pt}
\fstneg \phi_{\srfs{Smgr}}(\bfif\Fsub{1})
\end{align*}
(Recall that the formula $(\theta_{\srfs{U}})^{(\phi_{\trfs{Smgr}})}$
is the relativization of
$\theta_{\srfs{U}}$ to $\phi_{\srfs{Smgr}}(\bfix\kern1pt)$.
See Proposition~\ref{p2.7}.)

\underline{Explanation}:
{\thickmuskip=8mu \medmuskip=6mu \thinmuskip=6mu
$\abs{\rfs{Act}^1(C)} = \abs{\rfs{Act}(\rfs{Smgr}(C))} \stcup
\left(\rule{0pt}{12pt}
\abs{\rfs{Alg}(C)} \swsetminus \rfs{Smgr}(C)\right)$.}\break
The members of $\abs{\rfs{Act}(\rfs{Smgr}(C))}$
are represented by $k$-tuples
which fulfill the first disjunct of $\psi_{\srfs{U}}$,
and the members of $\abs{\rfs{Alg}(C)} \swsetminus \rfs{Smgr}(C)$
are represented by $k$-tuples fulfilling the second disjunct of
$\psi_{\srfs{U}}$.

More details: Let $p \stin \rfs{Act}(\rfs{Smgr}(C))$.
Then the representatives of $p$ in $\rfs{Alg}(C)$ are
the same as the representatives of $p$ in $\rfs{Alg}(\rfs{Smgr}(C))$.
That is, for every $\vecf \stin (\rfs{Alg}(C))^k$:
$\vecf$ represents $p$ in $\rfs{Alg}(C)$ iff 
$\vecf$ represents $p$ in $\rfs{Alg}(\rfs{smgr}(C))$.

Note that
$\abs{\rfs{Act}^1(C)} \swsetminus \abs{\rfs{Act}(\rfs{Smgr}(C))} =
\abs{\rfs{Alg}(C)} \swsetminus \rfs{Smgr}(C)$.
For $p$'s which belong to
$\abs{\rfs{Act}^1(C)} \swsetminus \abs{\rfs{Act}(\rfs{Smgr}(C))}$:
a $k$-tuple $\vecf \eqdf \seqn{f_1}{f_k} \stin (\rfs{Alg}(C))^k$ represents $p$
iff $f_1 = p$.

According to the above explanation, $\psi_{=}$ should be defined as follows:
\begin{align*}
\mbox{
$
\psi_{=}(\vecbfif,\vecbfig\kern1pt) \eqdf$}
&
\mbox{$
\left(\rule{0pt}{12pt}
\bigwedge_{i = 1}^k\phi_{\srfs{Smgr}}(\bfif\Fsub{i}) \stwedge
\bigwedge_{i = 1}^k\phi_{\srfs{Smgr}}(\bfig\fsub{i}) \stwedge
(\theta_{=})^{(\phi_{\trfs{Smgr}})}(\vecbfif,\vecbfig\kern1pt)
\right) \stvee$
}
\\&
\mbox{$
\left(\rule{0pt}{12pt}
\stneg\phi_{\srfs{Smgr}}(\bfif\fsub{1}) \stwedge (\bfif\fsub{1} = \bfig\fsub{1})
\right)
$.}
\vspace{-4mm}
\end{align*}
Define
\vspace{-2mm}
\begin{align*}
\psi_{\srfs{App}\nsub{1}}(\vecbfif,\vecbfig,\vecbfih\kern1pt) \eqdf &
\mbox{$\bigwedge_{i = 1}^k$}
\left(\rule{0pt}{12pt}\phi_{\srfs{Smgr}}(\bfif\fsub{i}) \stwedge
\phi_{\srfs{Smgr}}(\bfig\fsub{i}) \stwedge \phi_{\srfs{Smgr}}(\bfih\fsub{i})
\right)
\stwedge
\\&
(\theta_{\srfs{App}})^{(\phi_{\trfs{Smgr}})}
(\vecbfif,\vecbfig,\vecbfih\kern1pt).
\vspace{-2mm}
\end{align*}
\underline{Explanation}:
Note that
\begin{align}
\tag{$\star$}
(\rfs{App}\nsub{1})^{\srfs{Act}\nsup{1}(C)} =
(\rfs{App}\nsub{1})^{\srfs{Act}(C)} =
\rfs{App}^{\srfs{Act}(\srfs{Smgr}(C))}.
\end{align}
Since
$\theta_{\srfs{App}}(\vecbfif,\vecbfig,\vecbfih\kern1pt)$
is the formula which interprets
$\rfs{App}^{\srfs{Act}(\srfs{Smgr}(C))}$
in\break
$\rfs{Alg}(\rfs{Smgr}(C))$, its relativization to
$\phi_{\srfs{Smgr}}$, applied to members of $(\rfs{Smgr}(C))^k$,
interprets $\rfs{App}^{\srfs{Act}(\srfs{Smgr}(C))}$
in $\rfs{Alg}(C)$.
So by ($\star$), $\psi_{\srfs{App}\nsub{1}}(\vecbfif,\vecbfig,\vecbfih\kern1pt)$
is an interpreting formula for $(\rfs{App}\nsub{1})^{\srfs{Act}\nsup{1}(C)}$
in $\rfs{Alg}(C)$.
\vspace{2mm}

We now define $\psi_{\srfs{Cmp}\nsub{n,k}}$.
It is trivial that
$\psi_{\srfs{Cmp}_{n,k}}(\bfif,\bfig\fsub{1},\ldots,\bfig\fsub{n},\bfih\kern1pt)$
is the formula
$\bfih = \rfs{Cmp}_{n,k}(\bfif,\bfig\fsub{1},\ldots,\bfig\fsub{n})$.
(Recall that $\bfih = \rfs{Cmp}_{n,k}(\bfif,\bfig\fsub{1},\ldots,\bfig\fsub{n})$
stands for
$\rfs{Cmp}_{n,k}(\bfif,\bfig\fsub{1},\ldots,\bfig\fsub{n},\bfih\kern1pt)$.
So, in fact,
$\psi_{\srfs{Cmp}_{n,k}}(\bfif,\bfig\fsub{1},\ldots,\bfig\fsub{n},\bfih\kern1pt)$
is the formula
$\rfs{Cmp}_{n,k}(\bfif,\bfig\fsub{1},\ldots,\bfig\fsub{n},\bfih\kern1pt)$
itself.)
\vspace{2mm}

We have defined the FO-interpretation of
$\setm{\rfs{Act}^1(C)}{C \stin K_{\srfs{FT-cln}}}$ in\break
$\setm{\rfs{Alg}(C)}{C \stin K_{\srfs{FT-cln}}}$,
and it remains to define $\psi^*_{\srfs{Imap}}(\vecbfif,\bfig\kern1pt)$
and $\psi_{\srfs{Imap}}(\vecbfif,\bfig\kern1pt)$.
For $C \stin K_{\srfs{FT-cln}}$
let $\sigma_C$ denote the interpreting mapping from\\
$\theta_{\srfs{U}}[(\rfs{Alg}(\rfs{Smgr}(C)))^k]$ to
$\abs{\rfs{Act}(\rfs{Smgr}(C))}$.
Then by the definition of the FO-interpretation
$\fourtpl{\psi_{\srfs{U}}}{\psi_{=}}{\psi_{\srfs{App}\nsub{1}}}{\ldots}$,
the interpreting mapping $\tau_C$ from\break
$\psi_{\srfs{U}}[(\rfs{Alg}(C)^k]$ to $\abs{\rfs{Act}^1(C)}$
should be defined as follows.
Let $\vecf \eqdf \fseqn{f_1}{f_k}$.
Then
\begin{equation*}
\tau_C(\vecf) = 
\begin{cases}
\sigma_C(\vecf) & \text{ if } \vecf \stin
\theta_{\srfs{U}}[(\rfs{Alg}(\rfs{Smgr}(C)))^k]
\\
f_1 & \text{ if } \vecf \stin C^k
\text{ and } f_1 \stnot\stin \rfs{Smgr}(C).
\end{cases}
\end{equation*}

We define $\psi^*_{\srfs{Imap}}$.
The formula $\psi^*_{\srfs{Imap}}$ should have the propery that
for every $C \stin K$,
\vspace{-2mm}
$$
\psi^*_{\srfs{Imap}}[(\rfs{Act}^1(C))^{k + 1}] =
\setm{\fourtpl{f_1}{\ldots}{f_k}{\tau_C(\vecf\kern3pt)}}
{\vecf \stin \rfs{Dom}(\tau_C)},
\vspace{-2mm}
$$
where $\vecf = \fseqn{f_1}{f_k}$.

Let $\eta(\bfix\kern1pt)$ be the following $\calL_{\srfs{Act-Cln}}$-formula.
\vspace{-2mm}
$$
\eta(\bfix\kern1pt) \eqdf
\exists \bfiy\kern1pt \exists \bfiz\kern1pt
\left(\rule{0pt}{12pt}
\left(\rule{0pt}{10pt}\rfs{Cmp}_{1,1}(\bfiy,\bfix\kern1pt) = \bfiz\kern1pt\right)
\stvee
\left(\rule{0pt}{10pt}\rfs{Cmp}_{1,1}(\bfix,\bfiy\kern1pt) = \bfiz\kern1pt\right)
\right).
\vspace{-2mm}
$$
Clearly, for every clone $C$,
$\eta[\rfs{Act}^1(C)] = \abs{\rfs{Act}(\rfs{Smgr}(C)}$.
That is,\break
$\eta$ defines $\abs{\rfs{Act}(\rfs{Smgr}(C)}$ in $\rfs{Act}^1(C)$.
Let $\theta^*_{\srfs{Imap}}(\vecbfif,\bfih\kern1pt)$
be the formula participating in the FS-interpretation
of $\setm{\rfs{Act}(S)}{S \stin K_{\srfs{FT}}}$
in $\setm{\rfs{Alg}(S)}{S \stin K_{\srfs{FT}}}$.
Let
$$\mbox{
$\psi^*_{\srfs{Imap-1}}(\vecbfif,\bfih\kern1pt) \eqdf
\bigwedge_{i=1}^k\eta(\bfif\fsub{i}) \stwedge \eta(\bfih\kern1pt)
\stwedge
(\theta^*_{\srfs{Imap}})^{(\eta)}(\vecbfif,\bfih\kern1pt)$.
}$$
$\psi^*_{\srfs{Imap-1}}$
defines in $\rfs{Act}^1(C)$ the restriction of $\tau_C$ to
$\theta_{\srfs{U}}[(\rfs{Act}((\rfs{Smgr}(C))^k]$.
As was explained, for members $h$ of
$\abs{\rfs{Act}^1(C)} \swsetminus \abs{\rfs{Act}(\rfs{Smgr}(C))}$,
$\vecf$ represents $h$ iff $f_1 = h$.
So the formula
$$
\psi^*_{\srfs{Imap}}(\vecbfif,\bfih\kern1pt) \eqdf
\psi^*_{\srfs{Imap-1}}(\vecbfif,\bfih\kern1pt) \stvee
\left(\rule{0pt}{12pt}\stneg\eta(\bfif\fsub{1}\kern2.3pt) \stwedge
(\bfif\fsub{1} = \bfih\kern1pt)\right)
$$
is as needed.
\vspace{2mm}

We now define $\psi_{\srfs{Imap}}$.
Let $\theta_{\srfs{Imap}}(\vecbfif,\bfih\kern1pt)$
be the formula participating in the FS-interpretation
of $\setm{\rfs{Act}(S)}{S \stin K_{\srfs{FT}}}$
in $\setm{\rfs{Alg}(S)}{S \stin K_{\srfs{FT}}}$.
Set $\vecbfif \eqdf \fseqn{\bfif\fsub{1}}{\bfif\fsub{n}}$,
and define
\begin{align*}
\phi_{\srfs{Imap}}(\vecbfif,\bfih\kern1pt) \eqdf
&
\left(\rule{0pt}{12pt}\mbox{$\bigwedge_{i = 1}^k$}
\phi_{\srfs{Smgr}}(\bfif\fsub{i}) \stwedge
\phi_{\srfs{Smgr}}(\bfih\kern1pt) \stwedge
\theta_{\srfs{Imap}}^{(\phi_{\trfs{Smgr}})}(\vecbfif,\bfih\kern1pt)\right)
\stvee
\\&
\left(\rule{0pt}{12pt}\stneg\phi_{\srfs{Smgr}}(\bfif\kern2.3pt) \stwedge
(\bfif\fsub{1} = \bfih\kern1pt)\right).
\vspace{-2mm}
\end{align*}
Obviously, $\phi_{\srfs{Imap}}$ has its required propery.
That is, for every $C \stin K$,
{\thickmuskip=2mu \medmuskip=1mu \thinmuskip=1mu
\reqnomode
\begin{align*}
\tag*{\solidqed}
\psi_{\srfs{Imap}}[(\rfs{Alg}(C))^{k + 1}] =
\setm{\fourtpl{f_1}{\ldots}{f_k}{\tau_C(\vecf\kern3pt)}}
{\vecf \stin \rfs{Dom}(\tau_C) \stcap \abs{\rfs{Alg}(C)}}.\rule{5mm}{0pt}
\end{align*}
}
%
\vspace{-2mm}

\noindent
{\bf Proof of Theorem C in the introduction }
Theorem C says that $K^{\srfs{Act}}_{\srfs{FT-cln}}$\break
is FS-interpretable in
$K^{\srfs{Alg}}_{\srfs{FT-cln}}$.
Bi-definability is surely an equivalence relation. So by Corollary~\ref{c10.3}
and Observation~\ref{o10.4},
$(\dagger)$ $\setm{\rfs{Act}^1(C)}{C \stin K_{\srfs{FT-cln}}}$ and
$\setm{\rfs{Act}(C)}{C \stin K_{\srfs{FT-cln}}}$ are bi-definably-interpretable.
By Cor\-ollary~\ref{c10.5}, $(\ddagger)$
{\thickmuskip=4mu \medmuskip=3mu \thinmuskip=2mu
$\setm{\rfs{Act}^1(C)}{C \stin K_{\srfs{FT-cln}}}$ is FS-interp\-retable in
$\setm{\rfs{Alg}(C)}{C \stin K_{\srfs{FT-cln}}}$.}
By $(\dagger)$, $(\ddagger)$ and Observation~\ref{new-o2.6},
$K^{\srfs{Act}}_{\srfs{FT-cln}}$ is FS-interpretable in
$K^{\srfs{Alg}}_{\srfs{FT-cln}}$.
\hfill\solidqed
\vspace{2mm}

\noindent
{\bf Proof of Corollary D in the introduction }
The corollary follows from Theorem C and Theorem~\ref{t2.4-08-14}.
\hfill\solidqed

\kern 1mm

\section{Index}
\noindent
{\bf Symbol index by order of appearance}

\kern 1mm

\noindent
\indexentry{$\phi[M^n] \eqdf \setm{\vec{a} \stin |M|^n\ }{\ M\models \phi[\vec{a}]}$}{7}
\indexentry{$\phi[M^k,M^{\ell}] \eqdf \setm{\pair{\veca}{\vecb} \stin |M|^k \sttimes |M|^{\ell}}{M\models\phi[\veca,\vecb]}$}{8}
\indexentry{$\iso{f}{M}{N}$ means that $f$ is an isomorphism between $M$ and $N$}{9}
\indexentry{$\phi(K) \eqdf \setm{M \stin K}{M \models \phi}$}{10}
\indexentry{$\phi^{(\alpha)}$. The relativization of $\phi$ to $\alpha$}{11}
\indexentry{$f^h \eqdf h \circ f \circ h\inverse$}{12}
\indexentry{$\sim_f$. $a \sim_f b$ if $f(a) = f(b)$}{17}
\indexentry{$[a]_f \eqdf a \kern1.5pt\slash\kern-2.5pt \sim_f$}{17}
\indexentry{$\rvpair{a}{b}^A$. The transposition of $A$ which switches $a$ and $b$.\\\indent(Abbreviated by $\rvpair{a}{b}$)}{18}
\indexentry{$\rvtrpl{a}{b}{c}$}{27}

\kern 1mm

\noindent
{\bf Notation index by alphabetic order}

\kern 1mm

\baselineskip 16.5pt
\noindent
\indexentry{$\rfs{Act}(S) \eqdf \pair{S \stcup A}{\rfs{App}^S}$. ($S$ is a function semigroup)}{4}
\indexentry{$K^{\srfs{Act}} \eqdf \setm{\rfs{Act}(S)}{S \stin K}$}{14}
\indexentry{$\rfs{Act}(C) \eqdf \pair{C \stcup A}{\kern3pt\sngltn{\rfs{App}^C_n)}_{n \stin \bbN^+}}$. ($C$ is a clone)}{6}
\indexentry{$J^{\srfs{Act}} \eqdf \setm{\pair{\rfs{Act}(S)}{f}}{\pair{S}{f} \stin J}$}{15}
\indexentry{$\rfs{Act}^+(C) \eqdf \trpl{C \stcup A}{\sngltn{\rfs{App}^C_n}_{n \stin \bbN^+}}{\sngltn{\rfs{Cmp}^C_n}_{n \stin \bbN^+}}$}{60}
\indexentry{$\rfs{Act}^1(C) \eqdf \trpl{C \stcup A}{\rfs{App}^C_1}{\sngltn{\rfs{Cmp}^C_n}_{n \stin \bbN^+}}$}{60}
\indexentry{$J^{\srfs{ActGr}} \eqdf \setm{\pair{\rfs{ActGr}(S)}{f}}{\pair{S}{f} \stin J}$}{15}
\indexentry{$K^{\srfs{ActGr}} \eqdf \setm{\rfs{Act}^{\srfs{Gr}}(S)}{S \stin K}$}{13}
\indexentry{$\rfs{Act}^{\srfs{Gr}}(S) \eqdf \trpl{S \stcup A}{\circ^S}{\rfs{App}^S_{\srfs{Gr}}}$}{12}
\indexentry{$\rfs{Alg}(S) \eqdf \pair{S}{\circ^S}$. ($S$ is a function semigroup)}{3}
\indexentry{$K^{\srfs{Alg}} \eqdf \setm{\rfs{Alg}(S)}{S \stin K}$}{14}
\indexentry{$\rfs{Alg}(C) \eqdf \pair{C}{\kern3pt\sngltn{\rfs{Cmp}^C_{n,k})}_{n,k \stin \bbN^+}}$. ($C$ is a clone)}{6}
\indexentry{$J^{\srfs{Alg}} \eqdf \setm{\pair{\rfs{Alg}(S)}{f}}{\pair{S}{f} \stin J}$}{15}
\indexentry{$\alpha_{\srfs{Not-fxd}}$. Abbreviated by $\alpha$}{28}
\indexentry{$\rfs{App}^S$. The application function of a function semigroup $S$}{4}
\indexentry{$\rfs{App}^C_n$. The $n$-place application function of a clone $C$}{6}
\indexentry{$K^{\srfs{Aug}} \eqdf \setm{\pair{S}{f}}{S \stin K \mbox{ and } f \stin S}$}{15}
\indexentry{$\bfB(f) \eqdf \setm{[a]_f}{a \stin A \mbox{ and } \abs{[a_f]} \geq 2}$}{17}
\indexentry{$\beta_{\srfs{Not-fxd}}$. Abbreviated by $\beta$}{28}
\indexentry{$\rfs{Cmp}^C_{n,k}$. The $(n,k)$'th composition function of a clone $C$}{6}
\indexentry{$\rfs{cnst}(f)$. The constant value of a semi-constant function $f$}{20}
\indexentry{$\rfs{Cnst-dom}(f)$. The constant domain of a semi-constant function $f$}{20}
\indexentry{$\rfs{Dom}(R)$. The domain of a relation $R$}{11}
\indexentry{$\bfs{Dom}(S)$. If $S \subseteq A^A$, then $\bfs{Dom}(S) \eqdf A$}{3}
\indexentry{$\bfs{Dom}(C)$. If $C$ is a clone on $A$, then $\bfs{Dom}(C)$ denotes $A$}{6}
\indexentry{$\rfs{F}\Nsub{\srfs{\sbf{B}$\geq$3}}(S) \eqdf \setm{f \stin S}{\abs{\bfB(f)} \geq 3}$}{17}
\indexentry{$\rfs{F}\Nsub{\srfs{\sbf{B}$\leq$2}}(S) \eqdf \setm{f \stin S}{\abs{\bfB(f)} \leq 2}$}{17}
\indexentry{$\rfs{Fnc}_n(C) \eqdf C \stcap A^{A^n}$}{5}
\indexentry{$\phi_{\srfs{Oo-pair}}(\bfif,\bfix,\bfiu\kern1.5pt)$}{38}
\indexentry{$\phi_{\srfs{\sbf{B}$\geq$3}}(\bfif\kern2pt)$}{26}
\indexentry{$\phi_{\srfs{Cnst}}(\bfif\kern2pt)$. An $\calL^{\srfs{ActGr}}_{\srfs{FT}}$-formula which says that $\bfif$ is a constant function}{19}
\indexentry{$\phi_{\srfs{$\exists$-cnst}} \eqdf \exists \bfif\kern2.9pt \forall \bfig\kern1pt(\bfif \bfig = \bfif\kern2.3pt)$}{16}
\indexentry{$\phi_{\srfs{$\exists$-scnst}}$. A sentence which says that $S \stin K_{\srfs{$\exists$-scnst}}$}{19}
\indexentry{$\phi_{\srfs{Fxd-img}}^{\srfs{$\exists$-scnst}}(\bfif,\bfib\kern1pt)$}{22}
\indexentry{$\phi^{J_1}_{\srfs{Fxd-img}}$}{34}
\indexentry{$\phi^{J_2}_{\srfs{Fxd-img}} \eqdf \stneg\beta_{\srfs{Not-fxd}}$}{33}
\indexentry{$\phi^{J_3}_{\srfs{Fxd-img}}(\bfif,\bfix\kern1pt)$}{55}
\indexentry{$\phi^{J_4}_{\srfs{Fxd-img}}(\bfif,\bfix\kern1pt)$}{35}
\indexentry{$\phi_{\srfs{Fxd-img}}(\bfif,\bfix\kern1.5pt)$}{16}
\indexentry{$\phi_{\srfs{Gr}}(\bfif\kern2pt) \eqdf \exists \bfig\kern1pt(\phi_{\srfs{Id}}(\bfif \bfig\kern1pt) \stwedge \phi_{\srfs{Id}}(\bfig \bfif\kern2.3pt))$}{12}
\indexentry{$\phi_{\srfs{Id}}(\bfif\kern2pt) \eqdf \forall \bfig\kern1pt(\bfif \bfig = \bfig \bfif = \bfig\kern1pt)$}{12}
\indexentry{$\phi_{\srfs{In-rng}}(\bfif,\bfix\kern1.5pt)$. For $f$'s with at least three 1-1 pairs it says that\\\indent$\bfix \stin \rfs{Rng}(\bfif\kern2.6pt)$}{42}
\indexentry{$\phi_{\srfs{In-rng}}^1(\bfif,\bfix\kern1pt)$}{39}
\indexentry{$\phi^{J_1}$}{34}
\indexentry{$\phi^{J_2}(\bfif\kern2.3pt)$}{27}
\indexentry{$\phi^{J_3}(\bfif\kern2.3pt)$}{27}
\indexentry{$\phi^{J_4}(\bfif\kern2.3pt)$}{27}
\indexentry{$\phi_{\srfs{Mo-img}}(\bfif,\bfix\kern1.5pt)$}{45}
\indexentry{$\phi_{\srfs{Mo-pre}}(\bfif,\bfix\kern1.5pt)$}{36}
\indexentry{$\phi_{\srfs{Not-in-rng}}(\bfif,\bfix\kern2pt)$. For $f$'s with at least three 1-1 pairs it says that\\\indent$\bfix \stnot\stin \rfs{Rng}(\bfif\kern2.6pt)$}{42}
\indexentry{$\phi_{\srfs{Oo-img}}(\bfif,\bfix\kern1.5pt)$}{41}
\indexentry{$\phi_{\srfs{Prj}}(\bfif\kern2pt)$. An $\calL^{\srfs{ActGr}}_{\srfs{FT}}$-formula which says that $\bfif$ is a projection}{19}
\indexentry{$\phi_{\srfs{scnst}}(\bfif\kern2pt)$. An $\calL^{\srfs{ActGr}}_{\srfs{FT}}$-formula saying that $\bfif$ is a semi-constant function}{19}
\indexentry{$\phi_{\srfs{sim}}(\bfif,\bfix,\bfiy\kern1pt)$}{19}
\indexentry{$\phi_{\srfs{smpl-pair}}(\bfif,\bfix,\bfiu\kern1.5pt)$}{38}
\indexentry{$\rfs{Fxd}(f) \eqdf \setm{a \stin \rfs{Dom}(f)}{f(a) = a}$}{16}
\indexentry{$\rfs{Fxd-img}(f) \eqdf \setm{a \stin A}{f(a) \stin \rfs{Fxd}(f)}$}{16}
\indexentry{$\rfs{Gr}(S) \eqdf \setm{g \stin S \stcap \rfs{Sym}(A)}{g\inverse \stin S}$}{12}
\indexentry{$\rfs{Idp}(f) \eqdf \setm{a \stin A}{f\inverse[\sngltn{a}] = \sngltn{a}}$}{20}
\indexentry{$\rfs{Img}_f(B)$}{30}
\indexentry{$J_1 \eqdf (K_{\srfs{$\exists$-scnst}})^{\srfs{Aug}}$}{17}
\indexentry{$J_2$}{17}
\indexentry{$J_3$}{17}
\indexentry{$J_4$}{17}
\indexentry{$K^{\srfs{Act}}_{\srfs{FT-cln}} \eqdf \setm{\rfs{Act}(C)}{C \stin K_{\srfs{FT-cln}}}$}{6}
\indexentry{$K^{\srfs{Act}}_{\srfs{FT}} \eqdf \setm{\rfs{Act}(S)}{S \stin K_{\srfs{FT}}}$}{4}
\indexentry{$K^{\srfs{Act}\nsup{1}}_{\srfs{Sm-trnsv}}$}{60}
\indexentry{$K^{\srfs{Act}}_{\srfs{Sm-trnsv}}$}{60}
\indexentry{$K^{\srfs{Alg}}_{\srfs{FT-cln}} \eqdf \setm{\rfs{Alg}(C)}{C \stin K_{\srfs{FT-cln}}}$}{6}
\indexentry{$K_{\srfs{$\exists$-cnst}}$. The class of all $S \stin K_{\srfs{FT}}$ such that $S$ contains a constant\\\indent function}{15}
\indexentry{$K_{\srfs{$\exists$-scnst}}$. The class of all $S \stin K_{\srfs{FT}} \swsetminus K_{\srfs{$\exists$-cnst}}$ such that $S$ contains a\\\indent semi-constant function}{17}
\indexentry{$K^{\srfs{Alg}}_{\srfs{FT}} \eqdf \setm{\rfs{Alg}(S)}{S \stin K_{\srfs{FT}}}$}{3}
\indexentry{$K_{\srfs{FT-cln}}$. The class of all FT clones $C$ such that $\abs{\bfs{Dom}(C)} \neq 1,2,6$}{6}
\indexentry{$K_{\srfs{FT}}$. The class of all FT semigroups $S$ such that $\abs{\bfs{Dom}(S)} \neq 1,2,6$}{3}
\indexentry{$K_{\srfs{No-cnst}} \eqdf K_{\srfs{FT}} \swsetminus K_{\srfs{$\exists$-cnst}}$}{16}
\indexentry{$K_{\srfs{No-scnst}} \eqdf K_{\srfs{FT}} \swsetminus (K_{\srfs{$\exists$-scnst}} \stcup K_{\srfs{$\exists$-cnst}})$}{17}
\indexentry{$K_{\srfs{Sm-trnsv}}$}{60}
\indexentry{$\calL(K)$. The language of the class of structures $K$}{8}
\indexentry{$\calL(M)$. The language of the structure $M$}{7}
\indexentry{$\calL_{\srfs{Alg-cln}}$. The language of $\rfs{Alg}(C)$}{60}
\indexentry{$\calL_{\srfs{ActGr}}$}{13}
\indexentry{$\rfs{Mo-img}(f) \eqdf \setm{x \stin A}{\abs{f\inverse[\sngltn{x}]} \geq 2}$}{35}
\indexentry{$\rfs{Mo-pre}(f) \eqdf \bigcup \bfB(f)$}{35}
\indexentry{$\rfs{Oo-img}(f) \eqdf \setm{f(a)}{a \stin \rfs{Oo-pre}(f)}$}{35}
\indexentry{$\rfs{Oo-pre}(f) \eqdf \setm{a \stin A}{[a]_f = \sngltn{a}}$}{35}
\indexentry{$\rfs{Rng}(R)$. The the range of a relation $R$}{11}
\indexentry{$\rfs{Sym}(A)$. The symmetric group of $A$}{11}

\noindent
{\bf Index of definitions by alphabetic order}
\kern 1mm

\noindent
\indexentry{$f$-1-1 image. If $\pair{p}{q}$ is an $f$-1-1 pair, then $q$ is called an $f$-1-1 image}{37}
\indexentry{$f$-1-1 pair. $\pair{p}{q}$ is an $f$-1-1 pair if $f\inverse[\sngltn{q}] = \sngltn{p}$}{37}
\indexentry{$f$-1-1 preimage. If $\pair{p}{q}$ is an $f$-1-1 pair, then $p$ is called an $f$-1-1\\\indent preimage}{37}
\indexentry{augmented function semigroup}{15}
\indexentry{Bi-definably interpretable. $K$ and $K^*$ are bi-definably-interpretable\\\indent relative to $\dbltn{R}{R\inverse}$}{9}
\indexentry{Definably interpretable. $K^*$ is definably interpretable in $K$ relative to $R$}{9}
\indexentry{defining formula}{10}
\indexentry{defining sentence}{10}
\indexentry{First order interpretable}{7}
\indexentry{First order strongly interpretable}{8}
\indexentry{FO-definable subclass (first order definabe subclass)}{10}
\indexentry{FO-interpretable. Abbreviation of ``first order interpretable''}{7}
\indexentry{FO-interpretation}{8}
\indexentry{FS-interpretable}{8}
\indexentry{FS-interpretation}{9}
\indexentry{FT clone. A clone on a set $A$ which contains all transpositions of $A$}{6}
\indexentry{FT semigroup. A function semigroup on a set $A$ which contains all\\\indent transposions of $A$}{3}
\indexentry{Fully-transpositional semigroup. Abbreviated by FT semigroup}{3}
\indexentry{Function semigroup}{3}
\indexentry{Interpretating mapping}{8}
\indexentry{Relativization of a formula to another formula}{11}
\indexentry{Same-universe system}{9}
\indexentry{Semi-transitive clone}{59}
\indexentry{Semi-transitive semigroup}{59}
\indexentry{Semi-constant function}{17}
\indexentry{$f$-simple image. The $f$-1-1 image of a simple pair}{37}
\indexentry{$f$-simple pair. An $f$-1-1 pair $\pair{p}{q}$ such that $p \neq q$}{37}
\indexentry{$f$-simple preimage. The $f$-1-1 preimage of a simple pair}{37}
\indexentry{Subuniverse system}{8}
\indexentry{Transitive semigroup}{59}

\kern 1mm


%
%
%
%
%
\end{document}